\documentstyle{amsppt}
\def\deg{\hbox{\rm deg}}
\def\volum{\hbox{\rm vol}}
\def\mult{\hbox{\rm mult}}
\def\div{\hbox{\rm div}}
\def\Sing{\hbox{\rm Sing}}
\def\Ric{\hbox{\rm Ric}}
\def\Jac{\hbox{\rm Jac}}
\def\Pic{\hbox{\rm Pic}}
\def\Tr{\hbox{\rm Tr}}
\def\Td{\hbox{\rm Td}}
\def\e{\hbox{\rm e}}
\def\Im{\hbox{\rm Im}}
\def\Re{\hbox{\rm Re}}
\def\End{\hbox{\rm End}}
\def\Res{\hbox{\rm Res}}
\def\I{\hbox{\rm I}}
\def\II{\hbox{\rm II}}
\def\codim{\hbox{\rm codim}}
\def\Am{\hbox{\rm Am}}
\def\NE{\hbox{\rm NE}}
\def\Exc{\hbox{\rm Exc}}
\def\Sym{\hbox{\rm Sym}}
\NoRunningHeads
\overfullrule=0pt
\magnification=\magstep1
\document
\topmatter
\title
$\pmb{K3}$ Surfaces with Involution and Analytic Torsion
\endtitle
\author
Ken-ichi Yoshikawa
\endauthor
\address
Graduate School of Mathematics, Nagoya University,
Nagoya 464-8602, Japan
\endaddress
\email
yosikawa\@math.nagoya-u.ac.jp
\endemail
\endtopmatter

\beginsection 
$\S1$. Introduction

\par
In a series of works [Bo3-5], Borcherds developed a theory of modular 
forms over domains of type IV which admits an infinite product 
expansion. Such modular forms are said to be Borcherds's product in this 
paper. Among all Borcherds's products, Borcherds's $\Phi$-function 
([Bo4]) has an interesting geometric background; It is a modular form on 
the moduli space of Enriques surfaces characterizing the discriminant 
locus. In his construction, $\Phi$-function is obtained as the 
denominator function of one of the fake monster Lie superalgebras 
([Bo2, $\S14$]), although Enriques surface itself plays no role. 
After Borcherds, Jorgenson-Todorov ([J-T2,3]) and Harvey-Moore ([H-M])
discovered that {\it the Ray-Singer analytic torsion }\rm ([R-S]) of an
Enriques surface equipped with the normalized Ricci-flat K\"ahler metric
coincides with Borcherds's $\Phi$-function at its period point. The goal 
of this paper is to give a rigorous proof to their observation and
generalize it to an interesting class of $K3$ surfaces studied by Nikulin
([Ni4]). Let us briefly recall these surfaces.
\par
Let $(X,\iota)$ be a $K3$ surface with an anti-symplectic involution.
Let $M$ be a $2$-elementary hyperbolic lattice. The pair $(X,\iota)$ is
said to be a $2$-elementary $K3$ surface of type $M$ if the invariant 
sublattice of $H^{2}(X,\Bbb Z)$ with respect to the action of $\iota$ is 
isometric to $M$. Since $X/\iota$ is an Enriques surface when
$M\cong\II_{1,9}(2)$, $K3$ surfaces of this class are a kind of 
generalizations of Enriques surfaces. (We denote by $\I_{p,q}$ 
(resp. $\II_{p,q}$) the odd (resp. even) unimodular lattice of signature 
$(p,q)$.) By the Torelli theorem ([P-S-S]) and surjectivity of the period 
map ([To1]), the moduli space of $2$-elementary $K3$ surfaces of type $M$ 
is isomorphic to an arithmetic quotient of an open subset of the symmetric 
bounded domain of type IV via the period map. The period domain is denoted 
by $\Omega_{M}$, the modular group by $\Gamma_{M}$, and the moduli space 
by $\Cal M_{M}^{0}=(\Omega_{M}\backslash\Cal D_{M})/\Gamma_{M}$. Here, 
$\Cal D_{M}$ is the discriminant locus. Let $X^{\iota}$ be the fixed locus 
of $\iota$. Since only one of the irreducible components can have positive 
genus, $g(M)(\geq0)$, unless $M\cong\II_{1,9}(2)$ or 
$\II_{1,1}\oplus E_{8}(-2)$, one can define a map $j_{M}^{0}:
\Cal M_{M}^{0}\to\Cal A_{g(M)}$ to the Siegel modular variety by taking 
the period point of the Jacobi variety of $X^{\iota}$. 
(If $M\cong\II_{1,1}\oplus E_{8}(-2)$, $j_{M}$ takes its value in 
$S^{2}(\Cal A_{1})$, the second symmetric product of $\Cal A_{1}$.) 
Let $\Cal F_{g(M)}$ be the sheaf of Siegel modular forms of weight $1$ 
on $\Cal A_{g(M)}$. We define a $\Gamma_{M}$-invariant sheaf on (an 
open subset of) $\Omega_{M}$ by $\lambda_{M}:=j_{M}^{*}\Cal F_{g(M)}$.
A $\lambda_{M}^{\otimes q}$-valued modular form of weight $p$ 
(with a character of $\Gamma_{M}$) is said to be an automorphic form 
of weight $(p,q)$ in this paper whose Petersson norm is denoted by 
$\|\cdot\|$. 
\par
In this paper, we focus on the Ray-Singer analytic torsion regarded as 
a function on the moduli space of $K3$ surfaces. Since it coincides with 
(the Petersson norm of) the Jacobi $\Delta$-function for elliptic curves 
([R-S]), it is natural to expect that it may yield an interesting modular 
form for $K3$ surfaces ([J-T1]). With small modifications, it is the case 
at least for $2$-elementary $K3$ surfaces. We introduce the following 
function (see $\S3$) on the {\it polarized period domain }\rm of 
$2$-elementary $K3$ surfaces of type $M$;
$$
\tau_{M}(X,\iota,\kappa):=\volum(X,\kappa)^{\frac{14-r(M)}{8}}
\tau(X/\iota,\kappa)
\sqrt{\volum(X^{\iota},\kappa|_{X^{\iota}})\,
\tau(X^{\iota},\kappa|_{X^{\iota}})}.
\tag 1.1
$$
Here $\kappa$ is an $\iota$-invariant Ricci-flat K\"ahler metric of $X$,
$r(M)$ is the rank of $M$, $\tau(X/\iota,\kappa)$ is the Ray-Singer 
analytic torsion of orbifold $(X/\iota,\kappa)$ (or equivalently, the 
equivariant Ray-Singer analytic torsion of $(X,\iota,\kappa)$) and 
$\tau(X^{\iota},\kappa|_{X^{\iota}})$ is that of the fixed curve 
$(X^{\iota},\kappa|_{X^{\iota}})$. Note that
$\volum(X^{\iota},\kappa|_{X^{\iota}})$ and
$\tau(X^{\iota},\kappa|_{X^{\iota}})$ are defined multiplicatively with
respect to the irreducible decomposition of divisors.

\proclaim{Theorem 1.1 (Theorems 3.4, 7.2, 7.3)}
Suppose that $r(M)\leq17$.
\flushpar{$(1)$} 
$\tau_{M}$ is independent of a choice of Ricci-flat K\"ahler metrics and 
defines an invariant of a $2$-elementary $K3$ surface of type $M$.
It descents to a smooth $\Gamma_{M}$-invariant function on 
$\Omega_{M}^{0}:=\Omega_{M}\backslash\Cal D_{M}$. 
\flushpar{$(2)$} 
There exists an automorphic form, $\varDelta_{M}$, for $\Gamma_{M}$ of 
weight $(r(M)-6,4)$ such that $\tau_{M}=\|\varDelta_{M}\|^{-1/4}$ and
$\div(\varDelta_{M})=\Cal D_{M}$.
\flushpar{$(3)$} 
If $\delta$ is a root of $M^{\perp}$ and $\langle M\oplus\delta\rangle$ is 
the smallest $2$-elementary lattice containing $M$ and $\delta$, then
$\varDelta_{\langle M\oplus\delta\rangle}$ coincides with the regularized 
restriction of $\varDelta_{M}$ to $\Omega_{\langle M\oplus\delta\rangle}$.
\endproclaim

It may be worth mentioning that analogous statment is valid for theta
divisors of dimension $\geq1$ ([Yo2]).
By (3), any $\varDelta_{M}$ is the regularized restriction of such 
$\varDelta_{S}$ that $S$ is root free. Therefore, it is important to know 
an explicit formula for them. By Nikulin's classification, there are 
three such lattices up to isometry; 
$S\cong A_{1},\,\II_{1,1}(2),\,\II_{1,9}(2)$. 

\proclaim{Theorem 1.2 (Theorems 10.1, 10.2, 9.3)}
\flushpar{$(1)$} 
$\varDelta_{A_{1}}^{15}=\Delta_{6}^{8}/\Psi_{A_{1}^{\perp}}$ where 
$\Delta_{6}$ is the discriminant of plane sextic curves, and
$\Psi_{A_{1}^{\perp}}$ is Borcherds's product attached to the 
transcendental lattice 
and $\Theta_{E_{7}}(\tau)/\Delta(\tau)$. Here, $\Theta_{E_{7}}(\tau)$ is 
the vector valued theta series of $E_{7}$-lattice and $\Delta(\tau)$ is 
the Jacobi $\Delta$-function.
\flushpar{$(2)$}
$\varDelta_{\roman{II}_{1,1}(2)}^{17}=\Delta_{4,4}^{9}/
\Psi_{\roman{II}_{1,1}(2)^{\perp}}$ where $\Delta_{4,4}$ is the
discriminant of the linear system of bidegree $(4,4)$ on
$\Bbb P^{1}\times\Bbb P^{1}$ and $\Psi_{\roman{II}_{1,1}(2)^{\perp}}$ is 
Borcherds's product attached to the transcendental lattice 
and $\eta(\tau)^{-8}\eta(2\tau)^{-8}$. 
Here, $\eta(\tau)$ is the Dedekind $\eta$-function.
\flushpar{$(3)$}
$\varDelta_{\roman{II}_{1,9}(2)}$ coincides with the Borcherds 
$\Phi$-function associated to the fake monster Lie superalgebra of rank 
$10$ up to a constant, and $\varDelta_{\roman{II}_{1,9}(2)}^2$ is 
Borcherds's product attached to $\I_{2,10}(2)$ and
$\eta(\tau)^{-8}\eta(2\tau)^{8}\eta(4\tau)^{-8}$. 
\endproclaim

For the precise meaning of Borcherds's products used in this paper, see
Theorem 8.1 and Definition 8.1. Concerning (3), another product formula
can be found in [Bo5]. Jorgenson-Kramer ([J-K]) treat Borcherds's 
$\Phi$-function from the view point of Green currents.
Since Lefschetz fixed point formula in Alakelov theory is established by 
K\"ohler-Roessler recently ([K-R]), it will be possible to represent
$\varDelta_{M}$ by Green currents. 
\par
Among all the 2-elementary lattices, there are some interesting ones from 
the view point of Nikulin's classification. As in the case of 
$\II_{1,9}(2)$, $\varDelta_{M}$ for them is closely related to 
Borcherds's theory ([Bo1-5]).

\proclaim{Theorem 1.3 (Corollary 6.1, Theorems 8.2, 8.3, 9.4, 9.6)}
\flushpar{$(1)$} 
If $M\cong\II_{1,1}\oplus E_{8}(-2)$, one has $\varDelta_{M}^{3}=
C_{M}\,\Psi_{M^{\perp}}\otimes j_{M}^{*}(\Delta_{1}\Delta_{2})$
where $\Delta_{i}$ is the Jacobi $\Delta$-function in the $i$-th variable
on $\Cal A_{1}\times\Cal A_{1}$, $C_{M}$ is a constant, 
and $\Psi_{M^{\perp}}$ is the denominator function of a generalized 
Kac-Moody superalgebra. Moreover, $\Psi_{M^{\perp}}$ is Borcherds's 
product attached to the transcendental lattice and
$\Theta_{\Lambda_{16}}(\tau)/\Delta(\tau)$.  Here,
$\Theta_{\Lambda_{16}}(\tau)$ is the vector valued theta series of the
Barnes-Wall lattice. 
\flushpar{$(2)$} 
If $g(M)=0$, namely $M^{\perp}\cong\I_{2,20-r(M)}(2)$ 
($11\leq r(M)\leq 17$), $\varDelta_{M}$ is the denominator function of a 
generalized Kac-Moody superalgebra up to a constant. Moreover, 
$\varDelta_{M}^2$ is Borcherds's product attached to the transcendental 
lattice and $\eta(\tau)^{-8}\eta(2\tau)^{8}\eta(4\tau)^{-8}
\theta_{A_1}(\tau)^{r(M)-10}$ where $\theta_{A_1}(\tau)$ is the theta 
series of $A_1$-lattice.
\endproclaim

In view of Theorems 1.2, 1.3, it seems natural to conjecture that any 
$\varDelta_{M}$ in Theorem 1.1 is represented by Borcherds's product 
and discriminant of curves. Compared to fake monster Lie algebras, 
generalized Kac-Moody superalgebras appearing in Theorem 1.3 are not 
well understood. We shall show in Theorem 11.1 and Corollary 11.1 that
product of certain $10$ theta functions on the domain of type $I_{2,2}$
([Ma]) admits Borcherds's product which is analogous to 
Gritsenko-Nikulin's product formula for Igusa's modular form ([G-N1]). 
We remark that Gritsenko-Nikulin's $\Delta_{11}$ ([G-N3]) coincides 
with $\varDelta_{M}$ in (2) if $r(M)=17$. It is a mystery why 
the Ray-Singer analytic torsion is related to Borcherds's theory. 
(See [H-M] for a physical  explanation.) In view of works by 
Bershadski-Cecotti-Ooguri-Vafa ([B-C-O-V]) and Harvey-Moore ([H-M]), 
it may be more natural to study the Ray-Singer analytic torsion of 
Calabi-Yau $3$-folds associated to $2$-elementary $K3$ surfaces. 
(See [Bc], [V] for these Calabi-Yau $3$-folds.)
\par
Let us briefly describe the contents of this paper. In $\S2$, we recall
the theory of $2$-elementary $K3$ surfaces due to Nikulin ([Ni4]). 
In $\S3$, we recall the theory of Quillen metrics due to Bismut, Gillet,
Soul\'e ([B-G-S], [Bi]) and apply it to derive the variational formula 
of $\tau_{M}$ on $\Omega_{M}^{0}$. In $\S4$, to study the degeneration 
of Ricci-flat K\"ahler metrics along the discriminant locus, we recall 
apriori estimates of the Monge-Amper\`e equations due to Yau ([Ya]). 
Following Kobayashi's work ([Ko]), we apply them to construct an 
approximate Ricci-flat K\"ahler metrics for degenerations of type $I$.
Using it, we shall show in $\S5$ that Kronheimer's ALE instanton can be
used to compute the singularity of $\tau_{M}$. In $\S6$, we study the 
heat kernel of the ALE instanton and determine the singularity of 
$\tau_{M}$ along $\Cal D_{M}$. Theorem 1.1 shall be proved in $\S7$. 
We recall Borcherds's theory ([Bo1-5]) in $\S8,9$ and apply it to 
determine the structure of $\varDelta_{M}$ for some $M$. We compute 
some examples in $\S10,11$. 

\subsubhead
\bf{Acknowledgements}\rm
\endsubsubhead
The author is grateful to Professor R.E. Borcherds who kindly showed him 
how to construct $\Psi_{M^{\perp}}$ in $\S6.2$ from his $\Phi$-function 
of rank $26$, and modular forms of type $\rho_{M}$ from those of higher 
level, to Professor S. Kond\B o from whom he learned the moduli space of 
lattice polarized $K3$ surfaces, to Professor R. Kobayashi for teaching 
him degenerations of Ricci-flat K\"ahler metrics, and to 
Professor V.V. Nikulin and Doctor H. Imura for teaching him generalized 
Kac-Moody algebras. Without their helps, the present version would not 
exist. Together with above mentioned people, he is also grateful to 
Professors A. Fujiki, A. Matsuo, T. Ohsawa, K. Saito for discussions and
suggestions.

\beginsection 
$\S2$. $\pmb{2}$-Elementary $\pmb{K3}$ Surfaces and its Moduli Space

\subsubhead
\bf{2.1 Moduli Space of $\pmb{K3}$ Surfaces}\rm
\endsubsubhead
A compact complex surface $X$ is a $K3$ surface if and only if $q(X)=0$ 
and $K_{X}\cong\Cal O_{X}$. The cohomology lattice $H^{2}(X,\Bbb Z)$ 
together with the intersection product is isometric to 
$L_{K3}=\II_{1,1}^{3}\oplus E_{8}(-1)^{2}$ where $E_{8}$ is the positive 
definite $E_{8}$-lattice. (For a lattice $M$, $M(k)$ denotes the lattice 
with $\langle\,,\,\rangle_{M(k)}=k\langle\,,\,\rangle_{M}$,
$\Delta(M)=\{m\in M;\,\langle m,m\rangle_{M}=-2\}$ the root of $M$, 
and $M_{\Bbb R}$ (resp. $M_{\Bbb C}$) the vector space $M\otimes\Bbb R$
(resp. $M\otimes\Bbb C$). A sublattice $M'\subset M$ is said to be 
primitive if $M/M'$ is torsion free.) We denote by 
$\Pic_{X}:=H^{1,1}(X,\Bbb R)\cap H^{2}(X,\Bbb Z)$ the Picard lattice and 
by $\Delta(X)$ its root. By the Riemann-Roch theorem, 
$\Delta(X)=\Delta^{+}(X)\bigsqcup-\Delta^{+}(X)$ where $\Delta^{+}(X)$ 
consists of effective classes. The K\"ahler cone of $X$ is denoted by 
$C_{X}^{+}:=\{\kappa\in H^{1,1}(X,\Bbb R);\,
\langle\kappa,\delta\rangle>0,\,\forall\delta\in\Delta^{+}(X)\}$ and its 
closure in $H^{1,1}(X,\Bbb R)$ by $C_{X}$. 
Let $\phi_{X}:H^{2}(X,\Bbb Z)\to L_{K3}$ be an isometry. A pair 
$(X,\phi_{X})$ is said to be a marked $K3$ surface. Two marked $K3$ surfaces 
$(X,\phi_{X})$ and $(X',\phi_{X'})$ are isomorphic when there exists an 
isomorphism $f:X\to X'$ such that $\phi_{X}\circ f^{*}=\phi_{X'}$. There 
exist the fine moduli space of marked $K3$ surfaces. We briefly recall it
([P-S-S], [B-R], [B-P-V-V], [Be], [Mo], [Ni2,3] etc.). Define $\Omega$ by
$$
\Omega:=\{[\eta]\in\Bbb P(L_{K3,\Bbb C});\,
\langle\eta,\eta\rangle_{L_{K3}}=0,\,
\langle\eta,\bar{\eta}\rangle_{L_{K3}}>0\}.
\tag 2.1
$$
For $\eta\in\Omega$, let $C(\eta)=
\{x\in L_{K3,\Bbb R};\,\langle x,\eta\rangle_{L_{K3}}=0,\,
\langle x,x\rangle_{L_{K3}}>0\}$ be the positive cone over $\eta$ which
consists of two connected components; $C(\eta)=C(\eta)^{+}\bigsqcup
C(\eta)^{-}$. Let $S(\eta):=L_{K3}\cap\eta^{\perp}$ be the Picard lattice 
over $\eta$, $\Delta(\eta)$ the root of $S(\eta)$, and 
$h_{\delta}:=\{x\in L_{K3,\Bbb R};\,\langle x,\delta\rangle_{L_{K3}}=0\}$ 
the hyperplane orthogonal to $\delta$. Let 
$C(\eta)\backslash\bigcup_{\delta\in\Delta(\eta)}h_{\delta}=
\bigsqcup_{P}C_{P}$ be the decomposition into the connected components.
$P$ is parametrized by the pair $(C(\eta)^{+},\Delta_{P}^{+}(\eta))$ 
where $\Delta_{P}^{+}(\eta)$ is a partition 
$\Delta(\eta)=\Delta^{+}_{P}(\eta)\bigsqcup-\Delta^{+}_{P}(\eta)$ 
with the following property:
\newline
{\bf(P1)\rm} If $\delta_{1},\cdots,\delta_{k}\in\Delta^{+}_{P}(\eta)$ and 
$\delta=\sum_{i}n_{i}\delta_{i}$ ($n_{i}\geq 0$), then 
$\delta\in\Delta_{P}^{+}(\eta)$.
\newline
Set 
$K\Omega:=
\{(\eta,\kappa)\in\Omega\times L_{K3,\Bbb R};\,\kappa\in C(\eta)\}$
and
$$
K\Omega^{0}:=\{(\eta,\kappa)\in K\Omega;\,\kappa\in C(\eta)\backslash
\cup_{\delta\in\Delta(\eta)}h_{\delta}\},\quad
\tilde{\Omega}:=K\Omega^{0}/\sim.
\tag 2.2
$$ 
Here $(\eta,\kappa)\sim(\eta',\kappa')$ when $\eta=\eta'$ and 
$\kappa,\kappa'\in C_{P}(\eta)$ for some $P$. We denote by 
$p:K\Omega^{0}\to\tilde{\Omega}$ the projection map.
Let $\pi:\tilde{\Omega}\to\Omega$ be the projection to the first factor.
It is immediate by above construction that (1) $\tilde{\Omega}$ is
a non-separated smooth analytic space, (2) $\pi$ is \'etale and surjective,
(3) $\pi^{-1}(\eta)=\{P\}$; the set of all pairs 
$(C(\eta)^{+},\Delta(\eta)^{+})$ where $\Delta(\eta)^{+}$ satisfies 
{\bf(P1)\rm}.
\par
For $(\eta,\kappa)\in K\Omega^{0}$, put $C(\eta)^{+}_{\kappa}$ for the 
connected component containing $\kappa$, and $\Delta(\eta)^{+}_{\kappa}:=
\{\delta\in\Delta(\eta);\,\langle\kappa,\delta\rangle>0\}$. Then, 
$[(\eta,\kappa)]\in\tilde{\Omega}$ corresponds to 
$(\eta,C(\eta)^{+}_{\kappa},\Delta(\eta)^{+}_{\kappa})$.
\par
Let $\eta_{X}$ be a symplectic form; $H^{0}(X,K_{X})=\Bbb C\,\eta_{X}$.
Then, $\pi(X,\phi_{X}):=[\phi_{X}(\eta_{X})]\in\Bbb P(L_{K3,\Bbb C})$ is 
the period of $(X,\phi_{X})$ which lies on $\Omega$. Let $\kappa_{X}$ be 
a K\"ahler metric of $X$ and identify it with its class in 
$H^{2}(X,\Bbb R)$. The pair $(X,\kappa_{X})$ is said to be a polarized 
$K3$ surface and the triplet $(X,\phi_{X},\kappa_{X})$ a marked 
polarized $K3$ surface. The polarized period point of 
$(X,\phi_{X},\kappa_{X})$ is defined by
$\pi(X,\phi_{X},\kappa_{X}):=([\phi_{X}(\eta_{X})],\phi_{X}(\kappa_{X}))
\in K\Omega^{0}$ and the Burns-Rapoport period point by 
$[\pi(X,\phi_{X},\kappa_{X})]\in\tilde{\Omega}$. 
\proclaim{Theorem 2.1}
There exists the universal family of marked $K3$ surfaces over 
$\tilde{\Omega}$; $\Cal P:\Cal X\to\tilde{\Omega}$, such that the period 
map coincides with $\pi:\tilde{\Omega}\to\Omega$. 
\endproclaim

Let $M$ be a primitive hyperbolic sublattice of $L_{K3}$ with signature 
$(1,k)$. A marked $K3$ surface $(X,\phi_{X})$ is marked (ample)
$M$-polarized if $\phi_{X}(\Pic_{X})\supset M$ (and 
$\phi_{X}(C_{X}^{+})\cap M_{\Bbb R}\not=\emptyset$). As $M$ is hyperbolic, 
a marked $M$-polarized $K3$ surface is projective. Let 
$N:=M^{\perp}=\{l\in L_{K3};\,\langle l,M\rangle=0\}$ be the 
orthogonal compliment of $M$ whose signature is $(2,19-k)$. 
By definition, the period point of a marked $M$-polarized $K3$ surface is 
contained in the following subset of $\Omega$;
$$
\Omega_{M}:=\{[\eta]\in\Bbb P(N_{\Bbb C});\,\langle\eta,\eta\rangle_{N}=0,
\,\langle\eta,\bar{\eta}\rangle_{N}>0\}.
\tag 2.3
$$
Put $C(M):=\{x\in M_{\Bbb R};\,\langle x,x\rangle_{M}>0\}$ or equivalently
$C(M)=C(\eta)\cap M_{\Bbb R}$. Define
$$
K\Omega_{M}:=\{([\eta],\kappa)\in K\Omega;\,[\eta]\in\Omega_{M},\,
\kappa\in C(M)\},\quad K\Omega_{M}^{0}:=K\Omega^{0}\cap K\Omega_{M}.
\tag 2.4
$$
Set $\tilde{\Omega}^{0}_{M}:=K\Omega^{0}_{M}/\sim$ and 
$\tilde{\Omega}_{M}$ for the closure of $\tilde{\Omega}_{M}^{0}$ in 
$\tilde{\Omega}$. 

\proclaim{Lemma 2.1}
For any $\eta\in\Omega_{M}$, $M\cap C(M)\backslash
\bigcup_{\delta\in\Delta(\eta)\backslash\Delta(N)}h_{\delta}
\not=\emptyset$. 
\endproclaim

\demo{Proof}
If $M_{\Bbb Q}\cap C(M)\backslash
\bigcup_{\delta\in\Delta(\eta)\backslash\Delta(N)}h_{\delta}=\emptyset$,
then $C(M)\backslash
\bigcup_{\delta\in\Delta(\eta)\backslash\Delta(N)}h_{\delta}=\emptyset$ 
and thus there exists
$\delta\in\Delta(\eta)\backslash\Delta(N)$ such that
$M_{\Bbb R}\subset h_{\delta}$. This implies that $\delta\in M^{\perp}=N$
and contradicts the choice of $\delta$.\qed
\enddemo

Let $H_{l}:=\{[\eta]\in\Omega;\,\langle\eta,l\rangle_{L_{K3}}=0\}$ be the
hyperplane defined by $l\in L_{K3}$. We define the discriminant locus and 
an open subset by
$$
\Cal D_{M}:=\cup_{\delta\in\Delta(N)}H_{\delta},\quad
\Omega_{M}^{0}:=
\Omega_{M}\backslash\Cal D_{M}.
\tag 2.5
$$
Let $\pi_{M}:\tilde{\Omega}^{0}_{M}\to\Omega_{M}^{0}$ be the period map
restricted to $\tilde{\Omega}_{M}^{0}$. By Theorem 2.1 and Lemma 2.1,
we get the following. (See [Ni2,3] for the detail.)

\proclaim{Theorem 2.2}
$\tilde{\Omega}^{0}_{M}$ is the fine moduli space of marked ample 
$M$-polarized $K3$ surfaces whose period map is $\pi_{M}$. The universal 
family $\Cal P_{M}:\Cal X_{M}^{0}\to\tilde{\Omega}_{M}^{0}$ is obtained by 
putting $\Cal X_{M}^{0}:=\Cal X|_{\tilde{\Omega}_{M}^{0}}$ and 
$\Cal P_{M}:=\Cal P|_{\Cal X_{M}^{0}}$. 
\endproclaim

We are now interested in the boundary of $\tilde{\Omega}^{0}_{M}$ in 
$\tilde{\Omega}_{M}$ and the family of marked $M$-polarized $K3$ surfaces 
over it. For $\delta\in\Delta(N)$, put 
$H_{\delta}^{0}:=H_{\delta}\backslash
\bigcup_{d\in\Delta(N)\backslash\{\pm\delta\}}H_{d}$.

\proclaim{Lemma 2.2}
Let $\delta\in\Delta(N)$, $\eta\in H_{\delta}^{0}$ and
$\kappa\in M\cap C(M)\backslash
\bigcup_{d\in\Delta(\eta)\backslash\Delta(N)}h_{d}$. 
Then, $\Delta(\eta)=\Delta_{\kappa}^{+}(\eta)\bigsqcup
-\Delta_{\kappa}^{+}(\eta)\bigsqcup\{\pm\delta\}$ where 
$\Delta_{\kappa}^{+}(\eta)=\{d\in\Delta(\eta);\,\langle\kappa,d\rangle>0\}$.
We put 
$\Delta_{1}^{+}(\eta):=\Delta_{\kappa}^{+}(\eta)\bigsqcup\{\delta\}$ and 
$\Delta_{2}^{+}(\eta):=\Delta_{\kappa}^{+}(\eta)\bigsqcup\{-\delta\}$. 
Then, $\Delta_{i}^{+}(\eta)$ satisfies {\bf (P)}\rm.
\endproclaim

\demo{Proof}
Suppose that $d\in\Delta(\eta)\backslash\Delta(N)$. Then, either
$\langle\kappa,d\rangle>0$ or $\langle\kappa,d\rangle<0$ because 
$\kappa\not\in
\bigcup_{\delta'\in\Delta(\eta)\backslash\Delta(N)}h_{\delta'}$. 
This implies
$\Delta(\eta)\backslash\Delta(N)=\Delta_{\kappa}^{+}(\eta)\bigsqcup
-\Delta_{\kappa}^{+}(\eta)$. Suppose that $d\in\Delta(\eta)\cap\Delta(N)$
which means $d\in\Delta(N)$ and $\eta\in H_{d}$. By the choice of $\eta$,
$d=\pm\delta$ which prove the first assertion. 
Since the Weyl group $W(\eta)$ acts properly discontinuously on $C(\eta)$,
there exists a small neighborhood $K$ of $\kappa$ in $C(\eta)$ such that
$s_{d}(K)\cap K\not=\emptyset$ for $d\in\Delta(\eta)$ means $d=\pm\delta$
where $s_{d}$ is the reflection in $d$.
Thus there exists $\epsilon_{0}>0$ such that for any 
$0<\epsilon<\epsilon_{0}$, 
$\kappa+\epsilon\delta\in V(\eta)\backslash
\bigcup_{d\in\Delta(\eta)}h_{d}$. Clearly, 
$\Delta_{\kappa}^{+}(\eta)\bigsqcup\{-\delta\}=
\{d\in\Delta(\eta);\,\langle\kappa+\epsilon\delta,d\rangle>0\}$. Similarly, 
$\Delta_{1}^{+}(\eta)$ is the partition associated to
$\kappa-\epsilon\delta$.\qed
\enddemo

Let $D$ be the unit disc and $\gamma:D\to\Omega_{M}$ be a holomorphic curve 
which intersects transversally $H_{\delta}^{0}$ at $\eta_{0}=\gamma(0)$. 
For $\eta_{0}$, let $\kappa$ and
$P_{i}=(C(\eta_{0})^{+}_{\kappa},\Delta_{i}^{+}(\eta_{0}))$ be the same 
as in Lemma 2.2. Set $\kappa_{1}:=\kappa-\epsilon\delta$ and 
$\kappa_{2}:=\kappa+\epsilon\delta$. 
Let $\tau_{i}$ ($i=1,2$) be the points in $\tilde{\Omega}_{M}$ such that 
$\tau_{i}=(\eta_{0},C(\eta_{0})^{+},\Delta_{i}^{+}(\eta_{0}))$. Let $U$ be 
a small neighborhood of $\eta_{0}$ in $\Omega_{M}$, and $U_{i}$ a 
neighborhood of $\tau_{i}$ in $\tilde{\Omega}_{M}$ such that 
$\pi:U_{i}\to U$ is an isomorphism. Let
$\gamma_{i}:=(\pi|_{U_{i}})^{-1}\circ\gamma:D\to U_{i}$ be the lift of 
$\gamma$ to $U_{i}$ such that $\gamma_{i}(0)=\tau_{i}$. Let 
$(X^{i}_{t},\phi^{i}_{t})$ be the marked $K3$ surface corresponding to 
$\gamma_{i}(t)$ and $\eta_{t}=\pi(X^{i}_{t},\phi^{i}_{t})$ its period.
From construction and Burns-Rapoport's lemma ([B-R], [Mo, pp.306]),
it follows that $(\phi^{i}_{t})^{-1}(\kappa)$ is an ample class of 
$X^{i}_{t}$ and $\gamma_{1}(t)=\gamma_{2}(t)$ for $t\not=0$.

\proclaim{Lemma 2.3}
Let $C^{1}_{\delta}$ (resp. $C^{2}_{-\delta}$) be the effective cycle 
on $X^{1}_{0}$ (resp. $X^{2}_{0}$) such that 
$\phi_{0}^{1}(C^{1}_{\delta})=\delta$
(resp. $\phi_{0}^{2}(C^{2}_{-\delta})=-\delta$). Then, $C^{1}_{\delta}$
(resp. $C^{2}_{-\delta}$) is an irreducible $-2$ curve.
\endproclaim

\demo{Proof}
We only prove the case $i=1$. Since $\langle\kappa,d\rangle\geq0$ for any 
$d\in\Delta^{+}_{1}(\eta_{0})$ and $\langle\kappa,\kappa\rangle>0$,
$(\phi_{0}^{1})^{-1}(\kappa)$ is a pseudo-ample class ([Mo,$\S 5$]).
As $\langle\kappa,\delta\rangle=0$, $C^{1}_{\delta}$ consists of a chain of 
$-2$-curves by Mayer's theorem ([Mo, $\S 5$]). Since $C^{1}_{\delta}$ is 
effective, we can write $C^{1}_{\delta}=\sum_{j}m_{j}E_{j}$ where $E_{j}$ 
is an irreducible $-2$-curve such that 
$\langle\kappa,\phi_{0}^{1}(E_{j})\rangle=0$, and $m_{i}\in\Bbb Z_{+}$. 
Suppose $C_{\delta}^{1}$ is not irreducible. Then, 
$E_{1}\subsetneq C_{\delta}^{1}$. This implies 
$\phi_{0}^{1}(E_{1})\in\Delta(\eta_{0})\cap\kappa^{\perp}$. By Lemma 2.2 
and effectivity of $E_{1}$, $\phi_{0}^{1}(E_{1})=\delta$. Thus, 
$C_{\delta}^{1}=E_{1}$ and contradicts the assumption.\qed
\enddemo

Let $\pi_{i}:\gamma_{i}^{*}\Cal X\to D$ be the pullback of the 
universal family by $\gamma_{i}$. Recall that a compact complex surface 
with at most rational double points is said to be a generalized 
$K3$ surface when its minimal resolution is a $K3$ surface.

\proclaim{Proposition 2.1}
There exists a family of generalized $K3$ surface $\pi:\Cal Y\to D$, 
a contraction morphism $b_{i}:\gamma_{i}^{*}\Cal X\to\Cal Y$ which commutes
with the projections, and a birational map 
$e:\gamma_{1}^{*}\Cal X\dashrightarrow\gamma_{2}^{*}\Cal X$ such that
\flushpar{$(1)$} 
$b_{1}$ (resp. $b_{2}$) is the blow-down of $C^{1}_{\delta}$ (resp.
$C^{2}_{-\delta}$) to a point $o$,
\flushpar{$(2)$} 
$e:X_{0}^{1}\backslash C_{\delta}^{1}\to 
X_{0}^{2}\backslash C_{-\delta}^{2}$ is an isomorphism such that 
$b_{1}=b_{2}\circ e$,
\flushpar{$(3)$} 
$e$ is the identity map on $X_{t}$ for any $t\not=0$,
\flushpar{$(4)$} 
On $X_{1,0}$ and $X_{2,0}$, 
$\phi_{1}\circ e^{*}=w_{\delta}\circ\phi_{2}$ 
where $w_{\delta}$ is the reflection in $\delta$.
\flushpar{$(5)$} 
$(\Cal Y,o)$ and $(Y_{0},o)$ are nodes of dimension $3$ and $2$
respectively.
\endproclaim

\demo{Proof}
See [Mo, $\S 3$ Cor.2] and [Be, pp.143 Remarques].\qed
\enddemo

\proclaim{Lemma 2.4}
There exists an embedding $j:\Cal Y\to\Bbb P^{N}\times D$ such that
$(1)$ $\pi=pr_{2}\circ j$ 
$(2)$ $\phi_{t}^{i}(c_{1}(b_{i}^{*}\Cal O_{\Bbb P^{N}}(1)))=m\,\kappa$ for 
some $m\in\Bbb Z_{+}$ and any fiber of $\pi_{i}$.
\endproclaim

\demo{Proof}
Let $L_{i}\to\gamma_{i}^{*}\Cal X$ be the holomorphic line bundle such that
$\phi_{t}^{i}(c_{1}(L_{i}))=\kappa$ for any fiber of $\pi_{i}$. 
There exists $m\gg1$ such that the linear system $|L_{i}|$ is 
very ample on $X_{t}^{i}$ for $t\not=0$, and is base point free and an 
embedding modulo $C_{\delta}$; 
$\Phi_{|m\,L_{i}|}\times\pi_{i}:\gamma^{*}\Cal X\hookrightarrow
\Bbb P^{N}\times D$. 
By construction and Proposition 2.1, we get
$\Cal Y=(\Phi_{|m\,L_{i}|}\times\pi_{i})(\gamma_{i}^{*}\Cal X)$ and
$L_{i}^{m}=\Phi_{|mL_{i}|}^{*}\Cal O_{\Bbb P^{N}}(1)$.\qed
\enddemo

\subsubhead
\bf{2.2 $\pmb{2}$-Elementary $\pmb{K3}$ Surfaces and $\pmb{2}$-Elementary
Lattices}\rm
\endsubsubhead
For a lattice $M$, we denote by $M^{\lor}$ its dual lattice relative to the
quadratic form of $M$. $A(M):=M^{\lor}/M$ is called the discriminant group.
A lattice $M$ is said to be $2$-elementary when 
$A(M)\cong(\Bbb Z/2\Bbb Z)^{l(M)}$ for some $l(M)\in\Bbb Z_{\geq0}$. 
We denote by $r(M)$ the rank of $M$. Let $\delta(M)$ be the parity of the
discriminant form ([Ni1,4]). By Nikulin, the triplet $(r,l,\delta)$ 
determines the isometry class of 2-elementary lattices ([Ni1]). 
\par
Let $M$ be a primitive hyperbolic 2-elementary lattice of $L_{K3}$. $N$ is 
also 2-elementary and there exists a natural isomorphism between $A(M)$ and 
$A(N)$. Consider the sublattice $L':=M\oplus N$ and the involution $I_{M}$ 
on $L'$;
$$
I_{M}(x,y)=(x,-y)\quad(x\in M,\, y\in N).
\tag 2.6
$$
As $I_{M}$ uniquely extends to an involution of $L_{K3}$, identify
$I_{M}$ with the extended one. By construction, $M$ is the fixed part and 
$N$ is the anti-fixed part;
$$
M=\{l\in L;\,I_{M}(l)=l\},\quad N=\{l\in L;\,I_{M}(l)=-l\}.
\tag 2.7
$$
$I_{M}$ induces an involution on $\tilde{\Omega}$ by 
$I(X,\phi_{X},\kappa_{X})=(X,I_{M}\circ\phi_{X},\kappa_{X})$ where 
$(X,\phi_{X},\kappa_{X})$ is a marked polarized $K3$ surface. 
Let $\Phi:R^{2}\Cal P_{*}\Bbb Z\to L_{K3}$ be the global trivialization of 
the second cohomology group. (Thus $\Phi|_{(X,\phi_{X})}=\phi_{X}$.) 
By the universality, there exists an involution 
$\iota_{\Cal X}:\Cal X\to\Cal X$ such that the following diagrams commute;
$$
\CD
\Cal X @>\iota_{\Cal X}>> \Cal X\\
@V\Cal P VV @VV\Cal P V\\
\tilde{\Omega} @>>I> \tilde{\Omega},
\endCD
\quad\quad\quad\quad
\CD
R^{2}\Cal P_{*}\Bbb Z @>\iota_{\Cal X}^{*}>> R^{2}\Cal P_{*}\Bbb Z\\
@V\Phi VV @VV\Phi V\\
L_{K3} @>>I_{M}> L_{K3}
\endCD
\tag 2.8
$$
By construction, $\tilde{\Omega}_{M}^{0}$ is contained in the fixed 
locus of $I$. Put $\iota_{M}:=\iota_{\Cal X}|_{\Cal X_{M}}$. By (2.8), 
$\iota_{M}$ is an involution over $\Cal X_{M}^{0}$ which induces an 
involution on each fiber $(X,\phi_{X})$ such that 
$\phi_{X}\circ\iota^{*}_{M}=I_{M}\circ\phi_{X}$. In this way, associated 
to the 2-elementary primitive hyperbolic lattice $M$, any fiber of 
$\Cal X_{M}^{0}$ is a $K3$ surface with an anti-symplectic involution. 
Here, an involution $\iota:X\to X$ is anti-symplectic when
$\iota^{*}\eta_{X}=-\eta_{X}$. For a pair $(X,\iota)$ of $K3$ surface and
its anti-symplectic involution, put 
$H^{2}(X,\Bbb Z)_{\pm}=\{l\in H^{2}(X,\Bbb Z);\iota^{*}_{X}(l)=\pm l\}$.

\proclaim{Definition 2.1}
A pair $(X,\iota_{X})$ of a $K3$ surface and an anti-symplectic involution 
is said to be a $2$-elementary $K3$ surface of type $M$ if there exists a
marking $\phi_{X}$ such that 
$\phi_{X}\circ\iota_{X}^{*}\circ\phi_{X}^{-1}=I_{M}$, 
or equivalently $\phi_{X}(H^{2}(X,\Bbb Z)_{+})=M$. 
Two $2$-elementary $K3$ surfaces $(X,\iota)$, $(X',\iota')$ are 
isomorphic if there  exists an isomorphism $f:X\to X'$ such that
$f\circ\iota=\iota'\circ f$. Marking as above is said to be a marking of
$2$-elementary $K3$ surfaces of type $M$.
\endproclaim

The following is clear by Theorem 2.2 and above argument. (See [Ni2-4].)

\proclaim{Theorem 2.3}
$\tilde{\Omega}_{M}^{0}$ is the fine moduli space of marked $2$-elementary 
$K3$ surface of type $M$ and
$\Cal P_{M}:(\Cal X_{M}^{0},\iota_{M})\to\tilde{\Omega}_{M}^{0}$ is the 
universal family. 
\endproclaim

Let $\Gamma(M)$ (resp. $\Gamma_{M}$) be the following subgroup of 
$O(L_{K3})$ (resp. O(N));
$$
\Gamma(M):=\{g\in O(L),\,I_{M}\circ g=g\circ I_{M}\},\quad
\Gamma_{M}:=\{g|_{N};\,g\in\Gamma(M)\}.
\tag 2.9
$$
Then, $\Gamma(M)$ acts on $\tilde{\Omega}_{M}^{0}$ and $\Gamma_{M}$ acts 
on $\Omega_{M}^{0}$ by $g\cdot[(\eta,\kappa)]:=[(g\eta,g\kappa)]$ and
$(g|_{N})\eta=g\eta$ respectively. By the global Torelli theorem, two 
2-elementary $K3$ surface of type $M$ are isomorphic if their period 
points in $\tilde{\Omega}_{M}^{0}$ lie on the same $\Gamma(M)$-orbit. 
The following theorem can be proved in the same manner as the proof of
weak Torelli theorem and surjectivity theorem for the period map of 
Enriques surfaces ([Na, Corollary 4.14 and Theorem 7.1]).

\proclaim{Theorem 2.4}
Via the period map $\pi_{M}$, 
$\tilde{\Omega}_{M}^{0}/\Gamma(M)=\Omega_{M}^{0}/\Gamma_{M}$.
\endproclaim

Since $\Gamma_{M}$ is an arithmetic subgroup of $O(2,20-r(M))$, and 
$\Omega_{M}$ is two copies of the symmetric bounded domain of type IV, 
both $\Omega_{M}/\Gamma_{M}$ and $\Omega_{M}^{0}/\Gamma_{M}$ are 
quasi-projective algebraic variety by 
Baily-Borel. Define modular varieties by
$$
\Cal M_{M}:=\Omega_{M}/\Gamma_{M},\quad
\Cal M_{M}^{0}:=\Omega_{M}^{0}/\Gamma_{M}=\tilde{\Omega}_{M}^{0}/\Gamma(M).
\tag 2.10
$$
By Theorem 2.4, $\Cal M_{M}^{0}$ is the moduli space of 2-elementary $K3$ 
surfaces of type $M$. To see what happens on the involution along 
the generic point of discriminant locus, let 
$\pi_{i}:\gamma_{i}^{*}\Cal X\to D$ be the same as in Proposition 2.1, 
and $I$, $\iota_{\Cal X}$ as in (2.8). As $I\circ\gamma_{1}=\gamma_{2}$,
it follows from the universal property of 
$\Cal P:\Cal X\to\tilde{\Omega}$ that $\iota_{\Cal X}$ induces an 
isomorphism between $(\pi_{1},\gamma_{1}^{*}\Cal X,D)$ and 
$(\pi_{2},\gamma_{2}^{*}\Cal X,D)$. Namely, 
$I\circ\pi_{1}=\pi_{2}\circ\iota_{\Cal X}$. Thus, by Proposition 2.1, 
there exist two (rational) maps $e$ and $\iota_{\Cal X}$ between 
$\gamma_{1}^{*}\Cal X$ and $\gamma_{2}^{*}\Cal X$. 
Then, $\iota_{1}:=e^{-1}\circ\iota_{\Cal X}:
\gamma_{1}^{*}\Cal X\backslash C_{\delta}^{1}\to
\gamma_{1}^{*}\Cal X\backslash C_{\delta}^{1}$ extends to a rational 
involution of $\gamma_{1}^{*}\Cal X$ which commutes with the projection.

\proclaim{Proposition 2.2}
There exists a holomorphic involution $\iota_{\Cal Y}:\Cal Y\to\Cal Y$ which
commutes with the projection $\pi:\Cal Y\to U$ such that
$\iota_{\Cal X}=b_{1}^{-1}\circ\iota_{\Cal Y}\circ b_{1}$ on 
$\gamma_{1}^{*}\Cal X\backslash C_{\delta}^{1}$.
Namely, $\pi_{1}:\gamma_{1}^{*}\Cal X|_{D\backslash\{0\}}\to 
D\backslash\{0\}$ extends to the family $\pi:\Cal Y\to D$ by contracting
$C_{\delta}^{1}$. The central fiber $Y_{0}$ is a generalized $K3$ surface 
with one node $o$ on which $\iota_{\Cal Y}$ induces an anti-symplectic 
involution. Moreover, $o$ is a fixed point of $\iota_{\Cal Y}$.
\endproclaim

\demo{Proof}
Put $\iota_{\Cal Y}:=b_{1}\circ\iota_{1}\circ b_{1}^{-1}$. As $\iota_{1}$ 
is a regular involution over 
$\gamma_{1}^{*}\Cal X\backslash C_{\delta}^{1}$, 
$\iota_{\Cal Y}$ is an involution over $\Cal Y\backslash\{o\}$ where 
$o=b_{1}(C_{\delta}^{1})$. By putting $\iota_{\Cal Y}(o)=o$, 
$\iota_{\Cal Y}$ extends to an involution over $\Cal Y$.
Since $(\Cal Y,o)$ is normal by Proposition 2.1, $\iota_{\Cal Y}$ is 
regular.
\qed
\enddemo

For a 2-elementary $K3$ surface $(X,\iota)$ of type $M$, let
$X^{\iota}$ be the fixed locus; $X^{\iota}:=\{x\in X;\,\iota(x)=x\}$
consisting of disjoint union of finitely many smooth curves. 
Nikulin determined the topological type of $X^{\iota}$ ([Ni4]). 

\proclaim{Theorem 2.5}
Let $M$ be the lattice of Nikulin type 
$(r,l,\delta)$. Then, 
$$
X^{\iota}=
\cases
(1)\quad\emptyset\qquad\qquad\qquad\qquad\quad(r,l,\delta)=(10,10,0)\\
(2)\quad C^{(1)}_{1}+C^{(1)}_{2}\qquad\qquad\,\,\,\,
(r,l,\delta)=(10,8,0)\\
(3)\quad C^{(g(M))}+\sum_{i=1}^{k(M)}E_{i}\quad(r,l,\delta)\not=
(10,10,0),\,(10,8,0)
\endcases
$$ 
where $C^{(g)}$ is a smooth curve of genus $g$, and
$E_{i}$ is a smooth $-2$-curve. In $(3)$, $g(M)=(22-r-l)/2$ and
$k(M)=(r-l)/2$. 
\endproclaim 

Let $\pi:\Cal Y\to U$ and $\iota_{\Cal Y}$ be the degenerating family of 
2-elementary $K3$ surfaces of type $M$ considered in Proposition 2.2.
Set $Y_{0}$ for the central fiber. Since the embedding dimension of
$(Y_{0},o)$ is $3$ and $o$ is a fixed point of $\iota_{\Cal Y}$, 
$(\iota_{\Cal Y*})_{o}$ induces an involution on $\Bbb C^{3}$. 
As $\iota_{\Cal Y}$ is anti-symplectic, $\det(\iota_{\Cal Y*})_{o}=-1$ 
and $(\iota_{\Cal Y*})_{o}$ is expressed by the diagonal matrix with 
eigenvalues as follows in a suitable coordinates;
$$
\hbox{Type (0,3)}:\quad(-1,-1,-1),\quad
\hbox{Type (2,1)}:\quad(1,1,-1).
\tag 2.11
$$
As $(\Cal Y,o)$ is a node, by a careful look at the morsification 
procedure, we get the following.

\proclaim{Proposition 2.3}
Let $\Cal Z$ be the fixed locus of $(\Cal Y,\iota_{\Cal Y})$.
\flushpar{$(1)$} 
In case of type $(0,3)$, $o$ is an isolated point of $\Cal Z$ and there 
exists a neighborhood $V$ of $o$ such that $\pi(x,y,z,t)=t$ and
$\iota_{\Cal Y}(x,y,z,t)=(-x,-y,-z,t)$,
\newline
$(\Cal Y|_{V},o)=\{(x,y,z,t);\,xy-z^{2}-t^{2}=0\}$.
\flushpar{$(2)$} 
In case of type $(2,1)$, $\Cal Z$ has a unique node at $o$ and there 
exists a neighborhood $V$ of $o$ such that $\pi(x,y,z,t)=t$ and
$\iota_{\Cal Y}(x,y,z,t)=(-x,-y,-z,t)$,
\newline
$(\Cal Y|_{V},o)=\{(x,y,z,t);\,xy-z^{2}-t^{2}=0\}$
\endproclaim

For a 2-elementary primitive hyperbolic lattice $M$ and 
$\delta\in\Delta(N)$, set $\langle M\oplus\delta\rangle$ for the smallest 
2-elementary primitive hyperbolic lattice containing $M$ and $\delta$. 
It follows from definition that 
$\Delta(\langle M\oplus\delta\rangle^{\perp})=\delta^{\perp}\cap\Delta(N)$ 
and $H_{\delta}\cap\Omega_{M}=\Omega_{\langle M\oplus\delta\rangle}$.

\proclaim{Lemma 2.5}
For the family $\pi:\Cal Y\to D$ in Proposition $2.2$, 
$\iota_{0}^{1}:=b_{1}^{-1}(\iota_{\Cal Y}|_{Y_{0}})b_{1}$ extends to a
involution on the minimal resolution $X_{0}^{1}$, and 
$(X_{0},\iota_{0}^{1})$ becomes a $2$-elementary $K3$ surface of type 
$\langle M\oplus\delta\rangle$.
\endproclaim

\demo{Proof}
By Proposition 2.3, it is clear that $\iota_{0}^{1}$ extends to an 
involution on $X_{0}^{1}$. Let $\phi_{1}$ and $\phi_{2}$ be the marking 
as in Proposition 2.1. By Proposition 2.1 (4) and Definition 2.1, we get
$\phi_{1}^{-1}(\iota_{0}^{1})^{*}\phi_{1}=
(\phi_{1}^{-1}\iota_{M}^{*}\phi_{2})\circ(\phi_{2}^{-1}(e^{-1})^{*}\phi_{1})
=I_{M}\circ w_{\delta}$. Since $\delta\perp M$, we get
$I_{M}\circ w_{\delta}=w_{\delta}\circ I_{M}$ and 
$I_{M}\circ w_{\delta}=I_{\langle M\oplus\delta\rangle}$. 
\qed
\enddemo

Let $\Cal A_{g}:=\frak S_{g}/ Sp(2g;\Bbb Z)$ be the Siegel modular variety. 
(When $g=0$, $\Cal A_{g}$ is a point.) We denote by $\Cal A^{*}_{g}$ and 
$\Cal M_{M}^{*}$ the Satake-Baily-Borel compactification. 
Let $C=\sum_{i=1}^{l}C_{i}$ be a disjoint union of smooth 
curves and $\Jac(C_{i})$ the Jacobi variety of $C_{i}$. Put 
$[\Jac(C)]:=([\Jac(C_{1})],\cdots,[\Jac(C_{l})])\in
\Cal A_{g_{1}}\times\cdots\times\Cal A_{g_{l}}$ where $[\Jac(C_{i})]$ is the
period of $\Jac(C_{i})$. By Theorems 2.4 and 2.5, we can define 
a morphism $j_{M}^{0}:\Cal M_{M}^{0}\to\Cal A_{g(M)}$ as follows.

\proclaim{Definition 2.2}
$j_{M}^{0}:\Cal M_{M}^{0}\ni[(X,\iota_{X})]\longrightarrow[\Jac(X^{\iota})]
\in\Cal A_{g(M)}$. 
When $(r,l,\delta)=(10,10,0)$, $j_{M}^{0}$ is defined to be the constant 
map. When $(r,l,\delta)=(10,8,0)$, $j_{M}^{0}$ takes its value in 
$S^{2}(\Cal A_{1})$, the second symmetric product of $\Cal A_{1}$.
\endproclaim

\proclaim{Proposition 2.4}
$j_{M}^{0}$ extends to a rational map 
$j_{M}:\Cal M_{M}^{*}\dashrightarrow\Cal A_{g(M)}^{*}$ if $r\leq 17$.
\endproclaim

\demo{Proof}
Since $\Cal M_{M}^{*}$ and $\Cal A_{g(M)}^{*}$ are projective 
algebraic varieties, it is enough to show that $j_{M}$ extends to a 
morphism from an open subset $V$ of $\Cal M_{M}$ such that 
$\codim\,\Cal M_{M}\backslash V\geq2$ because the boundary components of 
$\Cal M_{M}^{*}$ have codimension $\geq2$ when $r(M)\leq 17$. Put 
$V=\Cal M_{M}^{0}\cup\bigcup_{\delta\in\Delta(N)}H_{\delta}^{0}/\Gamma_M$.
Take a point $\eta\in H_{\delta}^{0}$ and its small neighborhood $U$. By 
Proposition 2.2, we have a family of 2-elementary $K3$ surfaces of type 
$M$ whose fiber over $U\cap H_{\delta}$ is a generalized $K3$ surface with 
an anti-symplectic involution. Let $(X_{\eta},\iota_{\eta})$ be the 
corresponding generalized $K3$ surface with involution over $\eta$. 
Let $X^{\iota}_{\eta}$ be the fixed locus and $\hat{X}^{\iota}_{\eta}$ its
normalization. Define 
$j_{M}(X_{\eta},\iota_{\eta}):=[\Jac(\hat{X}^{\iota}_{\eta})]
\in\Cal A_{g(M)}^{*}$. 
Since $X^{\iota}_{\eta}$ has at most a node, it is well defined and 
gives an extension of $j_{M}$ to $U$. As $j_{M}$ is defined on 
$U\backslash H_{\delta}$, its extension to $U$ is  unique. This proves 
that $j_{M}^{0}$ extends to $V$.\qed
\enddemo

\proclaim{Proposition 2.5}
For $\delta\in\Delta(N)$, 
$j_{M}|_{H_{\delta}^{0}}=j_{\langle M,\delta\rangle}$.
\endproclaim

\demo{Proof}
Take a point $\eta_{0}\in H^{0}_{\delta}$ and its small neighborhood $U$ in
$\Omega_{M}$. Let $\tilde{U}\subset\tilde{\Omega}_{M}$ be a lift of $U$ 
such that $\tilde{\Omega}_{M}^{0}$ is dense in $\tilde{U}$. Let 
$\pi':\Cal Y\to\tilde{U}$ and $\pi':\Cal Z\to\tilde{U}$ be the family of
2-elementary $K3$ surfaces and its fixed locus. By the period map,
identify $U$ with $\tilde{U}$. Then, over $U\cap H^{0}_{\delta}$, 
a fiber of $\Cal Y$ is a generalized $K3$ surface  with an anti-symplectic
involution. Take $\tau\in U\cap H_{\delta}^{0}$. Let
$(Y_{\tau},\iota_{\tau})$ be the fiber over $\tau$ and $Z_{\tau}$ its  
fixed locus. Let $p:\tilde{Y}_{\tau}\to Y_{\tau}$ be the minimal 
resolution of $Y_{\tau}$. By Lemma 2.5, 
$\tilde{\iota}_{\tau}:=p^{-1}\circ\iota\circ p$ extends to an involution 
over $\tilde{Y}_{\tau}$. Set $\tilde{Z}_{\tau}$ for the fixed locus of 
$(\tilde{Y}_{\tau},\tilde{\iota}_{\tau})$. By Proposition 2.3, 
$\tilde{Z}_{\tau}=Z_{\tau}+E$ when $(\Cal Y,o)$ is of type $(0,3)$, and 
$\tilde{Z}_{\tau}=\hat{Z}_{\tau}$ is the normalization
when $(\Cal Y,o)$ is of type $(2,1)$. By definition, we find 
that $j_{\langle M,\delta\rangle}=j|_{H^{0}_{\delta}}$ on $U$.
\qed
\enddemo

\beginsection 
$\S3$. An Invariant of $\pmb{2}$-Elementary $\pmb{K3}$ Surfaces

\subsubhead
\bf{3.1 Determinant Bundles and Quillen Metrics}
\endsubsubhead
Let $(M,g)$ be a compact K\"ahler manifold on which acts a finite group $G$ 
holomorphically and isometrically. Let $\square_{0,q}^{G}$ be the Laplacian
restricted to the space of $G$-invariant $(0,q)$-forms and 
$\zeta_{0,q}^{G}(s)$ be the spectral zeta function of $\square^{G}_{0,q}$,
which is regular at $s=0$. $\zeta^{G}_{0,q}(s)$ is nothing but the spectral 
zeta function of the orbifold $(M/G,g)$.

\proclaim{Definition 3.1}
The Ray-Singer analytic torsion of $(M/G,g)$ is defined by
$$
\tau(M/G,g):=\prod_{q\geq 0}(\det\,\square_{0,q}^{G})^{(-1)^{q}q},\quad
\det\,\square_{0,q}^{G}:=
\exp\left(-\left.\frac{d}{ds}\right|_{s=0}\zeta_{0,q}^{G}(s)\right).
$$
\endproclaim

Let $\pi:X\to S$ be a proper smooth morphism of K\"ahler manifolds on which
acts the finite group $G$ holomorphically. Assume that $G$ preserves the 
fiber of $\pi$. An equivariant determinant of cohomology is a line bundle
on $S$ defined by $\lambda(X/S,G):=\bigotimes_{q\geq 0}
\left(\det\,R^{q}_{G}\pi_{*}\Cal O_{X}\right)^{(-1)^{q}}$.
Here the fiber at $s$ of the sheaf $R^{q}_{G}\pi_{*}\Cal O_{S}$ is the 
equivariant cohomology group 
$H^{q}(X_{s},\Cal O_{X_{s}})^{G}=H^{q}(X_{s}/G,\Cal O_{X_{s}/G})$. In the 
sequel, we assume that a $G$-invariant K\"ahler metrics $g_{X/S}$ is 
equipped on the relative tangent bundle $TX/S:=\ker\pi_{*}$. Namely, 
$g_{X/S}$ is a family of $G$-invariant K\"ahler metric 
$g_{t}:=g_{X/S}|_{X_{t}}$. Via the Hodge theory, the fiber of $\pi$ is 
identified with the determinant of $G$-invariant harmonic forms which 
induces the $L^{2}$-metric on $\lambda(X/S,G)$. We denoted it by 
$\|\cdot\|_{L^{2}}^{G}$. 

\proclaim{Definition 3.2}
The Quillen metric of $\lambda(X/S,G)$ relative to $g_{X/S}$ is defined 
by
$\|\cdot\|_{Q}^{2}(t):=\tau(X_{t}/G,g_{t})\cdot\|\cdot\|^{G}_{L^{2}}(t)$.
\endproclaim

Let $X_{g}=\{x\in X;\,g(x)=x\}$ be the fixed locus of $g$, $R_{TX_{g}/S}$
the curvature form of $(TX_{g}/S,g_{X/S}|_{TX_{g}/S})$, $R_{N_{X_{g}/X}}$
the curvature form of the relative normal bundle
$(N_{X_{g}/X},g_{X/S}|_{N_{X_{g}/X}})$, and $\exp(i\theta_{j})$ an 
eigenvalue of $g_{N_{X_{g}/X}}$. Let $\Td(A)$ and $e(A)$ be the Todd and
Euler genuses i.e., $\Td(A)=\det(A(1-e^{-A})^{-1})$ and $e(A)=\det(A)$.
The following theorems are due to Bismut ([Bi]) and Bismut-Gillet-Soul\'e 
([B-G-S]). (Although detailed proof of Theorem 3.1 is not written in 
[B-G-S], it can be proved in the similar way as in the case $G=\{1\}$. 
Details are left to the reader. We refer to K\"ohler-Roessler's 
paper ([K-R]) for further generalizations.)

\proclaim{Theorem 3.1}
The Chern form of $\lambda(X/S,G)_{Q}:=(\lambda(X/S,G),\|\cdot\|_{Q})$ is 
given by
$$
c_{1}(\lambda(X/S,G)_{Q})=
\frac{1}{|G|}\sum_{g\in G}\pi_{*}(\Td_{g}(TX/S,g_{X/S}))^{(1,1)}.
$$
Here $\Td_{g}(TX/S,g_{X/S})$ is the $g$-Todd genus of $(TX/S,g_{X/S})$;
\newline
$\Td_{g}(TX/S,g_{X/S})=
\Td\left(\frac{i}{2\pi}R_{TX_{g}/S}\right)\prod_{j=1}^{q}
\left(\Td/\e\right)\left(\frac{i}{2\pi}R_{N_{X_{g}/X}}+
i\,\theta_{j}\right)$.
\endproclaim

\proclaim{Theorem 3.2}
Let $g_{X/S}$ and $g_{X/S}'$ be $G$-invariant K\"ahler metrics on $TX/S$.
Let $\|\cdot\|_{Q}$ and $\|\cdot\|_{Q}'$ be the Quillen metrics of 
$\lambda(X/S,G)$ relative to $g_{X/S}$ and $g_{X/S}'$. Then, 
$$
\log\left(\frac{\|\cdot\|'_{Q}}{\|\cdot\|_{Q}}\right)^{2}=
\frac{1}{|G|}\sum_{g\in G}
\pi_{*}(\widetilde{\Td_{g}}(TX/S;g_{X/S},g'_{X/S}))^{(0,0)}
$$
where $\widetilde{\Td_{g}}(TX/S;g_{X/S},g'_{X/S})$ is the Bott-Chern 
secondary class associated to the $g$-Todd form and $g_{X/S}$, $g'_{X/S}$.
\qed
\endproclaim

\subsubhead
\bf{3.2 An Invariant via Analytic Torsion}\rm
\endsubsubhead
Let $(X,\iota,\kappa)$ be a 2-elementary $K3$ surface with 
$\iota$-invariant Ricci-flat K\"ahler metric $\kappa$ and symplectic form
$\eta_{X}$. The space of $L^2$-$(0,q)$-forms splits into 
$L^{2}(X,\wedge^{0,q})_{\pm}=\{f\in L^{2}(X,\wedge^{0,q});\,
\iota^{*}f=\pm f\}$. Let $\square^{0,q}_{\pm}$ be the Laplacian restricted 
to $L^{2}(X,\wedge^{0,q})_{\pm}$. By definition, the Ray-Singer analytic 
torsion of $(X/\iota,\kappa)$ is given by 
$\tau(X/\iota,\kappa)=\det\square_{+}^{0,2}/\det\square_{+}^{0,0}$
because $\zeta_{+}^{0,0}(s)-\zeta_{+}^{0,1}(s)+\zeta_{+}^{0,2}(s)\equiv 0$
where $\zeta_{\pm}^{0,q}(s)$ is the zeta function of $\square_{\pm}^{0,q}$. 
Since $\eta_{X}$ is a parallel section such that 
$\iota^{*}\eta_{X}=-\eta_{X}$, multiplication by $\eta_{X}$ induces an 
isometry between $(L^{2}(X,\wedge^{0,0})_{-},\Delta^{0,0}_{-})$ and 
$(L^{2}(X,\wedge^{0,2})_{+},\Delta^{0,2}_{+})$. Thus, we get the following.

\proclaim{Lemma 3.1}
$\tau(X/\iota,\kappa)=\det\square_{-}^{0,0}/\det\square_{+}^{0,0}$.
\endproclaim

Let $p:(X,\iota)\to S$ be a family of 2-elementary $K3$ surfaces of type
$M$ and $p:Z\to S$ be the fixed locus. By Theorem 2.5, $Z$ decomposes into
$Z=\sum_{i}C_{i}$ where $p:C_{i}\to S$ is a family of curves over $S$. Let
$g_{X/S}$ be a family of Ricci-flat K\"ahler metrics such that 
$\iota^{*}g_{X/S}=g_{X/S}$. Set $g_{C_{i}/S}:=g_{X/S}|_{C_{i}/S}$ for the 
induced metric on $TC_{i}/S$ and $g_{N_{C_{i}/X}}$ for the induced metric 
on the relative normal bundle $N_{C_{i}/X}$  which is defined by the 
following exact sequence;
$$
0\longrightarrow TC_{i}/S\longrightarrow \left.TX/S\right|_{C_{i}}
\longrightarrow N_{C_{i}/X}\longrightarrow0.
\tag 3.1
$$
We define four smooth functions over $S$ as follows;
$$
\aligned
\,
&\tau(X/\iota/S,g_{X/S})(s):=\tau(X_{s}/\iota,g_{s}),\quad
\tau(C_{i}/S,g_{C_{i}/S})(s):=\tau(C_{i,s},g_{C_{i,s}}),\\
&\volum(X/S,g_{X/S})(s):=\volum(X_{s},g_{s}),\quad
\volum(C_{i}/S)(s):=\volum(C_{i,s},g_{C_{i,s}}).
\endaligned
\tag 3.2
$$
By construction, the following is clear.

\proclaim{Lemma 3.2}
$$
c_{1}(X/S,g_{X/S})|_{C_{i}/S}=c_{1}(C_{i}/S,g_{C_{i}/S})+
c_{1}(N_{C_{i}/S},g_{N_{C_{i}/S}}).
$$
\endproclaim

\proclaim{Proposition 3.1}
$$
\left[\Td_{\iota}(X/S,g_{X/S})\right]^{(2,2)}=
\sum_{i}\left\{\frac{1}{8}\,c_{1}(X/S,g_{X/S})\,c_{1}(C_{i}/S,g_{C_{i}/S})-
\frac{1}{12}\,c_{1}(C_{i}/S,g_{C_{i}/S})^{2}\right\}.
$$
\endproclaim

\demo{Proof}
Put $x=c_{1}(C_{i}/S,g_{C_{i}/S})$ and 
$y=c_{1}(X/S,g_{X/S})-c_{1}(C_{i}/S,g_{C_{i}/S})$. By definition,
$(\Td/\e)(y+\pi\,i)=(1+\exp(-y))^{-1}$. Since
$$
\Td(x)=1+\frac{x}{2}+\frac{x^{2}}{12}+O(x^{3}),\quad
\frac{\Td}{\e}(y+\pi\,i)=\frac{1}{2}+\frac{y}{4}+O(y^{3}),
\tag 3.3
$$
we get the assertion by the definition of equivariant Todd form in 
Theorem 3.1 together with (3.3).\qed
\enddemo

Let $\eta_{t}$ be a symplectic form on $X_{t}$ depending holomorphically
in $t\in S$ and $\kappa_{t}$ be a Ricci-flat K\"ahler metric on $X_{t}$
depending smoothly in $t$. Put 
$\|\eta_{t}\|^{2}=\int_{X_{t}}\eta_{t}\wedge\overline{\eta}_{t}$.

\proclaim{Lemma 3.3}
$$
\kappa_{t}^{2}=f(t)\,\eta_{t}\wedge\overline{\eta}_{t},\quad
f(t)=\volum(X_{t},\kappa_{t})/\|\eta_{t}\|^{2}.
$$
\endproclaim

\demo{Proof}
Since $\kappa_{t}$ is Ricci-flat, there exists a constant $C_{t}$ such that
$\kappa_{t}^{2}=C_{t}\,\eta_{t}\wedge\overline{\eta}_{t}$. Integrating both
hand sides, we get the assertion.\qed
\enddemo

Let $\Omega_{M}$ be the period domain of 2-elementary $K3$ surfaces of 
type $M$. Let $\omega_{M}$ be its Bergman metric and $K_{M}(z,\bar{z})$ 
its Bergman kernel function;
$$
\omega_{M}(z)=\frac{i}{2\pi}\partial\bar{\partial}\log K_{M}(z,\bar{z}),
\quad K_{M}(z,\bar{z})=
\frac{\langle z,\bar{z}\rangle_{N}}{|\langle z,l_{M}\rangle_{N}|^{2}}
\tag 3.4
$$
where $l_{M}\in N_{\Bbb C}$ is a vector such that 
$H_{l_{M}}$ is the hyperplane at infinity of $\Bbb P(N_{\Bbb C})$. Let 
$\pi:S\to\Omega_{M}$ be the period map. Since
$\pi^{*}\omega_{M}=
\frac{i}{2\pi}\bar{\partial}_{M}\partial_{M}\log\|\eta_{t}\|^{2}$ by
Schumacher's formula ([Sch], [Ti], [To2]), we get the following by taking
$\frac{i}{2\pi}\partial_{M}\bar{\partial}_{M}$ of Lemma 3.3.

\proclaim{Proposition 3.2}
$$
c_{1}(X/S,g_{X/S})=-p^{*}\left\{\pi^{*}\omega_{M}+
\frac{i}{2\pi}\bar{\partial}_{M}\partial_{M}
\log\volum(X/S,g_{X/S})\right\}.
$$
\endproclaim

\proclaim{Proposition 3.3}
$$
\aligned
\,&\frac{1}{2}\left[p_{*}\Td(X/S,g_{X/S})+p_{*}\Td_{\iota}(X/S,g_{X/S})
\right]^{(1,1)}\\
&=\frac{r(M)-6}{8}\left\{-\pi^{*}\omega_{M}+
\frac{i}{2\pi}\bar{\partial}\partial\log\volum(X/S)\right\}
-\frac{1}{2}\sum_{i}p_{*}\Td(C_{i}/S,g_{C_{i}/S})^{(1,1)}.
\endaligned
$$
\endproclaim

\demo{Proof}
Since $\chi(X_{t})=24$, it follows from Proposition 3.2 and the projection
formula that
$$
\aligned
p_{*}\Td(X/S,g_{X/S})^{(1,1)}
&=\frac{1}{24}
\left[p_{*}c_{1}(X/S,g_{X/S})c_{2}(X/S,g_{X/S})\right]^{(1,1)}\\
&=-\pi^{*}\omega_{M}+
\frac{i}{2\pi}\bar{\partial}_{M}\partial_{M}\log\volum(X/S).
\endaligned
\tag 3.5
$$
As $\int_{C_{i,t}}c_{1}(C_{i,t})=1-g(C_{i,t})$, we get by Proposition 3.1 
and the projection formula,
$$
\aligned
\,&p_{*}\Td_{\iota}(X/S,g_{X/S})^{(1,1)}\\
&=p_{*}\left[\sum_{i}\frac{1}{8}c_{1}(C_{i}/S,g_{C_{i}/S})
p^{*}(-\pi^{*}\omega_{M}+
\frac{i}{2\pi}\bar{\partial}_{M}\partial_{M}\log\volum(X/S))
-\frac{1}{12}c_{1}(C_{i}/S,g_{C_{i}/S})^{2}\right]\\
&=\sum_{i}\frac{1-g(C_{i})}{8}\left\{-\pi^{*}\omega_{M}+
\frac{i}{2\pi}\bar{\partial}_{M}\partial_{M}\log\volum(X/S)\right\}-
\sum_{i}p_{*}\Td(C_{i}/S,g_{C_{i}/S})^{(1,1)}.
\endaligned
\tag 3.6
$$
It follows from Theorem 2.4 that 
$1/2+\sum_{i}(1-g(C_{i}))/8=(r(M)-6)/8$ for any $M$
which, together with (3.5-6), yields the assertion.\qed
\enddemo

\proclaim{Definition 3.3}
For a $2$-elementary $K3$ surface of type $M$ with an $\iota$-invariant 
Ricci-flat K\"ahler metric $\kappa$, we define $\tau_{M}$ by
$$
\tau_{M}(X,\iota,\kappa):=\langle\kappa,\kappa\rangle^{\frac{14-r(M)}{8}}
\tau(X/\iota,\kappa)
\prod_{i}\langle C_{i},\kappa\rangle^{\frac{1}{2}}
\tau(C_{i},\kappa|_{C_{i}})^{\frac{1}{2}}.
$$
\endproclaim

Let $\tau_{M/S}$ be the function defined by 
$\tau_{M/S}(s)=\tau_{M}(X_{s},\iota_{s},g_{X_{s}})$.
Let $j_{C_{i}/S}:S\ni s\to [\Jac(C_{i,s})]\in\Cal A_{g(C_{i})}$ be the 
period map. Let $\omega_{\Cal A_{g}}$ be the Bergman metric of $\Cal A_{g}$ 
and $K_{\Cal A_{g}}(\tau)$ the Bergman kernel of $\Cal A_{g}$;
$$
\omega_{\Cal A_{g}}(\tau)=
\frac{i}{2\pi}\bar{\partial}_{\tau}\partial_{\tau}\log K_{\Cal A_{g}}(\tau),
\quad K_{\Cal A_{g}}(\tau)=\det\Im\tau.
\tag 3.7
$$

\proclaim{Theorem 3.3}
$$
\frac{i}{2\pi}\bar{\partial}_{S}\partial_{S}\log\tau_{M/S}
=-\frac{r(M)-6}{8}\pi^{*}\omega_{M}
-\frac{1}{2}\sum_{i}j_{C_{i}/S}^{*}\omega_{\Cal A_{g(M)}}.
$$
\endproclaim

\demo{Proof}
Consider the family $p:X\to S$. Put $G=\{1,\iota\}$. Since 
$H^{q}(X_{s},\Cal O_{X_{s}})^{G}=0$ for $q>0$, it follows that 
$\lambda(X/S,\iota)=\Cal O_{S}\,1_{S}$ where $1_{S}$ is the canonical 
section of $R^{0}_{G}p_{*}\Cal O_{X}$ such that $1_{S}(s)=1$ in 
$H^{0}(X_{s},\Cal O_{X_{s}})^{G}$. Since the squared Quillen norm of 
$1_{S}$ is $\tau(X/\iota/S)\,\volum(X/S)$, it follows from Theorem 3.1 
and Proposition 3.3 that
$$
\aligned
\,&
\frac{i}{2\pi}\bar{\partial}\partial\log\tau(X/\iota/S)\,\volum(X/S)\\
&=\frac{r(M)-6}{8}\left\{\pi^{*}\omega_{M}+
\frac{i}{2\pi}\bar{\partial}\partial\log\volum(X/S)\right\}
-\frac{1}{2}\sum_{i}p_{*}\Td(C_{i}/S,g_{C_{i}/S})^{(1,1)}.
\endaligned
\tag 3.8
$$
Fix a marking of the family of curves $p:C_{i}\to S$. Let 
$T(\Jac(C_{i,s}))\in\frak S_{g(C_{i})}$ be the period matrix of 
$\Jac(C_{i,s})$ relative to this marking. Let $T(C_{i}/S)$ be a 
function with values in $\frak S_{g(C_{i})}$ defined by
$T(C_{i}/S)(s)=T(\Jac(C_{i,s}))$. Theorem 3.1 applied to $p:C_{i}\to S$
yields
$$
\frac{i}{2\pi}\bar{\partial}\partial\log\left\{\tau(C_{i}/S)\,
\volum(C_{i}/S)\det\Im\,T(J(C_{i}/S))\right\}=
p_{*}\Td(C_{i}/S,g_{C_{i}/S})^{(1,1)}.
\tag 3.9
$$
which, together with (3.8) and
$\frac{i}{2\pi}\bar{\partial}\partial\log\det\Im\,T(C_{i}/S)=
j^{*}\omega_{\Cal A_{g(C_{i})}}$, proves the assertion.\qed
\enddemo

\proclaim{Theorem 3.4}
$\tau_{M}(X,\iota,\kappa)$ does not depend on $\kappa$, and becomes 
an invariant of a $2$-elementary $K3$ surface $(X,\iota)$.
\endproclaim

\demo{Proof}
Let $\kappa_{0}$ and $\kappa_{1}$ be two $\iota$-invariant Ricci-flat
K\"ahler metrics of $(X,\iota)$, and set $T(X,\iota,\kappa):=
\volum(X,\kappa)^{(r(M)-6)/8}\tau_{M}(X,\iota,\kappa)$.
Let $1_{X}\in H^{0}(X,\Cal O_{X})=\lambda(X,\iota)$, 
$1_{C_{i}}\in H^{0}(C_{i},\Cal O_{C_{i}})$ and 
$\sigma_{i}\in\det H^{0}(C_{i},\Omega_{C_{i}}^{1})$ be the generators 
of each line. Let $\|\cdot\|_{Q,\kappa_{j}}$ be the Quillen metric
relative to $\kappa_{j}$. Since the $L^{2}$-metric on 
$\det H^{0}(C_{i},\Omega_{C_{i}}^{1})$ is independent of a choice
of K\"ahler metric, it follows from Theorem 3.2 that
$$
\aligned
\,&
\log\frac{T(X,\iota,\kappa_{0})}{T(X,\iota,\kappa_{1})}=
\log\frac{\|1_{X}\|_{Q,\kappa_{0}}^{2}}
{\|1_{X}\|_{Q,\kappa_{1}}^{2}}+\frac{1}{2}\sum
\log\frac{\|1_{C_{i}}\otimes\sigma_{i}\|_{Q,\kappa_{0}}^{2}}
{\|1_{C_{i}}\otimes\sigma_{i}\|_{Q,\kappa_{1}}^{2}}\\
&=\frac{1}{2}\int_{X}\widetilde{\Td}(X;\kappa_{0},\kappa_{1})+
\frac{1}{2}\sum\int_{C_{i}}
\widetilde{\Td}_{\iota}(X;\kappa_{0},\kappa_{1})+\sum
\int_{C_{i}}\widetilde{\Td}(C_{i};\kappa_{0},\kappa_{1}).
\endaligned
\tag 3.10
$$
Since
$\widetilde{\Td}(X)=\widetilde{c_{1}(X)\,c_{2}(X)}/24$,
$\widetilde{\Td}_{\iota}(X)=
\sum_{i}\widetilde{c_{1}(X)\,c_{2}(C_{i})}/8-
\widetilde{c_{1}(C_{i})^{2}}/12$, and
$\widetilde{\Td}(C_{i})=\widetilde{c_{1}(C_{i})^{2}}/12$,
we get
$$
\log\frac{T(X,\iota,\kappa_{0})}{T(X,\iota,\kappa_{1})}=
\frac{1}{48}\int_{X}\widetilde{c_{1}c_{2}}(X)(\kappa_{0},\kappa_{1})
+\frac{1}{16}\sum\int_{C_{i}}\widetilde{c_{1}(X)\,c_{1}(C_{i})}
(\kappa_{0},\kappa_{1}).
\tag 3.11
$$
By Yau's theorem ([Ya]), there exists a family of Ricci-flat K\"ahler 
metrics $\kappa_{t}$ $(0\leq t\leq 1)$ joining $\kappa_{0}$ and
$\kappa_{1}$. Let $\eta$ be a fixed symplectic form on $X$.
As in Lemma 3.3, put $f(t)=\volum(X,\kappa_{t})/\|\eta\|^{2}$.
Let $L_{t}=\kappa_{t}^{-1}\dot{\kappa}_{t}$ be a section of 
$\End(TX)$ such that 
$\kappa_{t}(L_{t}u,v)=\dot{\kappa}_{t}(u,v)$. By Lemma 3.3, we get
$\Tr\,\kappa_{t}^{-1}\dot{\kappa}_{t}=\partial_{t}\log f(t)$ which, 
together with the Bott-Chern formula ([B-C], [B-G-S, I, e)]), yields
$$
\aligned
\int_{X}\widetilde{c_{1}c_{2}}(X)(\kappa_{0},\kappa_{1})
&=\int_{X}\int_{1}^{0}dt\left.\frac{d}{d\epsilon}
\right|_{\epsilon=0}
c_{1}(R_{t}+\epsilon\kappa_{t}^{-1}\dot{\kappa}_{t})
c_{2}(R_{t}+\epsilon\kappa_{t}^{-1}\dot{\kappa}_{t})\\
&=\int^{0}_{1}dt\frac{d}{dt}\log f(t)\int_{X}c_{2}(R_{t})=
24\log\frac{\volum(X,\kappa_{0})}{\volum(X,\kappa_{1})},
\endaligned
\tag 3.12
$$
$$
\aligned
\int_{C_{i}}\widetilde{c_{1}(X)\,c_{1}(C_{i})}
(\kappa_{0},\kappa_{1})
&=\int_{C_{i}}\int_{1}^{0}dt\left.\frac{d}{d\epsilon}
\right|_{\epsilon=0}
c_{1}(R_{t}+\epsilon\kappa_{t}^{-1}\dot{\kappa}_{t})\,
c_{1}\left(R(C_{i})+\epsilon\,\frac{\dot{\kappa}_{t}|_{C_{i}}}
{\kappa_{t}|_{C_{i}}}\right)\\
&=\int_{1}^{0}dt\frac{d}{dt}\log f(t)\int_{C_{i}}
c_{1}(R_{t}(C_{i}))=\chi(C_{i})
\log\frac{\volum(X,\kappa_{0})}{\volum(X,\kappa_{1})}
\endaligned
\tag 3.13
$$
where $R_{t}$ is the curvature of $(X,\kappa_{t})$ and
$R_{t}(C_{i})$ of $(C_{i},\kappa_{t}|_{C_{i}})$. By (3.11-13), we get
$$
\log\frac{T(X,\iota,\kappa_{0})}{T(X,\iota,\kappa_{1})}=
\frac{r(M)-6}{8}\log
\frac{\volum(X,\kappa_{0})}{\volum(X,\kappa_{1})}
\tag 3.14
$$
which, together with the definition of $T(X,\iota,\kappa)$, yields the 
assertion.\qed
\enddemo

Let $\tau_{M}$ be the function on $\tilde{\Omega}_{M}^{0}$ defined by
$\tau_{M}(X,\phi,\iota):=\tau_{M}(X,\iota)$. By Theorems 3.3 and 3.4, 
$\tau_{M}$ is a smooth $\Gamma(M)$-invariant function on 
$\tilde{\Omega}_{M}^{0}$, and thus descends to a $\Gamma_{M}$-invariant 
smooth function (denoted by the same symbol $\tau_{M}$) on $\Omega_{M}^{0}$ 
by Theorem 2.4. Applying Theorem 3.3 to the universal family, we get the 
following.

\proclaim{Theorem 3.5}
$\tau_{M}$ descends to a smooth $\Gamma_{M}$-invariant function on 
$\Omega_{M}^{0}$ and satisfies the following variational formula: 
$$
\frac{i}{2\pi}\bar{\partial}\partial\log\tau_{M}
=-\frac{r(M)-6}{8}\omega_{M}-
\frac{1}{2}j_{M}^{*}\omega_{\Cal A_{g(M)}}.
$$
\endproclaim

\beginsection
$\S4$. Degeneration of $\pmb{K3}$ Surfaces and Monge-Amp\`ere Equation

\par
\subsubhead
\bf{4.1 Apriori Estimates for the Monge-Amp\`ere Equation}\rm
\endsubsubhead
Let $(X,\kappa)$ be a compact K\"ahler surface. Let $F\in C^{\infty}(X)$ be 
a given function and $\varphi\in C^{\infty}(X)$ satisfies the following 
complex Monge-Amp\`ere equation;
$$
\left(\kappa+\frac{i}{2\pi}\partial\bar{\partial}\varphi\right)^{2}=
e^{F}\kappa^{2},\quad\int_{X}\varphi\,\kappa^{2}=0,\quad
\kappa+\frac{i}{2\pi}\partial\bar{\partial}\varphi>0.
\tag 4.1
$$
Put $\kappa':=\kappa+\frac{i}{2\pi}\partial\bar{\partial}\varphi$ for a new
K\"ahler metric. Let $\Delta$ (resp. $\Delta'$) be the Laplacian relative 
to $\kappa$ (resp. $\kappa'$). Let $R=(R_{i\bar{j}k\bar{l}})$ be the 
holomorphic bisectional curvature of $(X,\kappa)$. Let $V$ be the volume of 
$(X,\kappa)$, $\lambda>0$ the first eigen value of $\Delta$, 
and $S$ the Sobolev constant of $(X,\kappa)$;
$\|f\|_{L^{4}}\leq S(\|df\|_{L^{2}}+\|f\|_{L^{2}})$ 
($\forall f\in C^{\infty}(X)$). For $f\in C^{0}(X)$, we denote by 
$|f|_{\infty}$ the sup-norm.

\proclaim{Proposition 4.1}
\flushpar{$(1)$} 
If $|F|_{\infty}\leq1$, 
$|\varphi|_{\infty}\leq C(S,\lambda^{-1},V)\,|F|_{\infty}$
where $C(x,y,z)$ is bounded if all of $x,y,z$ are bounded.
\flushpar{$(2)$} 
There exist 
$C_{i}=C_{i}(|\Delta F|_{\infty}/|R|_{\infty},|F|_{\infty},
|R|_{\infty}|F|_{\infty},S,\lambda^{-1},V)$ $(i=1,2)$ such that 
$C_{i}(x,y,z,w,s,t)$ is bounded from above and below if all of
$x,y,z,w,s,t$ are bounded, and the following inequality holds on $X$;
$C_{1}\,\kappa\leq\kappa'\leq C_{2}\,\kappa$.
\endproclaim

\demo{Proof}
See [Ko, pp.298-299] for (1) and [Ya, pp.350-351, pp.359 l.22-28],
[Ko, pp.300-302] for (2).\qed
\enddemo

\proclaim{Proposition 4.2}
Let $\varphi\in C^{\infty}(\Bbb B(2r))$ satisfies the Monge-Amp\`ere 
equation;
\newline
$\det\left(g_{i\bar{j}}+\frac{i}{2\pi}\varphi_{i\bar{j}}\right)=
e^{F}\det(g_{i\bar{j}})$,
$\varphi_{i\bar{j}}=
\partial^{2}\varphi/\partial z_{i}\partial\bar{z}_{j}$ on the ball of 
radius $2r$ in $\Bbb C^{2}$, and suppose $\lambda(\delta_{ij})\leq
(g_{i\bar{j}}+\frac{i}{2\pi}\varphi_{i\bar{j}})\leq\Lambda(\delta_{ij})$ 
over $\Bbb B(2r)$. Then, there exist constants
$\alpha=\alpha(\lambda,\Lambda,r)>0$ and 
$C=C(\lambda,\Lambda,r,|\partial\bar{\partial}\varphi|_{C^{0}(\Bbb B(2r))},
|F|_{C^{2}(\Bbb B(2r))},|g_{ij}|_{C^{2}(\Bbb B(2r))})\geq 0$ such that
$|\varphi|_{C^{2+\alpha}(\Bbb B(r))}\leq C$.
\endproclaim

\demo{Proof}
See [Si, Chap.2, $\S4$].\qed
\enddemo

\subsubhead
\bf{4.2 Construction of Approximate Ricci-Flat Metrics}\rm 
\endsubsubhead
Let $\pi:\Cal Y\to D$ be the degenerating family of $K3$ surfaces over the 
disc $D$ considered in Proposition 2.1 whose fiber is denoted by $Y_{t}$.
By Proposition 2.1 (5), there exists a coordinate neighborhood 
$(V,o)$ in $\Cal Y$ such that
$$
(V,o)\cong\{(x,y,z,t)\in\Bbb B(2);\,xy-z^{2}-t^{2}=0\},\quad
\pi(x,y,z,t)=t
\tag 4.2
$$
where $\Bbb B(r)$ is the ball of radius $r$ centered at $0$.
Let $L$ be a very ample line bundle over $\Cal Y$ as in Lemma 2.4 and put
$L_{t}:=L|_{Y_{t}}$. Since $L$ is very ample, we may assume that 
$\Cal Y$ is a closed subvariety of $\Bbb P^{N}\times D$ with
$\pi=pr_{2}\circ i$ and $L=\Cal O_{\Bbb P^{N}}(1)|_{\Cal Y}$. 
Let $H$ be a hyperplane of $\Bbb P^{N}$ such that $i(o)\not\in H$. 
Let $\sigma$ be the section of $L$ such that $(\sigma)_{0}=H\cap\Cal Y$. 
Let $h_{L}$ be a Hermitian metric of $L$ with the following properties; 
\newline
{\bf(P2)\rm}
there exist open subsets $W'\Subset W\Subset V$ such that
$h_{L}(\sigma,\sigma)\equiv 1$ on $W'$, and
$c_{1}(L,h_{L}):=-\frac{i}{2\pi}\partial\bar{\partial}
\log h_{L}(\sigma,\sigma)=\omega_{FS}$ on $X\backslash W$ where 
$\omega_{FS}$ is the restriction of the Fubini-Study form of $\Bbb P^{N}$.
\par
Let $\kappa_{Y_{t}}$ be the Ricci-flat K\"ahler metric on $Y_{t}$ 
cohomologous to $c_{1}(L_{t})$. By Kobayashi-Todorov ([K-T]), 
$\kappa_{Y_{0}}$ is a Ricci-flat K\"ahler metric on $Y_{0}$ 
in the sense of orbifold. Let
$\eta_{Y_{t}}\in H^{0}(Y_{t},\Omega^{2})$ be a symplectic form on 
$Y_{t}$ depending holomorphically in $t\in D$. By a suitable choice of 
$\eta_{Y_{t}}$, we may assume the following;
$$
\kappa_{Y_{t}}^{2}=h(t)\,\eta_{Y_{t}}\wedge\bar{\eta}_{Y_{t}},\quad
h(t)=\deg L/\|\eta_{Y_{t}}\|^{2},\quad h(0)=1.
\tag 4.3
$$
Fix an isomorphism 
$j:(\Bbb C^{2}/\pm 1\cap\Bbb B(2),0)\to (Y_{0}\cap V,o)$;
$$
j(z_{1},z_{2})=(z_{1}^{2},z_{2}^{2},z_{1}z_{2})\in Y_{0}.
\tag 4.4
$$
We denote by $(r,\sigma)$ the polar coordinates of $\Bbb C^{2}$; $r=\|z\|$
and $\sigma=z/\|z\|\in S^{3}$. By [K-T, Theorem 1], we get the following.

\proclaim{Lemma 4.1}
There exist $u_{0}\in C^{\infty}(Y_{0})$, $v\in C^{\omega}(\Bbb B(2))$ 
and $c>0$ such that 
\newline
$\kappa_{0}=c_{1}(L_{0},e^{-u_{0}}h_{L})$ and
$j^{*}u_{0}(z)=c(\|z\|^{2}+v(z))$. Here $v$ has the expansion;
\newline
$v(z)=\sum_{k=2}^{\infty}a_{2k}(\sigma)\,r^{2k}=
\sum a_{IJ}z^{I}\bar{z}^{J}$ where 
$a_{2k}(\sigma)=r^{-2k}\sum_{|I|+|J|=2k}a_{IJ}z^{I}\bar{z}^{J}$. 
\endproclaim

For simplicity, we assume $c=1$ in the sequel.( General case is easily
obtained by small modifications of this case.)
Let $D(\delta)$ be the disc of radius $\delta(\ll 1)$ and fix a smooth
trivialization over $D(\delta)$;
$I:(Y_{0}\backslash W)\times D(\delta)\cong\Cal Y\backslash W$ such that
$i_{t}:=I(\cdot,t):Y_{0}\backslash W\cong Y_{t}\backslash W$ is a 
diffeomorphism for any $t$. Set
$$
u_{t}:=(i_{t}^{-1})^{*}u_{0}\in C^{\infty}(Y_{t}\backslash W),\quad 
h_{t}:=e^{-u_{t}}\, h_{L}.
\tag 4.5
$$
Then, $c_{1}(L_{t},h_{t})$ is a (1,1)-form on $Y_{t}\backslash W$ 
approximating $\kappa_{0}$.

\proclaim{Proposition 4.3}
There exist $\epsilon,C>0$ with the following properties.
If $\delta>0$ is chosen small enough, then for any $|t|<\delta$, one has
\flushpar{$(1)$} 
$c_{1}(L_{t},h_{t})\geq\epsilon\,c_{1}(L_{t},h_{L})$ on 
$Y_{t}\backslash W$,
\flushpar{$(2)$} 
$|\star_{0}i_{t}^{*}c_{1}(L_{t},h_{t})^{2}-1|\leq C\,|t|$
on $Y_{0}\backslash W$.
\newline
Here, $\star_{t}$ denotes the Hodge $*$-operator relative to 
$c_{1}(L_{t},h_{L})=\omega_{FS}$.
\endproclaim

\demo{Proof}
By construction, there exist 
$\alpha,\beta,\gamma\in C^{\infty}(D(\delta)\times(Y_{0}\backslash W))$ 
such that
$$
i_{t}^{*}\star_{t}c_{1}(L_{t},h_{t})\wedge c_{1}(L_{t},h_{L})=\alpha,\quad
i_{t}^{*}\star_{t}c_{1}(L_{t},h_{t})^{2}=\beta,\quad
\star_{0}i_{t}^{*}c_{1}(L_{t},h_{t})^{2}=\gamma.
\tag 4.6
$$
By definition, we get $\alpha(z,0)\geq 3\epsilon_{0}$,
$\beta(z,0)\geq9\epsilon_{0}^{2}$ and $\gamma(z,0)=1$. 
As $\alpha,\beta,\gamma$ are continuous in $t$, if $|t|\ll1$,
we get $\alpha(t,z)\geq 2\epsilon_{0}$, $\beta(t,z)\geq\epsilon_{0}^{2}$ 
and $|\gamma(t,z)-1|\leq C\,|t|$, from which the assertion follows.\qed
\enddemo

\remark{Remark}
In the sequel of this section, we denote by $C(>0)$ a constant 
independent of $s\in[0,1]$, $t\in D$, and $x\in Y_{t}$, though its value
may change in each estimate.
\endremark
Our next task is to construct an approximating family of K\"ahler metrics
on $W$. Put $\Cal X:=\{(x,y,z,t)\in\Bbb C^{4};\,xy-z^{2}-t^{2}=0\}$ and
$\pi(x,y,z,t)=t$. The following is due to Kronheimer ([Kr]).

\proclaim{Proposition 4.4}
There exists $q\in C^{0}(\Cal X)\cap C^{\infty}(\Cal X\backslash\{0\})$
such that
\flushpar{$(1)$} 
$j^{*}q_{0}=\partial\bar{\partial}\|z\|^{2}$ where 
$j:\Bbb C^{2}/\pm 1\to\Cal X_{0}$ is the map defined by $(4.4)$,
\flushpar{$(2)$} 
$\partial\bar{\partial}q_{t}$ is a Ricci-flat ALE metric on 
$\Cal X_{t}$ where $q_{t}:=q|_{\Cal X_{t}}$, $\Cal X_{t}:=\pi^{-1}(t)$
\flushpar{$(3)$} 
$q$ is homogeneous of order $1$; $q(sx,sy,sz,st)=|s|\,q(x,y,z,t)$.
\endproclaim

\remark{Remark}
A Riemannian manifold $(M,g)$ is ALE if there exists a compact subset $K$ 
of $M$, a finite group $\Gamma\subset O(m,\Bbb R)$ and a diffeomorphism 
$f:\Bbb R^{m}\backslash\Bbb B(R)/\Gamma\to M\backslash K$ such that
$f^{*}g=\delta+O(r^{-2-k})$ for some $k>0$ where $\delta$ is the 
Euclidean metric.
\endremark

By (4.3), we identify a neighborhood $V$ of $o$ in $\Cal Y$ with a 
neighborhood around $0$ of $\Cal X$. Thus 
$Y_{t}\cap V=\Cal X_{t}\cap\Bbb B(2)$ and we regard $q_{t}$ to be a 
function on $Y_{t}\cap V$. Put $S:=\{x\in\Cal X;\,q(x)=1\}$ and
$S_{t}:=S\cap\Cal X_{t}$ for the level set of $q$, and
$S_{<\delta}:=\bigcup_{|t|<\delta}S_{t}$ for the sublevel set.
We consider $S_{t}$ to be a subset in $\Bbb C^{3}$.
As $\pi:S_{<\delta}\to D(\delta)$ is of maximal rank when $\delta\ll 1$, 
we can construct a trivialization of $S_{<\delta}$ by integrating 
vector fields $\xi,\zeta$ on $S_{<\delta}$ such that 
$\pi_{*}\xi=\partial/\partial u$ and $\pi_{*}\zeta=\partial/\partial v$
where $t=u+i\,v$ is the coordinate of $D(\delta)$.

\proclaim{Lemma 4.2}
There exists a trivialization 
$\Psi:S_{0}\times D(\delta)\to S_{<\delta}$ with the property that 
$\Psi(\cdot,t):S_{0}\to S_{t}$ is a diffeomorphism and that
$\Psi(\cdot,0)$ is the identity map on $S_{0}$. We put 
$\psi_{t}(\cdot):=\Psi(\cdot,t)$ and 
$\phi_{t}=\psi_{t}^{-1}:S_{t}\to S_{0}$ for the inverse map.
\endproclaim
 
With this identification of $S_{t}$ with $S_{0}$, we can construct a
good deformation of the K\"ahler potential by using the polar coordinates
in Lemma 4.1. Namely, if 
$u_{0}=q_{0}+\sum_{k>2}a_{2k}(\sigma)\,q_{0}^{k}$ is the expansion in the
polar coordinates, then
$$
u_{t}=q_{t}+\sum_{k>2}(\phi_{t}^{*}a_{2k})\,q_{t}^{2k}
\tag 4.7
$$
will be a good approximation of $u_{0}$. Let us verify it in the sequel.
\par
Set $Y_{t,a}:=\{x\in Y_{t}\cap V;\,q_{t}(x)\geq a\}$.
We extend $\phi_{t}:Y_{t,\delta^{-1}|t|}\to Y_{0,\delta^{-1}|t|}$ and 
$\psi_{t}:Y_{0,\delta^{-1}|t|}\to Y_{t,\delta^{-1}|t|}$ by using the
$\Bbb R_{+}$-action. Namely, we define
$$
\phi_{t}(x):=q_{t}(x)\cdot\phi_{q_{t}(x)^{-1}t}(q_{t}(x)^{-1}x),\quad
\psi_{t}(x):=q_{0}(x)\cdot\psi_{q_{0}(x)^{-1}t}(q_{0}(x)^{-1}x)
\tag 4.8
$$
where $\lambda\cdot x=(\lambda x,\lambda y,\lambda z)$ if
$x=(x,y,z)\in\Bbb C^{3}$. Put $\mu_{\lambda}(x):=\lambda\cdot x$. 
Let $J_{t}$ be the complex structure of $Y_{t}$.
Let $\bar{\partial}_{Y_{t}}$ be the $\bar{\partial}$-operator of $Y_{t}$. 
Via $\psi_{t}$, identify $\bar{\partial}_{Y_{t}}$ with the 
$\bar{\partial}$-operator of 
$(Y_{0,\delta^{-1}|t|}\cap V,\psi_{t}^{*}J_{t})$ and similarly 
$\partial_{Y_{t}}$ with the $\partial$-operator. We denoted them by 
$\bar{\partial}_{t}$ and $\partial_{t}$ respectively;
$$
\partial_{t}:=\psi_{t}^{*}\circ\partial_{Y_{t}}\circ\phi_{t}^{*},\quad
\bar{\partial}_{t}:=
\psi_{t}^{*}\circ\bar{\partial}_{Y_{t}}\circ\phi_{t}^{*},\quad
\partial_{t}\bar{\partial}_{t}=\psi_{t}^{*}\circ\partial_{Y_{t}}
\bar{\partial}_{Y_{t}}\circ\phi_{t}^{*}.
\tag 4.9
$$
Take $f\in C^{\infty}(Y_{0}\cap V)$ and identify it with an even function 
on $\Bbb C^{2}\cap\Bbb B(2)$ via $j$. Let $(z_{1},z_{2})$ be the complex
coordinates of $\Bbb C^2$, and $z_{1}=x_{1}+i\,x_{2}$,
$z_{2}=x_{3}+i\,x_{4}$ be the real coordinates. 
For $z\in Y_{0,\delta^{-1}|t|}\cap V$, we can write
$$
\partial_{t}\bar{\partial}_{t}f(z)
=\sum a_{ij}^{kl}(t,z)\,\partial_{kl}f(z)\,dx_{i}\wedge dx_{j}
+\sum b_{ij}^{k}(t,z)\,\partial_{k}f(z)\,dx_{i}\wedge dx_{j}
\tag 4.10
$$
where $a_{ij}^{kl}$ and $b^{k}_{ij}$ restricted to $D(\delta)\times S_0$ 
are $C^{\infty}$-functions, and $\partial_{i}f=\partial f/\partial x_{i}$ 
etc. 

\proclaim{Lemma 4.3}
Putting $\|z\|=r$, one has the following in the polar coordinates:
$$
\partial_{t}\bar{\partial}_{t}f(z)=
\sum\left\{a_{ij}^{kl}(r^{-2}t,\sigma)\,\partial_{kl}f(z)+
r^{-1}b_{ij}^{k}(r^{-2}t,\sigma)\,\partial_{k}f(z)\right\}
dx_{i}\wedge dx_{j}.
$$
\endproclaim

\demo{Proof}
As $\mu_{\lambda}j(z)=j(\sqrt{\lambda}z)$ by (4.5), the action of 
$\Bbb R_{+}$ on $Y_{0}\cap V$ is expressed by 
$\mu_{\lambda}(z)=(\sqrt{\lambda}z_{1},\sqrt{\lambda}z_{2})$.
Since $\phi_{\lambda t}\circ\mu_{\lambda}=\mu_{\lambda}\circ\phi_{t}$
on $\Bbb R_{+}\times Y_{t,\delta^{-1}|t|}$ by definition, we get
$\partial_{\lambda t}\bar{\partial}_{\lambda t}=
\mu_{\lambda^{-1}}^{*}\circ\partial_{t}\bar{\partial}_{t}\circ
\mu_{\lambda}^{*}$ which yields
$a_{ij}^{kl}(\lambda^{-1}t,\lambda^{-1/2}z)=a_{ij}^{kl}(t,z)$ and
$\lambda^{-1/2}b_{ij}^{k}(\lambda^{-1}t,\lambda^{-1/2}z)=
b_{ij}^{k}(t,z)$. Putting $\lambda=r$, we get the assertion.\qed
\enddemo

Let $\omega_{Y_{t}}$ be the Kronheimer's Ricci-flat ALE metric on 
$Y_{t}\cap V$; 
$$
\omega_{Y_{t}}:=\partial_{Y_{t}}\bar{\partial}_{Y_{t}}q_{t}.
\tag 4.11
$$
As $\psi_{t}^{*}q_{t}=q_{0}$, we find that 
$\omega_{t}:=\partial_{t}\bar{\partial}_{t}q_{0}$ is the K\"ahler form on 
$(Y_{0,\delta^{-1}|t|}\cap V,\psi_{t}^{*}J_{t})$.
Taking $f(z)=\|z\|^{2}$ in (4.10), it follows that, for any 
$z\in Y_{0,\delta^{-1}|t|}\cap V$,
$$
\omega_{t}(z)
=\sum g_{ij}(r^{-2}t,\sigma)\,dx_{i}\wedge dx_{j}
\tag 4.12
$$
where $g_{ij}$ is a smooth functions on $D(\delta)\times S_{0}$. Let $\rho$ 
be a cut-off function such that $\rho(s)\equiv0$ for $s\leq2\delta^{-1}$, 
$\rho(s)\equiv1$ for $s\geq4\delta^{-1}$,
$0\leq\rho'(s)\leq C_{0}\delta$ and $|\rho''(s)|\leq C_{0}\delta^{2}$. 
Let $v\in C^{\omega}(\Bbb B(1))$ be the error term of $u_{0}$ 
appeared in Lemma 4.1. Set
$$
v_{t}(x):=\rho_{t}(x)\,\phi_{t}^{*}v(x),\quad
\rho_{t}(x):=\rho(|t|^{-1}q_{t}(x)),\quad
\tilde{v}_{t}:=\rho(|t|^{-1}r^{2})\,v=\psi_{t}^{*}v_{t}.
\tag 4.13
$$
Then, $v_{t},\rho_{t}\in C^{\infty}(Y_{t})$ and 
$\tilde{v}_{t}\in C^{\infty}(Y_{0})$. By (4.8), we get
$\psi_{t}^{*}\partial_{Y_{t}}\bar{\partial}_{Y_{t}}\,v_{t}=
\partial_{t}\bar{\partial}_{t}\tilde{v}_{t}$.

\proclaim{Lemma 4.4}
Let $*_{t}$ be the Hodge $*$-operator relative to $\omega_{t}$.
There exist $B_{i}\in C^{\infty}(D(\delta)\times [0,1]\times S_{0})$ 
$(i=1,2)$ such that, for any $z\in Y_{0,\delta^{-1}|t|}\cap V$, 
$$
*_{t}\left(\partial_{t}\bar{\partial}_{t}\tilde{v}_{t}
\wedge\omega_{t}\right)(z)
=r^{2}\,B_{1}(r^{-2}|t|,r^{2},\sigma),\quad
*_{t}(\partial_{t}\bar{\partial}_{t}\tilde{v_{t}})^{2}(z)
=r^{4}\,B_{2}(r^{-2}|t|,r^{2},\sigma).
$$
\endproclaim

\demo{Proof}
As $v\in C^{\omega}([0,1]\times S^{3})$ by Lemma 4.1,
it follows from (4.12) and
$$
\align
r^{-1}\partial_{i}\{\rho(|t|^{-1}r^{2})\,v\}
&=2|t|^{-1}r^{-1}x_{i}\rho'(|t|^{-1}r^{2})\,v
+r^{-1}\rho(|t|^{-1}r^{2})\partial_{i}v,\\
\partial_{ij}\{\rho(|t|^{-1}r^{2})\,v\}
&=4|t|^{-2}\rho''(|t|^{-1}r^{2})\,x_{i}x_{j}\,v+
2|t|^{-1}r^{2}\rho'(|t|^{-1}r^{2})\delta_{ij}\,v\\
\,&\quad+2|t|^{-1}\rho'(|t|^{-1}r^{2})
(x_{j}\,\partial_{i}v+x_{i}\,\partial_{j}v)+
\rho(|t|^{-1}r^{2})\,\partial_{ij}v
\tag 4.14
\endalign
$$
that there exists 
$A_{i}\in C^{\infty}(D(\delta)\times [0,1]\times S_{0})$ such that
$\omega_{t}^{2}(z)=A_{0}(r^{-2}|t|,r^{2},\sigma)\,dV$,
$\partial_{t}\bar{\partial}_{t}\tilde{v}_{t}\wedge\omega_{t}(z)=
r^{2}\,A_{1}(r^{-2}|t|,r^{2},\sigma)\,dV$ and
$(\partial_{t}\bar{\partial}_{t}\tilde{v}_{t})^{2}(z)
=r^{4}\,A_{2}(r^{-2}|t|,r^{2},\sigma)\,dV$ where $dV$ is the volume form 
of $(Y_{0},\omega_{0})$. Since $\omega_{0}^{2}=dV$, we get 
$A_{0}(0,r,\sigma)\equiv 1$ for any $(r,\sigma)\in [0,1]\times S_{0}$. 
By the compactness of $[0,1]\times S_{0}$, there exists $C_{0}$ such that
$0<C_{0}^{-1}\leq A_{0}(s,r,\sigma)\leq C_{0}<\infty$ 
for $|s|<\delta$ and $(r,\sigma)\in[0,1]\times S_{0}$.
Since $*_{t}F=F/\omega_{t}^{2}$ for a $4$-form $F$, we get
$B_{1}=A_{1}/A_{0}$ and $B_{2}=A_{2}/A_{0}$.\qed
\enddemo

We define $\Omega_{Y_{t}}$ and $G_{Y_{t}}$ as follows;
$$
\Omega_{Y_{t}}:=
\omega_{Y_{t}}+\partial_{Y_{t}}\bar{\partial}_{Y_{t}}\,v_{t},
\quad G_{Y_{t}}:=\Omega_{Y_{t}}^{2}/\omega_{Y_{t}}^{2}.
\tag 4.15
$$
As before, set $\Omega_{t}:=\psi_{t}^{*}\Omega_{Y_{t}}=
\omega_{t}+\partial_{t}\bar{\partial}_{t}\tilde{v}_{t}$ and
$G_{t}:=\psi_{t}^{*}G_{Y_{t}}$. 

\proclaim{Proposition 4.5}
If $|t|\ll1$, $\Omega_{Y_{t}}$ becomes a K\"ahler metric on
$Y_{t}\cap V$.
\endproclaim

\demo{Proof}
When $q_{t}(x)\leq\delta^{-1}|t|$, $v_{t}\equiv 0$ by (4.13) and 
$\Omega_{Y_{t}}>0$ because $\Omega_{Y_{t}}=\omega_{Y_{t}}$ in this case.
Thus, it is enough to show that
$*_{t}\Omega_{t}\wedge\omega_{t}(z)>0$ and
$*_{t}\Omega_{t}^{2}(z)>0$ for any $z\in Y_{0,\delta^{-1}|t|}\cap V$
if $|t|\ll1$. By Lemma 4.4, we get
$$
*_{t}\Omega_{t}\wedge\omega_{t}(z)=
1+r^{2}\,B_{1}(|t|/r^{2},r^{2},\sigma),
\,\,
*_{t}\Omega_{t}^{2}=1+2r\,B_{1}(|t|/r^{2},r^{2},\sigma)+
r^{2}\,B_{2}(|t|/r^{2},r^{2},\sigma).
\tag 4.16
$$
As $1+r^{2}\,B_{1}(0,r^{2},\sigma)\geq C>0$
for any $z\in\Bbb B(1)=Y_{0}\cap V$ because $\Omega_{0}$ is a K\"ahler 
metric, choosing $t$ small enough, the right hand side of the first 
formula of (4.16) is greater than $C/2$. Similarly, the rest inequality 
can be shown.
\qed
\enddemo

Let $\tau_{Y_{t}}$ be the holomorphic family of symplectic forms on 
$Y_{t}\cap V$ such that
\newline
$\tau_{Y_{t}}\wedge\bar{\tau}_{Y_{t}}=\omega_{Y_{t}}^{2}$.
Under the identification of $Y_{t}\cap V$ with $\Cal X_{t}\cap\Bbb B(2)$, 
we get $\omega_{Y_{t}}=g_{Y_{t}}(I_{t}\cdot,\cdot)$ and
$\tau_{Y_{t}}=g_{Y_{t}}(J_{t}\cdot,\cdot)+
\sqrt{-1}g_{Y_{t}}(K_{t}\cdot,\cdot)$ where
$g_{Y_{t}}$ is the underlying Riemannian metric and $(I_{t},J_{t},K_{t})$ 
are complex structures defining the hyper-K\"ahler structure. Since $(V,o)$ 
is normal, there exists a holomorphic function $f_{V}\in\Cal O(V)$ 
such that $\eta_{Y_{t}}|_{V}=f_{Y_{t}}\cdot\tau_{Y_{t}}$ and
$f_{Y_{t}}:=f_{V}|_{Y_{t}}$. As $f_{Y_{t}}$ has no zero for any 
$t\in B(\delta)$, $f_{V}$ has no zero on $V$. Thus 
there exists the lower bound on $V$; $|f_{V}(x)|\geq C>0$. By (4.3), (4.15) 
and the definition of $f_{Y_{t}}$, we get
$$
|f_{Y_{0}}|^{2}=(\eta_{Y_{0}}\wedge\bar{\eta}_{Y_{0}})/
(\tau_{Y_{0}}\wedge\bar{\tau}_{Y_{0}})=\Omega_{0}^{2}/\omega_{0}^{2}=G_{0}.
\tag 4.17
$$
Let $u_{t}$ be the function on $Y_{t}\backslash W$ as in (4.6) and $h_{L}$
be the Hermitian metric of $L$ with the property {\bf(P2)\rm} as before.
Let $\chi\geq0$ be a cut-off function on $\Cal Y$ such that 
$\chi(x)\equiv 0$ for $x\in W$ and $\chi(x)\equiv 1$ for 
$x\in\Cal Y\backslash V$. We shall use the following (1,1)-form 
$\tilde{\kappa}_{Y_{t}}$ as a background metric in approximating 
$\kappa_{Y_{t}}$; 
$$
\tilde{\kappa}_{Y_{t}}:=c_{1}(L_{t},e^{-\theta_{t}}h_{L})=
\frac{i}{2\pi}\partial_{Y_{t}}\bar{\partial}_{Y_{t}}\{\theta_{t}-
\log h_{L}(\sigma,\sigma)\},\quad
\theta_{t}:=\chi\,u_{t}+(1-\chi)(q_{t}+v_{t}).
\tag 4.18
$$

\proclaim{Proposition 4.6}
If $|t|\ll1$, $\tilde{\kappa}_{Y_{t}}$ is a K\"ahler metric with the 
following estimate;
$$
\left|\frac{\eta_{Y_{t}}\wedge\bar{\eta}_{Y_{t}}}
{\tilde{\kappa}_{Y_{t}}^{2}}-1\right|_{\infty}\leq C\,|t|.
$$
\endproclaim

\demo{Proof}
(1) Suppose $x=\psi_{t}(z)\in Y_{t,\delta^{-1}|t|}\cap V$.
By definition, one has
$$
\eta_{Y_{t}}\wedge\bar{\eta}_{Y_{t}}/\Omega_{Y_{t}}^{2}=
(\eta_{Y_{t}}\wedge\bar{\eta}_{Y_{t}}/\omega_{Y_{t}}^{2})\cdot
(\Omega_{Y_{t}}^{2}/\omega_{Y_{t}}^{2})^{-1}=
|f_{Y_{t}}|^{2}G_{Y_{t}}^{-1}.
\tag 4.19
$$
As $G_{t}(z)=1+r^{2}B(\frac{|t|}{r^{2}},r^{2},\sigma)$ 
where $B(s,r^{2},\sigma):=B_{1}(s,r^{2},\sigma)+r^{2}B_{2}(s,r^{2},\sigma)$,
$$
|G_{t}(z)-G_{0}(z)|=r^{2}|B(s,r^{-2}|t|,\sigma)-B(s,0,\sigma)|\leq 
|dB|_{\infty}\,|t|\leq C\,|t|.
\tag 4.20
$$ 
Let $\|\cdot\|$ be the Euclidean norm of $\Bbb C^{3}$. Then, we get
$$
||f_{V}(\psi_{t}(z))|^{2}-|f_{Y_{0}}(z)|^{2}|
\leq|f_{V}|_{\infty}\,|df_{V}|_{\infty}\,\|\psi_{t}(j(z))-j(z)\|
\leq C\,|t|
\tag 4.21
$$ 
which, together with (4.19-20), yields that 
$$
\aligned
\left|\psi_{t}^{*}\frac{\eta_{Y_{t}}\wedge\bar{\eta}_{Y_{t}}}
{\Omega_{Y_{t}}^{2}}(z)-1\right|
=\left|\frac{\psi_{t}^{*}|f_{V}|^{2}-G_{t}}{G_{t}}\right|
&\leq\frac{\left|\psi_{t}^{*}|f_{V}|^{2}-|f_{Y_{0}}|^{2}\right|+
|G_{t}-G_{0}|}{|G_{t}|}\\
&\leq\frac{C\,|t|}{|G_{0}|-C\,|t|}.
\endaligned
\tag 4.22
$$
Since both $\Omega_{Y_{0}}$ and $\omega_{Y_{0}}$ are K\"ahler metrics on
$Y_{0}\cap V$, there exists $C>0$ such that 
$|G_{0}|-C\,|t|\geq C$ on $Y_{0}\cap V$, which, 
together with (4.22), yields the assertion. 
\par{(2)}
If $q_{t}(x)\leq\delta^{-1}|t|$, it follows from (4.13) and (4.15) that 
$\Omega_{Y_{t}}(x)=\omega_{Y_{t}}(x)$ and
$$
\eta_{Y_{t}}\wedge\bar{\eta}_{Y_{t}}(x)/\Omega_{Y_{t}}^{2}(x)=
\eta_{Y_{t}}\wedge\bar{\eta}_{Y_{t}}(x)/\omega_{Y_{t}}^{2}(x)=
|f_{Y_{t}}(x)|^{2}.
\tag 4.23
$$
By (4.17), we get $f_{Y_{t}}(0)=1$ because $\Omega_{0}(0)=\omega_{0}(0)$ by
construction which, together with (4.23), yields
$$
\left|\frac{\eta_{Y_{t}}\wedge\bar{\eta}_{Y_{t}}}{\Omega_{Y_{t}}^{2}}(x)-1
\right|=\left||f_{V}(x)|^{2}-|f_{V}(0)|^{2}\right|\leq
|f_{V}|_{\infty}\,|df_{V}|_{\infty}\|x\|.
\tag 4.24
$$
By Proposition 4.4 (3), one has
$C^{-1}(\|x\|+|t|)\leq q_{t}(x)\leq C(\|x\|+|t|)$ 
which, together with $q_{t}(x)\leq\delta^{-1}|t|$, implies 
$\|x\|\leq C\,|t|$ and yields the assertion.
\par{(3)}
Consider $Y_{t}\backslash W$.
By construction (Lemma 4.1 and (4.15)), we get $u_{0}=q_{0}+v_{0}$ on 
$V$. Therefore, $\tilde{\kappa}_{Y_{0}}=\kappa_{Y_{0}}$ and there is a 
smooth function $w$ on $V\backslash W$ such that 
$u_{t}-(q_{t}+v_{t})|_{V\backslash W}=t\,w|_{Y_{t}\cap(V\backslash W)}+
\bar{t}\,\bar{w}|_{Y_{t}\cap(V\backslash W)}$. Thus
there exist (1,1)-forms $\xi$ and $\xi'$ on $V\backslash W$ such that
$\tilde{\kappa}_{Y_{t}}|_{V\backslash W}=c_{1}(L_{t},h_{t})+t\,\xi+
\bar{t}\,\xi'$
which, together with Proposition 4.4, implies that $\tilde{\kappa}_{Y_{t}}$ 
is positive on $Y_{t}\cap(V\backslash W)$. As $\tilde{\kappa}_{Y_{t}}>0$ 
on $(Y_{t}\backslash V)\cup(Y_{t}\cap W)$ by Propositions 4.4 and 4.5, 
it becomes a K\"ahler metric on $Y_{t}$. Now, the desired estimate 
follows from Propositions 4.3 and 4.6.\qed
\enddemo

If $|t|\ll1$, the following extension of $\omega_{Y_{t}}$ becomes a 
K\"ahler metric on $Y_{t}$;
$$
\omega_{Y_{t}}:=c_{1}(L_{t},e^{-w_{t}}h_{L})=
\frac{i}{2\pi}\partial_{Y_{t}}\bar{\partial}_{Y_{t}}\{w_{t}-
\log h_{L}(\sigma,\sigma)\},\quad
w_{t}:=\chi\,u_{t}+(1-\chi)\,q_{t}.
\tag 4.25
$$
By construction, $\omega_{Y_{t}}|_{Y_{t}\cap W}=
\partial_{Y_{t}}\bar{\partial}_{Y_{t}}q_{t}$. Let $\tilde{\kappa}_{t}(s)$ 
$(0\leq s\leq 1)$ be the homotopy of K\"ahler metrics on $Y_{t}$ 
joining $\tilde{\kappa}_{Y_{t}}$ and $\omega_{Y_{t}}$ defined by
$$
\tilde{\kappa}_{t}(s):=s\,\tilde{\kappa}_{Y_{t}}+
(1-s)\,\omega_{Y_{t}}. 
\tag 4.26
$$

\proclaim{Proposition 4.7}
Let $\Ric(\tilde{\kappa}_{t}(s))=
\partial_{Y_{t}}\bar{\partial}_{Y_{t}}\log
(\tilde{\kappa}_{t}(s)^{2}/\eta_{Y_{t}}\wedge\bar{\eta}_{Y_{t}})$ be the
Ricci curvature of $(Y_{t},\tilde{\kappa}_{t}(s))$. Then, one has
$|\Ric(\tilde{\kappa}_{t}(s))|_{\infty}\leq C$.
\endproclaim

\demo{Proof}
As the assertion is clear for $x\in Y_{t}\backslash W$, we suppose that 
$x\in Y_{t}\cap W$. When $q_{t}(x)\leq\delta^{-1}|t|$, 
$\tilde{\kappa}_{Y_{t}}(s)=\omega_{Y_{t}}$ by construction, which is 
followed by the assertion as $\omega_{Y_{t}}$ is Ricci-flat. 
When $x=\psi_{t}(r,\sigma)\in Y_{t,\delta^{-1}}\cap W$, 
it follows from (4.19) that
$$
\psi_{t}^{*}\Ric(\tilde{\kappa}_{t}(s))=\partial_{t}\bar{\partial}_{t}
\log\{1+r^{2}B'(s,r^{-2}|t|,r^{2},\sigma)\}
\tag 4.27
$$
because $\partial_{Y_{t}}\bar{\partial}_{Y_{t}}\log|f_{Y_{t}}|^{2}=0$
where $B'(s,u,r^{2},\sigma):=2s\,B_{1}(u,r^{2},\sigma)
+s^{2}\,r^{2}B_{2}(u,r^{2},\sigma)$. As
$|r^{-1}\partial_{i}(r^{2}B')|_{\infty}+
|\partial_{ij}(r^{2}B')|_{\infty}\leq C$, we get the assertion.\qed
\enddemo

Let $\tilde{R}_{t}(s)$ be the curvature of 
$(Y_{t},\tilde{\kappa}_{t}(s))$.

\proclaim{Proposition 4.8}
For any $(t,x)\in D(\delta)\times Y_t\cap W$, one has
$|\tilde{R}_{t}(s)(x)|\leq C\,q_{t}(x)^{-1}$.
\endproclaim

\demo{Proof}
By (4.10-11) and (4.15), there exists
$g_{ij}\in C^{\infty}(D(\delta)\times[0,1]\times S_{0}\times[0,1])$ 
such that $\psi_{t}^{*}\tilde{\omega}_{t}(s)(z)=
\sum_{ij}g_{ij}(r^{-1}|t|,r^{2},\sigma,s)\,dx_{i}dx_{j}$
and $(g_{ij})$ is an uniformly positive definite matrix. As
$|\partial^{k}g_{ij}(r^{-1}|t|,r^{2},\sigma,s)|\leq C_{k}\,r^{-k}$,
we get the assertion.\qed
\enddemo

\proclaim{Theorem 4.1}
For any $(t,s)\in D(\delta)\times[0,1]$, one has
$C^{-1}\,\tilde{\kappa}_{t}(s)\leq
\kappa_{Y_{t}}\leq C\,\tilde{\kappa}_{t}(s)$.
\endproclaim

\demo{Proof}
Since $\tilde{\kappa}_{t}(s)
=\omega_{Y_{t}}+s\partial_{Y_{t}}\bar{\partial}_{Y_{t}}v_{t}$ on 
$Y_{t}\cap V$, one has
$C^{-1}\,\omega_{Y_{t}}\leq\tilde{\kappa}_{t}(s)\leq C\,\omega_{Y_{t}}$
by (4.16). Thus, it is enough to show the inequality when $s=0$. 
Consider the Monge-Amp\`ere equation (4.1) with $X=Y_{t}$, 
$\kappa=\tilde{\kappa}_{Y_{t}}$, $\phi=\phi_{t}$ and
$F=F_{t}=
\log(h(t)\eta_{Y_{t}}\wedge\bar{\eta}_{Y_{t}}/\tilde{\kappa}_{Y_{t}}^{2})$.
Here $h(t)$ is defined in (4.3). Since $h(t)=1+O(|t|)$ as $t\to 0$, 
it follows from Proposition 4.7 that $|F_{t}|_{\infty}\leq C\,|t|$. 
Thus, by Propositions 4.7-8, all of 
$|\Delta_{t}F_{t}|_{\infty}/|R_{t}|_{\infty}$, $|F_{t}|_{\infty}$,
$|R_{t}|_{\infty}|F_{t}|_{\infty}$, $V_{t}$ are uniformly bounded in
$t\in D(\delta)$ where the subscript $t$ means that 
these quantities are considered relative to $\tilde{\kappa}_{Y_{t}}$. 
As $\tilde{\kappa}_{Y_{t}}$ is quasi-isometric to $\omega_{Y_{t}}$ on
$Y_{t}\cap V$, the Sobolev constant $S_{t}$ is uniformly bounded. 
As $\lambda_{t}$, the first eigenvalue of the Laplacian of 
$(Y_{t},\tilde{\kappa}_{Y_{t}})$ is continuous in $t$ 
([Yo1, Theorem 5.1]), we also get an uniform bound $\lambda_{t}^{-1}<C$ 
because $\lambda_{0}>0$. The assertion follows from Proposition 4.1 (2) 
together with these bounds.\qed
\enddemo

Consider the following Monge-Amp\`ere equation on $Y_{t}$;
$$
\aligned
\,&
\left(\tilde{\kappa}_{Y_{t}}+\frac{i}{2\pi}\partial\bar{\partial}
\varphi_{t}(s)\right)^{2}
=e^{sF_{t}+a_{t}(s)}\,\tilde{\kappa}_{Y_{t}}^{2},\quad
\int_{Y_{t}}\varphi_{t}(s)\,\tilde{\kappa}_{Y_{t}}^{2}=0
\endaligned
\tag 4.28
$$
where $F_{t}=\log(h(t)\,\eta_{Y_{t}}\wedge\bar{\eta}_{Y_{t}}/
\tilde{\kappa}_{Y_{t}}^{2})$ and 
$a_{t}(s)=\log(\deg L/\int_{Y_{t}}e^{sF_{t}}\,\tilde{\kappa}_{Y_{t}}^{2})$.
Set $\kappa_{t}(s)$ for the smooth homotopy of K\"ahler metrics joining
$\tilde{\kappa}_{Y_{t}}$ and $\kappa_{Y_{t}}$;
$$
\kappa_{t}(s):=
\tilde{\kappa}_{Y_{t}}+\frac{i}{2\pi}\partial\bar{\partial}\varphi_{t}(s)>0,
\quad
\kappa_{t}(0)=\tilde{\kappa}_{Y_{t}},\quad\kappa_{t}(1)=\kappa_{Y_{t}}.
\tag 4.29
$$

\proclaim{Theorem 4.2}
For any $t\in D$, $s\in[0,1]$ and $x\in Y_{t}$, one has
$$
|F_{t}|_{\infty}+|a_{t}(s)|\leq C\,|t|,\quad
C^{-1}\,\tilde{\kappa}_{Y_{t}}(x)\leq\kappa_{t}(s)(x)\leq 
C\,\tilde{\kappa}_{Y_{t}}(x), \quad 
|\Ric(\kappa_{t}(s))|_{\infty}\leq C.
$$
\endproclaim

\demo{Proof}
The first inequality follows from Proposition 4.8, the third one from 
Proposition 4.7 and (4.28). The second one is similarly proved as 
Theorem 4.1.\qed
\enddemo

Let $R_{t}(s)$ be the curvature of
$(Y_{t},\kappa_{t}(s))$ and $\nabla_{t}(s)$ its covariant derivative.

\proclaim{Theorem 4.3}
On $Y_{t}\cap W$, one has
$|\nabla_{t}(s)^{k}\,R_{t}(s)(x)|\leq C\,q_{t}(x)^{-1-\frac{k}{2}}$.
\endproclaim

\demo{Proof} 
For $p\in S_{s}$, let $B(p,r)$ be the metric ball of radius $r$ 
centered at $p$ relative to the metric $\omega_{Y_{s}}$, 
and $w=(w_{1},w_{2})$ its holomorphic normal coordinates. 
Fix $\epsilon>0$ small enough so that $(B(p,\epsilon),w)$ becomes 
a coordinate neighborhood for any $p\in S_{s}$ ($s\in D(\delta)$).   
Take $x_{0}\in Y_{t,\delta^{-1}|t|}\cap W$. Put 
$y_{0}:=q_{t}(x_{0})^{-1}x_{0}\in S_{q_{t}(x_{0})^{-1}t}$ and
$B=B(y_{0},\epsilon)$. For $y\in B$, put $x:=q_{t}(x_{0})y\in Y_{t}$
which makes $(B,w)$ to be a coordinate neighborhood of $Y_{t}$ at 
$x_{0}$. As $y=y_{0}+w$, $v(z)=O(\|z\|^{4})$ and 
$q_{t}(x_{0})^{-1}(q_{t}(x)+v_{t}(x))=
q_{q_{t}(x_{0})^{-1}t}(y)+
q_{t}(x_{0})^{-1}v(q_{t}(x_{0})\phi_{q_{t}(x_{0})^{-1}t}(y))$, we get
$|q_{t}(x_{0})^{-1}(q_{t}+v_{t})|_{C^{k}(B)}\leq C(k,|v|_{C^k})$.
Thus, by Propositions 4.5 and 4.8, if we write
$q_{t}(x_{0})^{-1}\tilde{\kappa}_{Y_{t}}=
\sum_{i,j}\tilde{g}_{i\bar{j}}(t,w)dw_{i}\wedge d\bar{w}_{j}$ on $B$,
we get $\lambda\,I\leq(\tilde{g}_{i\bar{j}})\leq\Lambda\,I$ and
$|\tilde{g}_{i\bar{j}}|_{C^{k}(B)}\leq C_{k}$ where 
$\lambda,\Lambda,C_{k}>0$ are independent of $x_{0}$ and $t$.
Putting 
$\bar{\kappa}_{t}:=q_{t}(x_{0})^{-1}\,\tilde{\kappa}_{Y_{t}}$ and
$\bar{\varphi}_{t}(s):=q_{t}(x_{0})^{-1}\,\varphi_{t}(s)$, 
the Monge-Amp\`ere equation (4.28) becomes as follows on $B$;
$$
\det\left(\tilde{g}_{i\bar{j}}+\frac{i}{2\pi}\partial_{i\bar{j}}
\bar{\varphi}_{t}(s)\right)^{2}=
e^{sF_{t}+a_{t}(s)}\,\det(\tilde{g}_{i\bar{j}}),
\quad\partial_{i\bar{j}}\bar{\varphi}_{t}(s)=
\frac{\partial^{2}\bar{\varphi}_{t}(s)}
{\partial w_{i}\bar{\partial}w_{j}}.
\tag 4.30
$$
We find $|\partial_{i\bar{j}}\bar{\varphi}_{t}(s)|_{C^{0}(B)}\leq C$ and
$|s\,F_{t}+a_{t}(s)|_{C^{2}(B)}\leq C$ by Theorem 4.2 and Proposition 4.7, 
which, together with Proposition 4.2, yields that there exists
$\alpha>0$ independent of $t$ and $x$ such that
$|\bar{\varphi}_{t}(s)|_{C^{2+\alpha}(B')}\leq C$ where 
$B'=B(y_{0},\epsilon/2)$. Thus, 
$\tilde{g}_{i\bar{j}}+\frac{i}{2\pi}\partial_{i\bar{j}}\bar{\varphi}_t(s)$ 
is uniformly elliptic and H\"older continuous on $B'$.
From the same argument as in [Si, Chap.II, (1.4)], it follows that 
arbitrary $C^{k}$-norm of $\bar{\varphi}_{t}(s)$ on $B'$ is 
uniformly bounded and thus so is the covariant derivative of curvature of 
$q_{t}(x_{0})^{-1}\,\kappa_{t}(s)$. Rewriting this in terms of the metric
$\kappa_{t}(s)$, we get the assertion. 
When $q_{t}(x_{0})\leq\delta^{-1}|t|$, as
$|t|^{-1}\mu_{|t|}^{*}\tilde{\kappa}_{\Cal X_{t}}=\omega_{\Cal X_{1}}$
by the map $\mu_{|t|}(x)=|t|\cdot x$, we also get the assertion by the same
argument as above because the injectivity radius of 
$(\Cal X_{1},\omega_{\Cal X_{1}})$ is positive.\qed 
\enddemo

\beginsection
$\S5$. Reduction to the ALE Instanton

\par
\subsubhead
\bf{5.1 Estimates of Anomaly}\rm
\endsubsubhead
Let us consider the situation in $\S 4.2$ and keep notations there. 
In this section, we assume that $\pi:\Cal Y\to D$ is a degenerating family 
of 2-elementary $K3$ surfaces with involution $\iota$ and 
$L$ is $\iota$-invariant. Thus, $\kappa_{Y_{t}}$ is also $\iota$-invariant 
by the uniqueness of Ricci-flat K\"ahler metric. As $o$ is a fixed point 
of $\iota$, according to the type in (2.11), we fix the local coordinates 
$(x,y,z,t)$ as in Proposition 2.3. Relative to this local coordinates, 
we can construct $\tilde{\kappa}_{Y_{t}}$ as in (4.18) and 
$\omega_{Y_{t}}$ as in (4.25). By taking the average of the action of 
$\iota$ if necessary, we may assume that both $\tilde{\kappa}_{Y_{t}}$ 
and $\omega_{Y_{t}}$ are $\iota$-invariant, and so is 
$\tilde{\kappa}_{t}(s)$.
\par
Let $\widetilde{\Td}(Y_{t};\kappa_{Y_{t}},\tilde{\kappa}_{Y_{t}})$ be the 
Bott-Chern secondary class associated to the Todd genus and 
$\kappa_{Y_{t}}$, $\tilde{\kappa}_{Y_{t}}$. By definition 
([B-C], [B-G-S, I (e)]), we get
$$
\aligned
\widetilde{\Td}(Y_{t};\kappa_{Y_{t}},\tilde{\kappa}_{Y_{t}})^{(2,2)}
&=\frac{1}{24}\int_{1}^{0}\left(\Tr\,\kappa_{t}(s)^{-1}\dot{\kappa}_{t}(s)
\right)\,c_{2}(R_{t}(s))\,ds\\
&\quad+\frac{1}{24}\int_{1}^{0}c_{1}(R_{t}(s))\cdot
\left.\frac{d}{d\epsilon}\right|_{\epsilon=0}c_{2}\left(R_{t}(s)+
\epsilon\,\frac{\dot{\kappa}_{t}(s)}{\kappa_{t}(s)}\right)\,ds
\endaligned
\tag 5.1
$$
where $\dot{\kappa}_{t}(s)=\partial_{s}\kappa_{t}(s)$. 
As in $\S4.2$, we denote by $C$ a constant independent of $s\in[0,1]$,
$t\in D$ and $x\in Y_{t}$, though its value may change in each estimate.
For a norm $|\cdot|$, $|\cdot|_{\kappa}$ means that it uses the metric 
$\kappa$ in measuring.

\proclaim{Lemma 5.1}
For any $t\in D\backslash\{0\}$, one has
$$
\left|\int_{Y_{t}}\widetilde{\Td}(Y_{t};\kappa_{Y_{t}},
\tilde{\kappa}_{X_{t}})^{(2,2)}\right|\leq
C\,\sup_{s\in[0,1]}
\left|\kappa_{t}(s)^{-1}\dot{\kappa}_{t}(s)
\right|_{\infty,\kappa_{X_{t}}}.
$$
\endproclaim

\demo{Proof}
In the proof, every norms and volumes are those relative to 
$\kappa_{t}(s)$. By (5.1),
$$
\left|\int_{Y_{t}}\widetilde{\Td}(Y_{t};\kappa_{Y_{t}},
\tilde{\kappa}_{Y_{t}})^{(2,2)}\right|\leq\sup_{s\in[0,1]}
|\kappa_{t}(s)^{-1}\dot{\kappa}_{t}(s)|_{\infty}\int_{0}^{1}ds
\int_{Y_{t}}|R_{t}(s)|^{2}\,dv_{t}(s).
\tag 5.2
$$
Let $\tau_{t}(s)$ the scalar curvature of $(Y_{t},\kappa_{t}(s))$. 
As is well known (cf. [G]),
$$
24=\chi(Y_{t})
=\int_{Y_{t}}\frac{1}{32\pi^{2}}\left(|R_{t}(s)|^{2}-
4|\Ric(\kappa_{t}(s))|^{2}+\tau_{t}(s)^{2}\right)dv_{t}(s).
\tag 5.3
$$
Let $\theta_{1},\theta_{2}$ be a local unitary frame and write 
$\Ric(\kappa_{t}(s))=i(\rho_{1}\theta_{1}\bar{\theta}_{1}+
\rho_{2}\theta_{2}\bar{\theta}_{2})$ $(\rho_{1},\rho_{2}\in\Bbb R)$. 
Then, $\tau_{t}(s)=\rho_{1}+\rho_{2}$ and 
$\Ric(\kappa_{t}(s))^{2}=2\rho_{1}\rho_{2}\,dv_{t}(s)$. Thus,
$$
\int_{Y_{t}}\tau_{t}(s)^{2}dv_{t}(s)
=\int_{Y_{t}}|\Ric(\kappa_{t}(s))|^{2}\,dv_{t}(s)+
\int_{X_{t}}c_{1}(Y_{t})^{2}
=\int_{Y_{t}}|\Ric(\kappa_{t}(s))|^{2}\,dv_{t}(s)
\tag 5.4
$$
because $c_{1}(Y_{t})$ is cohomologous to zero which,
together with (5.3), yields
$$
\int_{Y_{t}}|R_{t}(s)|^{2}dv_{t}(s)=
32\pi^{2}\left(24+3\int_{Y_{t}}|\Ric(\kappa_{t}(s))|^{2}dv_{t}(s)\right).
\tag 5.5
$$
The assertion follows from (5.2) and (5.5) together with Theorem 4.2.\qed
\enddemo

\proclaim{Lemma 5.2}
For any $t\in D\backslash\{0\}$, one has
$$
\left|\int_{Y_t}\widetilde{\Td}(Y_t;\kappa_{Y_t},\tilde{\kappa}_{Y_t})
\right|\leq C.
$$
\endproclaim

\demo{Proof}
Let $\Delta_{t}(s)$ be the Laplacian of $(Y_{t},\kappa_{t}(s))$. 
Differentiating (4.28) by $s$, we get
$$
\Delta_{t}(s)\dot{\varphi}_{t}(s)=F_{t}+\dot{a}_{t}(s),\quad
\int_{Y_{t}}\dot{\varphi}_{t}(s)\,\kappa_{t}(0)^{2}=0.
\tag 5.6
$$
Multiplying $|\dot{\varphi}_{t}(s)|^{p-2}\dot{\varphi}_{t}(s)$ to the both
hand sides, the integration by parts yields 
$$
\int_{X_{t}}\left|d|\dot{\varphi}_{t}(s)|^{\frac{p}{2}}
\right|^{2}_{\kappa_{t}(s)}dv_{t}(s)\leq
\frac{p^{2}}{4(p-1)}\,\int_{X_{t}}|s\,F_{t}+a_{t}(s)|\,
|\dot{\varphi}_{t}(s)|^{p-1}\,dv_{t}(s)
\tag 5.7
$$
for any $p\geq2$. By Theorem 4.2, (5.7) is also valid after changing 
the metric from $\kappa_{t}(s)$ to $\kappa_{t}(0)$. Applying Moser's 
iteration argument ([Ko, pp.298-299]), we get
$$
|\dot{\varphi}_{t}(s)|_{\infty}\leq 
C\,\|\dot{\varphi}_{t}(s)\|_{L^{2},\kappa_{t}(s)}
\tag 5.8
$$
where $C$ depends only on the constants in Theorem 4.2, 
the Sobolev constant of $(Y_{t},\kappa_{t}(0))$ and 
$|s\,F_{t}+a_{t}(s)|_{\infty}$. By Theorems 4.1 and 4.2, $C$ is 
uniformly bounded because so is the Sobolev constant of 
$(Y_{t},\kappa_{t}(0))$. As in the proof of Theorem 4.1, we find that 
$\lambda_{t}$, the first eigenvalue of $\Delta_{t}(0)$, is also 
uniformly bounded away from zero; $\lambda_{t}\geq\lambda>0$ which, 
together with (5.7), $p=2$ and Theorem 4.2, yields 
$$
\|\dot{\varphi}_{t}(s)\|_{L^{2}}^{2}\leq\lambda^{-1}\
\|d\dot{\varphi}_{t}(s)\|_{L^{2}}^{2}\leq
C\,\lambda^{-1}\|F_{t}+\dot{a}_{t}(s)\|_{L^{2}}\,
\|\dot{\varphi}_{t}(s)\|_{L^{2}}\leq 
C\,|t|\,\|\dot{\varphi}_{t}(s)\|_{L^{2}}
\tag 5.9
$$
where norms are those relative to $\kappa_{t}(0)$.
Comparing (5.8) and (5.9), we get
$$
|\dot{\varphi}_{t}(s)|_{\infty}\leq C\,|t|.
\tag 5.10
$$
Take $x\in Y_{t}\cap W$, and let $(B,w)$ be the coordinate neighborhood 
centered at $x$ as in the proof of Theorem 4.3. 
Put $\hat{\kappa}_{t}(s):=q_{t}(x)^{-1}\kappa_{t}(s)$ for the 
rescaled metric on $B$. Let 
$\hat{\Delta}_{t}(s)=q_{t}(x)\,\Delta_{t}(s)$ be the Laplacian of 
$(X_{t},\hat{\kappa}_{t}(s))$. Then, (5.6) becomes
$$
\hat{\Delta}_{t}(s)\dot{\varphi}_{t}(s)=q_{t}(x)\,(F_{t}+\dot{a}_{t}(s)).
\tag 5.11
$$
By Proposition 4.7, Theorem 4.2 and (5.10), the Schauder estimate 
([G-T, Theorem 6.2]) applied to (5.11) yields (one can apply it because
the geometry of $(B(0,r_{0}),\hat{\kappa}_{t}(s))$ is bounded as shown
in the proof of Theorem 4.3.)
$$
\aligned
|\dot{\varphi}_{t}(s)|_{C^{2+\alpha}(B(0,r_{0}))}
&\leq
C(|\dot{\varphi}_{t}(s)|_{\infty}+|q_{t}(x)\,(F_{t}+\dot{a}_{t}(s))|
_{C^{\alpha}(B(0,2r_{0}))})\\
&\leq C(|t|+q_{t}(x))\leq C\,q_{t}(x).
\endaligned
\tag 5.12
$$
As $\kappa_{t}(s)$ converges smoothly in $t$ to a K\"ahler metric on 
$Y_{0}$ outside of $W$ and thus 
$|\dot{\varphi}_{t}(s)|_{C^{2}(Y_{t}\backslash W)}\leq C$, it follows
from (5.12) that, for any $x\in Y_{t}\cap W$,
$$
\left|\partial_{Y_{t}}\bar{\partial}_{Y_{t}}\dot{\varphi}_{t}(s)(x)
\right|_{\hat{\kappa}_{t}(s)}\leq C\,q_{t}(x).
\tag 5.13
$$
Rewriting (5.13) in terms of $\kappa_{t}(s)$, we get the assertion by
Lemma 5.1 because
$$
\left|\kappa_{t}(s)^{-1}\dot{\kappa}_{t}(s)(x)\right|_{\kappa_{t}(s)}=
|\partial_{Y_{t}}\bar{\partial}_{Y_{t}}
\dot{\varphi}_{t}(s)(x)|_{\kappa_{t}(s)}
\leq C\,q_{t}(x)^{-1}\,q_{t}(x)\leq C.\qed
\tag 5.14
$$
\enddemo

\proclaim{Lemma 5.3}
For any $t\in D\backslash\{0\}$, one has
$$
\left|\int_{Y_t}\widetilde{\Td}(Y_t;\tilde{\kappa}_{Y_t},\omega_{Y_t})
\right|\leq C.
$$
\endproclaim

\demo{Proof}
Since (5.2) and (5.5) are also valid for $\tilde{\kappa}_{t}(s)$, we get
$$
\left|\int_{Y_{t}}\widetilde{\Td}(Y_{t};\tilde{\kappa}_{Y_{t}},
\omega_{Y_{t}})\right|\leq
C\,\sup_{s\in[0,1]}|\tilde{\kappa}_{t}(s)^{-1}\dot{\kappa}_{t}(s)|_{\infty}
\tag 5.15
$$
As
$\partial_{s}\tilde{\kappa}_{t}(s)=\tilde{\kappa}_{Y_{t}}-\omega_{Y_{t}}$
by (4.26), we get 
$\sup_{s\in[0,1]}|\tilde{\kappa}_{t}(s)^{-1}\dot{\tilde{\kappa}}_{t}(s)
|_{\infty}\leq C$ by Theorem 4.1 which, together with (5.15), implies the
assertion.\qed
\enddemo

\proclaim{Proposition 5.1}
For any $t\in D\backslash\{0\}$, one has
$$
\left|\int_{Y_{t}}\widetilde{\Td}(Y_{t};\kappa_{Y_{t}},\omega_{Y_{t}})
\right|\leq C.
$$
\endproclaim

\demo{Proof}
Clear by Lemmas 5.2 and 5.3 because the Bott-Chern secondary class
does not depend on a path joining two metrics ([B-C], [B-G-S I (e)]).
\qed
\enddemo

Let $\pi:Z\to D$ be the family of fixed curves of $\iota$; 
$Z=\{x\in\Cal Y;\,\iota(x)=x\}$. Let $\kappa_{Z_{t}}(s)$ be the restriction 
of $\kappa_{t}(s)$ to $Z_{t}$; $\kappa_{Z_{t}}(s)=\kappa_{t}(s)|_{Z_{t}}$. 
By definition,
$$
\aligned
\,
&\int_{Z_{t}}\widetilde{c_{1}(Y_{t})\,c_{1}(Z_{t})}
(\kappa_{Y_{t}},\tilde{\kappa}_{Y_{t}})\\
&=\int_{1}^{0}ds\int_{Z_{t}}
\left(\Tr\,\kappa_{t}(s)^{-1}\dot{\kappa}_{t}(s)\right)\,
c_{1}(Z_{t},\kappa_{Z_{t}}(s))+\kappa_{Z_{t}}(s)^{-1}\dot{\kappa}_{Z_{t}}(s)
\,c_{1}(Y_{t},\kappa_{t}(s)).
\endaligned
\tag 5.16
$$

\proclaim{Lemma 5.4}
For any $t\in D\backslash\{0\}$, one has
$$
\left|\int_{Z_{t}}\widetilde{c_{1}(Y_{t})\,c_{1}(Z_{t})}
(\kappa_{Y_{t}},\tilde{\kappa}_{Y_{t}})\right|\leq C.
$$
\endproclaim

\demo{Proof}
When $o=\Sing\,\Cal Y\not\in Z$, the assertion is obvious. Thus, we assume 
that $o\in Z_{0}$. By Proposition 2.5, $Z_{t}$ is a degenerating family of
curves such that $Z_{0}$ has only one node at $o$. By (5.16), it is enough 
to show the followings;
$$
\align
(1)\quad|\Tr\,\kappa_{t}(s)^{-1}\dot{\kappa}_{t}(s)|_{\infty}\leq C\,|t|,
&\quad(2)\quad |c_{1}(Z_{t},\kappa_{Z_{t}}(s))|_{\infty,\kappa_{Z_{t}}}
\leq C\,|t|^{-1},
\\
(3)\quad |\kappa_{Z_{t}}(s)^{-1}\dot{\kappa}_{Z_{t}}(s)|_{\infty}\leq C,
\quad\,
&\quad(4)\quad |c_{1}(X_{t},\kappa_{t}(s))|_{\infty,\kappa_{t}(s)}\leq C.
\endalign
$$
As $\Tr\,\kappa_{t}(s)^{-1}\dot{\kappa}_{t}(s)=F_{t}+\dot{a}_{t}(s)=O(|t|)$
by Theorem 4.2, we get (1). (4) is proved in Theorem 4.2. 
(3) follows from (5.14) because, for any $v\in TZ_{t}$,
$$
-\left|\kappa_{t}(s)^{-1}\dot{\kappa}_{t}(s)\right|_{\infty}\leq
\kappa_{t}(s)(v,v)^{-1}\partial_{s}\kappa_{t}(s)(v,v)\leq
\left|\kappa_{t}(s)^{-1}\dot{\kappa}_{t}(s)\right|_{\infty}.
\tag 5.17
$$
As $\kappa_{t}(s)$ is $\iota$-invariant, $R_{t}(s)$ restricted to $Z_{t}$ 
splits as follows ([B-G-V, $\S 6.3$ (6.1)]);
$$
\left.R_{t}(s)\right|_{Z_t}=R(Z_t,\omega_{t}(s))\oplus R(N_{Z_t/Y_t})
\tag 5.18
$$
where $R(N_{Z_{t}/Y_{t}})$ is the curvature of the normal bundle of $Z_{t}$
relative to the induced metric, and the splitting is orthogonal. Therefore,
by Theorem 4.3, we get
$$
|c_{1}(Z_{t},\kappa_{Z_{t}}(s))|_{\infty}=
|R(Z_{t},\kappa_{Z_{t}}(s))|_{\infty}
\leq|R_{t}(s)|_{Z_{t}}|_{\infty}\leq|R_{t}(s)|_{\kappa_{t}(s)}\leq 
C\,|t|^{-1}.\qed
\tag 5.19
$$
\enddemo

\proclaim{Lemma 5.5}
For any $t\in D\backslash\{0\}$, one has
$$
\left|\int_{Z_{t}}\widetilde{c_{1}(Y_{t})\,c_{1}(Z_{t})}
(\tilde{\kappa}_{Y_{t}},\omega_{Y_{t}})\right|\leq C.
$$
\endproclaim

\demo{Proof}
In the same way as the proof of Lemma 5.4, we may assume that 
$o\in Z_{0}$. Put
$\tilde{\kappa}_{Z_{t}}(s):=\tilde{\kappa}_{t}(s)|_{Z_{t}}$. By (5.16), 
it is enough to prove the following; for any $x\in Y_{t}\cap W$,
$$
\align
(1)\quad|\Tr\,\tilde{\kappa}_{t}(s)^{-1}\dot{\tilde{\kappa}}_{t}(s)|(x)
\leq C\,q_{t}(x),
&\quad(2)\quad |c_{1}(Z_{t},\tilde{\kappa}_{Z_{t}}(s))|(x)
\leq C\,q_{t}(x)^{-1},\\
(3)\quad |\tilde{\kappa}_{Z_{t}}(s)^{-1}\dot{\tilde{\kappa}}_{Z_{t}}(s)|(x)
\leq C,\quad\quad\,\,
&\quad(4)\quad |c_{1}(Y_{t},\tilde{\kappa}_{t}(s))|(x)\leq C
\endalign
$$
where length is measured relative to $\tilde{\kappa}_{Z_{t}}$ and
$\tilde{\kappa}_{Y_{t}}$. As 
$\tilde{\kappa}_{t}(s)^{-1}\dot{\tilde{\kappa}}_{t}(s)=
\tilde{\kappa}_{Y_{t}}-\omega_{Y_{t}}$, we get (3) by Theorem 4.1
and (1) by Lemma 4.4. (4) follows from 
Proposition 4.7. Since $\tilde{\kappa}_{t}(s)$ is $\iota$-invariant,  
by the same argument as the proof of (2) in Lemma 5.2, we get (2) by 
Proposition 4.8.
\qed
\enddemo

\proclaim{Proposition 5.2}
For any $t\in D\backslash\{0\}$, one has
$$
\left|\int_{Z_{t}}\widetilde{c_{1}(Y_{t})\,c_{1}(Z_{t})}
(\kappa_{Y_{t}},\omega_{Y_{t}})\right|\leq C.
$$
\endproclaim

\demo{Proof}
Clear by Lemmas 5.4 and 5.5\qed
\enddemo

Define $\tilde{\tau}_{M}(Y_{t},\iota_{t},\omega_{Y_{t}})$ in the same 
manner as Definition 3.3. By Propositions 5.1 and 5.2, we get the 
following.

\proclaim{Theorem 5.1}
For any $t\in D\backslash\{0\}$, one has
$$
|\log\tau_{M}(Y_{t},\iota_{t})-
\log\tilde{\tau}_{M}(Y_{t},\iota_{t},\omega_{Y_{t}})|\leq C.
$$
\endproclaim

\subsubhead
\bf{5.2 Singularity of Type $\pmb{(0,3)}$ and Asymptotics of 
$\pmb{\tau_{M}}$}\rm
\endsubsubhead
In this subsection, we assume that $(Y_{0},\iota,o)$ is a singularity of 
type (0,3) in the sense of (2.11). Put 
$\tau_{M}(t):=\tau_{M}(Y_{t}/\iota_{t},\omega_{Y_{t}})$.

\proclaim{Theorem 5.2}
As $t\to 0$, one has the following asymptotic formula:
$$
\log\tau_{M}(t)=-\frac{1}{8}\log|t|^{2}+o(\log|t|).
$$
\endproclaim

\demo{Proof}
Let $\Delta_{t}$ be the Laplacian of $(Y_{t},\omega_{Y_{t}})$. As
$\omega_{Y_{t}}$ is $\iota$-invariant, $\Delta_{t}$ splits into 
$\Delta_{t}^{\pm}$ as in Lemma 3.1. Let $K_{t}^{\pm}(s,x,y)$ be 
the heat kernel of $\Delta_{t}^{\pm}$ and $K_{t}(s,x,y)$ that of 
$\Delta_{t}$. Since $K_{t}^{\pm}(s,x,y)=
K_{t}(s,x,y)/2\pm K_{t}(s,x,\iota y)/2$ and thus
$$
\Tr\,e^{-s\Delta_{t}^{+}}-\Tr\,e^{-s\Delta_{t}^{-}}
=\int_{X_{t}}K_{t}(s,x,\iota_{t}x)\,dv_{X_{t}},
\tag 5.20
$$
it follows from [B-G-V, Theorem 6.11] that there exists
$a_{i}(z,t)\in C^{\infty}(Z_{t})$ such that
$$
\Tr\,e^{-s\Delta_{t}^{+}}-\Tr\,e^{-s\Delta_{t}^{-}}
\sim\int_{Z_{t}}\left(\frac{a_{1}(z,t)}{s}+a_{0}(z,t)\right)dv_{Z_{t}}+
O(s)\quad(s\to0).
\tag 5.21
$$
Put $a_{i}(t)=\int_{Z_{t}}a_{i}(z,t)dv_{Z_{t}}$. From (5.20-21) and 
Lemma 3.1, it follows that
$$
\aligned
\tau(Y_{t}/\iota_{t},\omega_{Y_{t}})
&=\log\det\Delta_{t}^{-}-\log\det\Delta_{t}^{+}\\
&=\int_{0}^{1}\frac{ds}{s}\left(\Tr\,e^{-s\Delta_{t}^{+}}-
\Tr\,e^{-s\Delta_{t}^{-}}-\frac{a_{-1}(t)}{s}-a_{0}(t)\right)\\
&\quad+\int_{1}^{\infty}\frac{ds}{s}
\left(\Tr\,e^{-s\Delta_{t}^{+}}-\Tr\,e^{-s\Delta_{t}^{-}}-1\right)+
a_{-1}(t)-\Gamma'(1)\,a_{0}(t)+\Gamma'(1)\\
&=\int_{0}^{1}\frac{ds}{s}\left\{\int_{Y_{t}}K_{t}(s,x,\iota x)\,dv_{Y_{t}}-
\int_{Z_{t}}\left(\frac{a_{1}(z,t)}{s}+a_{0}(z,t)\right)dv_{Z_{t}}\right\}\\
&\quad+\int_{1}^{\infty}\frac{ds}{s}\left\{\int_{Y_{t}}K_{t}(s,x,\iota x)
dv_{Y_{t}}-1\right\}+a_{-1}(t)-\Gamma'(1)\,a_{0}(t)+\Gamma'(1).
\endaligned
\tag 5.22
$$
Let $(X,\omega_{X})$ be the ALE instanton (see (6.1) below). By (4.25), 
we have 
\newline
$(Y_{t}\cap W,\omega_{Y_{t}})\cong(X\cap B(|t|^{-1/2}),|t|\,\omega_{X})$.
Let $K(s,x,y)$ be the heat kernel of $(X,\omega_{X})$. As $\omega_{Y_{t}}$ 
is Ricci-flat on $Y_{t}\cap W$, it follows from [C-L-Y, $\S 3,4$] that there 
exist constants $C,\gamma>0$ such that
$$
0<K_{t}(s,x,y)\leq C\,s^{-2}e^{-\frac{\gamma d(x,y)}{s}},\quad
|d_{x}K_{t}(s,x,y)|\leq C\,s^{-\frac{5}{2}}e^{-\frac{\gamma d(x,y)}{s}}
\tag 5.23
$$
for any $t\in D\backslash\{0\}$, $s\in(0,1]$ and $x,y\in Y_{t}\cap W$,
and one has the same estimates for $K(s,x,y)$ for any $s>0$ and $x,y\in X$. 
By Duhamel's principle (cf. [B-B, pp.63-67]) together with (5.23), there 
exist $c,C>0$ such that
$$
\left|\int_{Y_{t}\cap W}K_{t}(s,x,\iota x)\,dv_{Y_{t}}-
\int_{X\cap B(|t|^{-\frac{1}{2}})}K(|t|^{-1}s,x,\iota x)\,dv_{X}\right|\leq
C\,e^{-\frac{c}{s}}
\tag 5.24
$$
for any $t\in D\backslash\{0\}$ and $s\in(0,1]$. Put 
$J(T):=\int_{0}^{T}s^{-1}ds\int_{B(\sqrt{T})}K(s,x,\iota x)\,dv_{X}$.
Since $W\cap Z_{t}=\emptyset$, there appears no contribution from $Z_{t}$ to
the divergence of $\tau(Y_{t}/\iota_{t},\omega_{Y_{t}})$ which, together 
with (5.24), yields
$$
\int_{0}^{1}\frac{ds}{s}\int_{Y_{t}\cap W}K_{t}(s,x,\iota x)\,dv_{Y_{t}}
=J(|t|^{-1})+O\left(\int_{0}^{1}\frac{ds}{s}e^{-\frac{c}{s}}\right)
=J(|t|^{-1})+O(1).
\tag 5.25
$$
As $\omega_{Y_{t}}$ is a smooth family of metrics outside of 
$Y_{t}\cap W$, one has
$$
\left|\int_{0}^{1}\frac{ds}{s}\left(\int_{Y_{t}\backslash W}
K_{t}(s,x,\iota x)dv_{Y_{t}}-
\int_{Z_{t}}\left(\frac{a_{1}(z,t)}{s}+a_{0}(z,t)\right)dv_{Z_{t}}\right)
\right|\leq C.
\tag 5.26
$$
Let $\lambda_{k}(t)$ be the $k$-th eigenvalue of $\Delta_{t}$. Since the 
Sobolev constant of $(Y_{t},\omega_{Y_{t}})$ is uniformly bounded 
away from zero, it follows from  [C-L] that there exists $N\geq 0$ 
such that, for any $t\in D$ and $k\geq N$, one has
$\lambda_{k}(t)\geq C\,k^{1/2}$. By [Yo1, Theorem 5.1], 
$\lambda_{k}(t)$ is a continuous function on $D$ for any $k\geq 0$ and 
$\lambda_{1}(0)>0$.  Thus, we get the bound
$\lambda_{k}(t)\geq C\,k^{1/2}$ and
$0\leq\Tr\,e^{-s\Delta_{t}}-1\leq C\,e^{-\lambda_{1}(0)s/2}$
for $s\geq1$, which yields
$$
\aligned
\int_{1}^{\infty}\frac{ds}{s}\left|\int_{Y_{t}}K(s,x,\iota x)\,dv_{Y_{t}}-1
\right|
&=\int_{1}^{\infty}\frac{ds}{s}\left|\Tr\,e^{-s\Delta_{t}^{+}}-
\Tr\,e^{-s\Delta_{t}^{-}}-1\right|\\
&\leq\int_{1}^{\infty}\frac{ds}{s}
\left(\Tr\,e^{-s\Delta_{t}}-1\right)\leq C.
\endaligned
\tag 5.27
$$
Let $K(s,x,y;\lambda\omega_{X})$ be the heat kernel and 
$B(r;\lambda\omega_{X})$
the metric ball of radius $r$ of $(X,\lambda\omega_{X})$. Let 
$ds^{2}_{\Bbb C^{2}}$ be the Euclidean metric of $\Bbb C^{2}$ and
$K_{\Bbb C^{2}/\pm1}(s,x,y)$ be the heat kernel of 
$(\Bbb C^{2}/\pm 1,ds^{2}_{\Bbb C^{2}})$. Since
$(X,\lambda\omega_{X},\iota,o)$ converges to 
$(\Bbb C^{2}/\pm 1,ds^{2}_{\Bbb C^{2}},\iota',0)$ as $\lambda\to 0$ 
where $\iota'(z_{1},z_{2})=(\sqrt{-1}z_{1},\sqrt{-1}z_{2})$, we get
$$
\aligned
\,&
\lim_{T\to\infty}\frac{J(T)}{\log T}
=\lim_{T\to\infty}\int_{0}^{1}d\sigma
\int_{B(T^{\frac{1-\sigma}{2}};T^{-\sigma}\omega_{X})}
K(1,x,\iota x;T^{-\sigma}\omega_{X})\,dv_{T^{-\sigma}\omega_{X}}\\
&=\int_{0}^{1}d\sigma\int_{\Bbb C^{2}/\pm 1}
K_{\Bbb C^{2}/\pm1}(1,x,\iota'x)\,dv_{\Bbb C^{2}}
=\frac{1}{2}\int_{\Bbb C^{2}}\left(
\frac{e^{-\frac{|z-iz|^{2}}{4s}}}{(4\pi s)^{2}}
+\frac{e^{-\frac{|z+iz|^{2}}{4s}}}{(4\pi s)^{2}}\right)\,dv_{\Bbb C^{2}}
=\frac{1}{4}.
\endaligned
\tag 5.28
$$
As $\log\tau_{M}(t)=J(|t|^{-1})+O(1)\quad(t\to 0)$ by (5.22), (5.26-27) 
and (1.1), the assertion follows from (5.28) because $\pi:Z\to D$ 
is a smooth morphism and thus $\log\tau(Z_{t},\omega_{Z_{t}})$ is 
bounded on $D$.\qed
\enddemo

\beginsection
$\S6$. Singularity of $\pmb{\tau_{M}}$ along the Discriminant Locus

\par
\subsubhead
\bf{6.1 Heat Kernel on the ALE Instanton}\rm
\endsubsubhead
Consider the following affine quadric with involution and the fixed locus;
$$
\aligned
\,&
X:=\{(z_{1},z_{2},z_{3})\in\Bbb C^{3};\,z_{1}z_{2}-z_{3}^{2}=1\},\quad
\iota(z_{1},z_{2},z_{3})=(z_{1},z_{2},-z_{3}),\\
\,&
X^{\iota}=Z:=\{(z_{1},z_{2})\in\Bbb C^{2};\,z_{1}z_{2}=1\}.
\endaligned
\tag 6.1
$$ 
Let $\omega_{X}=\partial_{X}\bar{\partial}_{X}q$ be the Kronheimer's ALE 
metric on $X$ ([Kr]) as in $\S4.2$, and put 
$\omega_{Z}:=\omega_{X}|_{Z}$ for its restriction to $Z$. Fix $o\in X$ and 
put $r(x):=\hbox{dist}(o,x)$. Let $B(\rho)$ be the metric ball of radius 
$\rho$ centered at $o$. By Kronheimer [Kr], there exist $c,C_{\alpha}>0$ 
with the following properties; 
(1) For any $y\in X$, the injectivity radius at
$y$ is greater than $j_{y}=c(1+r(y))$. (2) On the metric ball 
$B(y,j_{y})$ with the normal coordinates $x=(x_{1},\cdots,x_{4})$, 
the metric tensor $g(x)=\sum g_{ij}\,(x)dx_{i}dx_{j}$ satisfies
$$
\sup_{x\in B(y,j_{y})}\left|\partial^{\alpha}(g_{ij}(x)-\delta_{ij})\right|
\leq C_{\alpha}(1+r(y))^{-|\alpha|}.
\tag 6.2
$$
For simplicity, we assume $c=2$ by considering $(2c)^{-1}\omega_{X}$ 
if necessary. Let $K(t,x,y)$ be the heat kernel of $(X,\omega_{X})$ and 
define its parametrix ([B-G-V, Theorem 2.26]);
$$
p(t,x,y):=
(4\pi t)^{-\frac{m}{2}}e^{-\frac{d(x,y)^{2}}{4t}}
\{u_{0}(x,y)+t\,u_{1}(x,y)\}
\tag 6.3
$$
where $u_{i}(\cdot,y)\in C^{\infty}(B(y,j_{y}))$. (In [B-G-V], $\Phi_{i}$
is used instead of $u_{i}$.)

\proclaim{Lemma 6.1}
For any $y\in X$ and $(x,t)\in B(y,j_{y}/2)\times[0,1+|y|^{2}]$, one has
$$
|K(t,x,y)-p(t,x,y)|\leq C(c,C_{\alpha})\,(1+r(y)^{2})^{-2}
e^{-\frac{\gamma\,d(x,y)^{2}}{t}}.
$$
\endproclaim

\demo{Proof}
Let $B(r)$ be the ball of radius $r<1$ in $\Bbb R^{4}$. 
Let $(M,g)$ be a Riemannian 4-manifold such that $B(1)$ is 
embedded into $M$ and $g=\sum g_{ij}\,dx_{i}dx_{j}$ on $B(1)$. Suppose 
that $|\partial^{k}(g_{ij}(x)-\delta_{ij})|\leq A_{k}$ for any $k\geq0$ 
and $x\in B(1)$. Let $L(t,x,y)$ be the heat kernel of $(M,g)$. 
By [C-L-Y, $\S3,4$], there exists $r,C,\gamma>0$ depending only on $A_{k}$ 
($k\leq 6$) and $r$ such that, for any $(x,y)\in B(r)\times B(r)$ and 
$t\in(0,1]$, 
$$
0<L(t,x,y)\leq C\,t^{-2}e^{-\frac{\gamma d(x,y)^{2}}{t}},\quad
|d_{x}L(t,x,y)|\leq C\,t^{-\frac{5}{2}}e^{-\frac{\gamma d(x,y)^{2}}{t}}.
\tag 6.4
$$
Let $L'(t,x,y)$ be the Dirichlet heat kernel of $(B(1),g)$ which satisfies 
the same estimates (6.4). By Duhamel's principle, for any 
$(x,y)\in B(r)\times B(r)$ and $t\in(0,1]$,
$$
|L(t,x,y)-L'(t,x,y)|\leq C\,e^{-\frac{c_{0}}{t}}
\tag 6.5
$$
where $c_{0},C>0$ depend only on $C_{k}$. Since $C^{k}$-norm of $g_{ij}$ is 
bounded by $A_{k}$ ($k=6$ is enough), it follows from the usual asymptotic
expansion of the heat kernel that there exits a constant $C,\gamma>0$ 
depending only on $A_{k}$ ($k\leq 6$) and $r$ such that, for any 
$(x,y)\in B(r)\times B(r)$ and $t\in(0,1]$,
$$
|L'(t,x,y)-p_{L'}(t,x,y)|\leq C\,e^{-\frac{\gamma d(x,y)^{2}}{t}}
\tag 6.6
$$
where $p_{L'}(t,x,y)$ is defined by the formula (6.3) for $L'$ and $d(x,y)$ 
is the distance relative to the metric $g$. Comparing (6.5) and (6.6), 
we find 
$$
|L(t,x,y)-p_{L}(t,x,y)|\leq C\,e^{-\frac{\gamma d(x,y)^{2}}{t}}.
\tag 6.7
$$
For any $y\in X$, consider the rescaled space $(X,j_{y}^{-2}\omega_{X})$ 
and its heat kernel $L(t,x,y)$. As $j_{y}^{4}\,K(j_{y}^{2}t,x,y)=L(t,x,y)$ 
and $j_{y}^{4}\,p(j_{y}^{2}t,x,y)=p(t,x,y)$, it follows from (6.7) that
$$
|K(j_{y}^{2}t,x,y)-p(j_{y}^{2}t,x,y)|\leq 
C\,(1+r(y))^{-4}\exp
\left(-\frac{\gamma d(x,y;\omega_{X})^{2}}{j_{y}^{2}t}\right),
\quad(t\in(0,1]).
\tag 6.8
$$
The assertion follows from (6.8) by putting $s=j_{y}^{2}t$.\qed
\enddemo

Let $\Cal N_{Z/X}$ be the normal bundle of $Z$ in $X$ and $v$ its fiber 
coordinate. Set
\newline
$\Omega:=\{(z,v)\in\Cal N_{Z/X};\,\|v\|_{z}\leq r(z)+1\}$,
$\Omega(a):=\{(z,v)\in\Omega;\,r(exp_{z}(v))\leq a\}$.
\newline
Via the exponential map, identify $\Omega$ with a tubular neighborhood of 
$Z$ in $X$. Let $dv_{Z}$ be the volume form of $(Z,\omega_{Z})$. Then, the 
volume form of $(X,\omega_{X})$ can be written on $\Omega$ as follows;
$dv_{X}(x)=J(z,v)\,dv\,d\bar{v}\,dv_{Z}(z)$ where $x=\exp_{z}(v)$. Set
$$
q_{i}(z,v):=J(z,v)\,u_{i}(x,\iota(x)),\quad 
q_{i}^{(k)}(z,v):=\frac{1}{k!}\partial_{t}^{k}|_{t=0}q_{i}(z,tv).
\tag 6.9
$$
In the sequel, $C$ denotes a constant which depends only 
on $C_{\alpha}$ ($|\alpha|\leq6$), $c$ and $\gamma$. By definition, 
we get $q_{0}(z,0)=1$. By (6.2), the following is clear.

\proclaim{Lemma 6.2}
For any $(z,v)\in\Omega$, one has the following estimates;
\flushpar{$(1)$} 
$|q_{0}(z,v)-q_{0}(z,0)|\leq C\,|v|(1+r)^{-1}$,
$(2)$ $|\partial^{\alpha}_{v}q_{i}(z,v)|\leq C\,(1+r)^{-(2i+|\alpha|)}$,
\flushpar{$(3)$} 
$|q_{0}(z,v)-q_{0}(z)-q_{0}^{(1)}(z,v)-q_{0}^{(2)}(z,v)|
\leq C\,|v|^{3}(1+r)^{-3}$,
\flushpar{$(4)$} 
$|q_{1}(z,v)-q_{1}(z,0)|\leq C\,|v|(1+r)^{-3}$.
\newline
In $(2)$, $|\alpha|>0$ when $i=0$ and $|\alpha|\geq 0$ when $i=1$.
\endproclaim

Let $dx\,dy=\frac{\sqrt{-1}}{2}dv\,d\bar{v}$ be the volume form on the fiber
of $\Cal N_{Z/X}$. Set
\newline
$I_{1}(s,z):=\int_{|v|\leq r(z)+1}s^{-2}e^{-|v|^{2}/s}
\left(q_{0}(z,v)-q_{0}(z)-q_{0}^{(1)}(z,v)-q_{0}^{(2)}(z,v)\right)dx\,dy$.

\proclaim{Lemma 6.3}
For any $T>1$, one has
$$
\left|\int_{0}^{T}\frac{ds}{s}\int_{r\leq\sqrt{T}}I_{1}(s,z)dv_{Z}\right|
\leq C\,(\log T+1).
$$
\endproclaim

\demo{Proof}
Suppose $r\geq\sqrt{s}$. It follows from Lemma 6.2 (3) that
$$
|I_{1}(s,z)|
\leq\int_{|v|\leq r+1}C\,e^{-\frac{|v|^{2}}{s}}\frac{|v|^{3}}{s^{2}}
\frac{dx\,dy}{(1+r)^{3}}
\leq C\,\int_{\rho\leq r+1}e^{-\frac{\rho^{2}}{s}}
\frac{\rho^{4}d\rho}{s^{2}(1+r)^{3}}
\leq C\,\frac{s^{\frac{1}{2}}}{(1+r)^{3}}
\tag 6.10
$$
which yields
$$
\aligned
\left|\int_{1}^{T}\frac{ds}{s}\int_{\sqrt{s}\leq r\leq\sqrt{T}}
I_{1}(z,s)dv_{Z}\right|
&\leq C\int_{1}^{T}\frac{ds}{s}\int_{\sqrt{s}}^{\infty}
s^{\frac{1}{2}}(1+r)^{-3}r\,dr\leq C\,\log T.
\endaligned
\tag 6.11
$$
Suppose $r\leq\sqrt{s}$. It follows from Lemma 6.2 (1) that
$$
\left|\int_{|v|\leq r+1}e^{-\frac{|v|^{2}}{s}}
\left(q_{0}(z,v)-q_{0}(z)\right)\frac{dx\,dy}{s^{2}}\right|
\leq C\,\int_{|v|\leq r+1}e^{-\frac{|v|^{2}}{s}}
\frac{|v|\,dx\,dy}{s^{2}(1+r)}\leq C\,\frac{(1+r)^{2}}{s^{2}}.
\tag 6.12
$$
Similarly, using Lemma 6.2 (2), we get
$$
\left|\int_{|v|\leq r+1}e^{-\frac{|v|^{2}}{s}}
q_{0}^{(2)}(z,v)\,\frac{dx\,dy}{s^{2}}\right|
\leq C\,\int_{|v|\leq r+1}e^{-\frac{|v|^{2}}{s}}
\frac{|v|^{2}\,dx\,dy}{s^{2}(1+r)^{2}}
\leq C\,\frac{(1+r)^{2}}{s^{2}}.
\tag 6.13
$$
Since $\int_{|v|\leq r+1}e^{-\frac{|v|^{2}}{s}}v^{k}\,dx\,dy=0$ $(k>0)$, 
we get
$$
\int_{|v|\leq r+1}s^{-2}e^{-\frac{|v|^{2}}{s}}q_{0}^{(1)}(z,v)dx\,dy=0
\tag 6.14
$$
which, together with (6.12-13), yields
$$
\left|\int_{1}^{T}\frac{ds}{s}\int_{r\leq\sqrt{s}}I_{1}(z,s)dv_{Z}\right|
\leq C\int_{1}^{T}\frac{ds}{s}s^{-2}\int_{r\leq\sqrt{s}}(1+r)^{2}
r\,dr\leq C\,\log T.
\tag 6.15
$$
When $s\leq 1$, $1+r(z)\geq s$ for any $z\in Z$. Thus, by (6.10), we get
$$
\aligned
\left|\int_{0}^{1}\frac{ds}{s}\int_{r\leq\sqrt{T}}I_{1}(z,s)dv_{Z}\right|
&\leq C\int_{0}^{1}\frac{ds}{s}\sqrt{s}\int_{r\leq\sqrt{T}}(1+r)^{-3}
r\,dr\leq C.
\endaligned
\tag 6.16
$$
The assertion follows from (6.11) and (6.15-16).\qed
\enddemo

Set
$I_{2}(s,z):=\int_{|v|\leq r(z)+1}s^{-1}e^{-|v|^{2}/s}
\left(q_{1}(z,v)-q_{1}(z,0)\right)dx\,dy$.

\proclaim{Lemma 6.4}
For any $T>1$, one has
$$
\left|\int_{0}^{T}\frac{ds}{s}\int_{r\leq\sqrt{T}}I_{2}(s,z)dv_{Z}\right|
\leq C\,(\log T+1).
$$
\endproclaim

\demo{Proof}
Suppose $r\geq\sqrt{s}$. It follows from Lemma 6.2 (4) that
$$
|I_{2}(s,z)|
\leq C\int_{|v|\leq r+1}e^{-\frac{|v|^{2}}{s}}\frac{|v|}{s}
\frac{dx\,dy}{(1+r)^{3}}
\leq C\int_{\rho\leq r+1}e^{-\frac{\rho^{2}}{s}}
\frac{\rho^{2}}{s}\frac{d\rho}{(1+r)^{3}}\leq 
C\frac{\sqrt{s}}{(1+r)^{3}}
\tag 6.17
$$
which yields
$$
\left|\int_{1}^{T}\frac{ds}{s}\int_{\sqrt{s}\leq r\leq\sqrt{T}}
I_{2}(s,z)dv_{Z}\right| 
\leq C\int_{1}^{T}\frac{ds}{s}s^{\frac{1}{2}}\int_{\sqrt{s}}^{\infty}
(1+r)^{-3}r\,dr\leq C\,\log T.
\tag 6.18
$$
Suppose $r\leq\sqrt{s}$. It follows from the same lemma that
$$
|I_{2}(s,z)|
\leq C(1+r)^{-3}\int_{\rho\leq r+1}
s^{-1}e^{-\frac{\rho^{2}}{s}}\rho^{2}d\rho
\leq C\,s^{-\frac{1}{2}}(1+r)^{-1}
\tag 6.19
$$
which yields
$$
\left|\int_{1}^{T}\frac{ds}{s}\int_{r\leq\sqrt{s}}I_{2}(s,z)dv_{Z}\right| 
\leq C\int_{1}^{T}\frac{ds}{s}s^{-\frac{1}{2}}
\int_{0}^{\sqrt{s}}(1+r)^{-1}r\,dr\leq C\,\log T.
\tag 6.20
$$
When $s\leq 1$, it follows from (6.17) that
$$
\left|\int_{0}^{1}\frac{ds}{s}\int_{r\leq\sqrt{T}}I_{2}(s,z)dv_{Z}\right| 
\leq C\int_{0}^{1}\frac{ds}{s}s^{\frac{1}{2}}
\int_{0}^{\infty}(1+r)^{-3}r\,dr\leq C.
\tag 6.21
$$
The assertion follows from (6.17), (6.19) and (6.21).\qed
\enddemo

\proclaim{Lemma 6.5}
For any $T>1$, one has
$$
\left|\int_{0}^{T}\frac{ds}{s}\int_{\Omega(\sqrt{T})}
\left\{K(s,x,\iota x)-p(s,x,\iota x)\right\}dv_{X}\right|\leq
C\,(\log T+1).
$$
\endproclaim

\demo{Proof}
Put $E(s,x):=K(s,x,\iota x)-p(s,x,\iota x)$,
$I_{3}(s)=\int_{\Omega(\sqrt{T})\backslash\Omega(\sqrt{s})}E(s,x)\,dv_{X}$,
\newline
and $I_{4}(s)=\int_{\Omega(\sqrt{s})}E(s,x)\,dv_{X}$.
Suppose $1\leq s\leq T$. By Lemma 6.1, we get
$$
\aligned
|I_{3}(s)|
&\leq C\int_{\Omega(\sqrt{T})\backslash\Omega(\sqrt{s})}(1+r^{2})^{-2}
e^{-\frac{\gamma|v|^{2}}{s}}dv_{Z}dx\,dy\\
&\leq C\int_{\sqrt{s}}^{\sqrt{T}}\frac{r\,dr}{(1+r^{2})^{2}}
\int_{|v|\leq r}e^{-\frac{\gamma |v|^{2}}{s}}dx\,dy\leq
C\,\int_{\sqrt{s}}^{\sqrt{T}}\frac{dr}{r^{3}}\leq C.
\endaligned
\tag 6.22
$$
As $(X,\omega_{X})$ is Ricci-flat, we get a bound 
([L-Y, Theorem 3.2]); $K(t,x,y)\leq C\,t^{-2}$ for any $t>0$ and 
$x,y\in X$. Then, we have
$$
\aligned
|I_{4}(s)|
&\leq C\int_{r\leq\sqrt{s}}dv_{Z}\int_{|v|\leq r}\left(\frac{1}{s^{2}}+
\frac{1}{s(1+r^{2})}\right)dx\,dy\\
&\leq C\frac{1}{s^{2}}\int_{r\leq\sqrt{s}}r^{2}dv_{Z}+
C\int_{r\leq\sqrt{s}}dv_{Z}
\leq C\frac{1}{s^{2}}\int_{r\leq\sqrt{s}}r^{3}dr+
C\frac{1}{s}\int_{r\leq\sqrt{s}}r\,dr\leq C.
\endaligned
\tag 6.23
$$
Suppose that $s\leq 1$. By Lemma 6.1, we get
$$
|I_{3}(s)+I_{4}(s)|
\leq C\int_{\Omega(\sqrt{T})}
\frac{e^{-\frac{\gamma |v|^{2}}{s}}dv_{X}}{(1+r^{2})^{2}}
\leq C\,s\int_{0}^{\sqrt{T}}\frac{r\,dr}{(1+r^{2})^{2}}
\int_{\Bbb C}e^{-\frac{\gamma |v|^{2}}{s}}\frac{dx\,dy}{s}\leq C\,s.
\tag 6.24
$$
Together with (6.23-25), we get
$$
\left|\int_{0}^{T}\frac{ds}{s}(I_{3}(s)+I_{4}(s))\right|
\leq C\int_{1}^{T}\frac{ds}{s}+C\int_{0}^{1}\frac{ds}{s}s
\leq \,C(\log T+1).\qed
\tag 6.25
$$
\enddemo

\proclaim{Lemma 6.6}
Put $q_{0;1\bar{1}}(z,0)=\partial_{v}\partial_{\bar{v}}|_{v=0}q_{0}(z,v)$. 
For any $T>1$, one has
$$
\aligned
\,&\int_{B(\sqrt{T})}K(s,x,\iota x)dv_{X}-
\frac{1}{16\pi s}\int_{r\leq\sqrt{T}}dv_{Z}-
\frac{1}{16\pi}\int_{r\leq\sqrt{T}}q_{0;1\bar{1}}(z,0)dv_{Z}\\
&=\int_{B(\sqrt{T})\backslash\Omega(\sqrt{T})}K(s,x,\iota x)dv_{X}+I_{3}(s)
+I_{4}(s)-\frac{1}{16\pi s}\int_{r\leq\sqrt{T}}e^{-\frac{r^{2}}{s}}dv_{Z}\\
&\quad-\frac{1}{16\pi}\int_{r\leq\sqrt{T}}q_{0;1\bar{1}}(z,0)\,
e^{-\frac{r^{2}}{s}}\left(1+\frac{r^{2}}{s}\right)dv_{Z}+
\int_{r\leq\sqrt{T}}(I_{1}(s,z)+I_{2}(s,z))\,dv_{Z}.
\endaligned
$$
\endproclaim

\demo{Proof}
By the definition of $p(s,x,y)$, $J(z,v)$ and $q_{i}^{(k)}(z,v)$, we get
$$
\aligned
\,&
(4\pi s)^{2}\,p(s,x,\iota x)\,J(z,v)
=e^{-\frac{|v|^{2}}{s}}\left\{q_{0}(z,0)+q_{0}^{(1)}(z,v)+q_{0}^{(2)}(z,v)
\right\}+s\,e^{-\frac{|v|^{2}}{s}}q_{1}(z,0)\\
&\quad+e^{-\frac{|v|^{2}}{s}}\left\{q_{0}(z,v)-q_{0}(z,0)-q_{0}^{(1)}(z,v)-
q_{0}^{(2)}(z,v)\right\}
+s\,e^{-\frac{|v|^{2}}{s}}\left\{q_{1}(z,v)-q_{1}(z,0)\right\}
\endaligned
\tag 6.26
$$
which, together with $q_{0}(z,0)=1$, $q_{1}(z,0)=0$ ($\Ric(\omega_{X})=0$),
(6.14) and the definition of $I_{1}$ and $I_{2}$, yields
$$
\aligned
\,&\int_{\Omega(\sqrt{T})}p(s,x,\iota x)dv_{X}
-\int_{r\leq\sqrt{T}}(I_{1}(s,z)+I_{2}(s,z))\,dv_{Z}\\
&=\int_{\Omega(\sqrt{T})}\left\{(4\pi s)^{-2}e^{-\frac{|v|^{2}}{s}}
+(4\pi s)^{-2}e^{-\frac{|v|^{2}}{s}}q_{0;1\bar{1}}(z,0)\,|v|^{2}\right\}
dv_{X}\\
&=\int_{r\leq\sqrt{T}}dv_{Z}\int_{|v|\leq r+1}e^{-\frac{r^{2}}{s}}
\frac{dx\,dy}{(4\pi s)^{2}}+
\int_{r\leq\sqrt{T}}q_{0;1\bar{1}}(z,0)dv_{Z}\int_{|v|\leq r+1}
e^{-\frac{r^{2}}{s}}r^{2}\frac{dx\,dy}{(4\pi s)^{2}}\\
&=\frac{1}{16\pi s}\int_{r\leq\sqrt{T}}dv_{Z}-
\frac{1}{16\pi s}\int_{r\leq\sqrt{T}}e^{-\frac{(r+1)^{2}}{s}}dv_{Z}\\
&\quad+\frac{1}{16\pi}\int_{r\leq\sqrt{T}}q_{0;1\bar{1}}(z,0)dv_{Z}-
\frac{1}{16\pi}\int_{r\leq\sqrt{T}}q_{0;1\bar{1}}(z,0)\,
e^{-\frac{(r+1)^{2}}{s}}\left(1+\frac{(r+1)^{2}}{s}\right)dv_{Z}.
\endaligned
\tag 6.27
$$
Therefore, we get
$$
\aligned
\,&\int_{\Omega(\sqrt{T})}p(s,x,\iota x)dv_{X}-
\frac{1}{16\pi s}\int_{r\leq\sqrt{T}}dv_{Z}-
\frac{1}{16\pi}\int_{r\leq\sqrt{T}}q_{0;1\bar{1}}(z,0)dv_{Z}\\
&=-\frac{1}{16\pi s}\int_{r\leq\sqrt{T}}e^{-\frac{(r+1)^{2}}{s}}dv_{Z}-
\frac{1}{16\pi}\int_{r\leq\sqrt{T}}q_{0;1\bar{1}}(z,0)\,
e^{-\frac{(r+1)^{2}}{s}}\left(1+\frac{(r+1)^{2}}{s}\right)dv_{Z}\\
&\quad+\int_{r\leq\sqrt{T}}(I_{1}(s,z)+I_{2}(s,z))dv_{Z}
\endaligned
\tag 6.28
$$
which, together with the definition of $I_{3},I_{4}$ and the following, 
yields the assertion;
$$
\aligned
\,&\int_{B(\sqrt{T})}K(s,x,\iota x)dv_{X}-
\frac{1}{16\pi s}\int_{r\leq\sqrt{T}}dv_{Z}-
\frac{1}{16\pi}\int_{r\leq\sqrt{T}}q_{0;1\bar{1}}(z,0)dv_{Z}\\
&=\int_{B(\sqrt{T})\backslash\Omega(\sqrt{T})}K(s,x,\iota x)dv_{X}+
\int_{\Omega(\sqrt{T})}\{K(s,x,\iota x)-p(s,x,\iota x)\}dv_{X}\\
&\quad+\int_{\Omega(\sqrt{T})}p(s,x,\iota x)dv_{X}-
\frac{1}{16\pi s}\int_{r\leq\sqrt{T}}dv_{Z}-
\frac{1}{16\pi}\int_{r\leq\sqrt{T}}q_{0;1\bar{1}}(z,0)dv_{Z}.\qed
\endaligned
\tag 6.29
$$
\enddemo

\proclaim{Lemma 6.7}
For any $s\in(0,1]$, one has
$$
\left|\int_{B(\sqrt{T})}K(s,x,\iota x)dv_{X}-
\frac{1}{s}\int_{r\leq\sqrt{T}}\frac{dv_{Z}}{16\pi}-
\int_{r\leq\sqrt{T}}q_{0;1\bar{1}}(z,0)\frac{dv_{Z}}{16\pi}\right|\leq 
C\,\sqrt{s}.
$$
\endproclaim

\demo{Proof}
Put
$I_{5}(s):=
\int_{B(\sqrt{T})\backslash\Omega(\sqrt{T})}K(s,x,\iota x)dv_{X}$,
$I_{6}(s):=\int_{r\leq\sqrt{T}}s^{-1}e^{-\frac{r^{2}}{s}}dv_{Z}$, and 
\newline
$I_{7}(s):=\int_{r\leq\sqrt{T}}q_{0;1\bar{1}}(z,0)\,
e^{-\frac{r^{2}}{s}}\left(1+\frac{r^{2}}{s}\right)dv_{Z}$.
For any $x\in B(\sqrt{T})\backslash\Omega(\sqrt{T})$, we get
$2d(x,\iota x)\geq r(x)+1$ which, together with (6.23), 
yields that, for any $s>0$,
$$
\aligned
\,&
|I_{5}(s)|\leq C\int_{X}s^{-2}e^{-\frac{r^{2}+1}{s}}dv_{X}\leq 
C\,e^{-\frac{1}{s}},\quad
|I_{6}(s)|\leq C\int_{Z}s^{-1}e^{-\frac{r^{2}+1}{s}}dv_{Z}\leq
C\,e^{-\frac{1}{s}},\\
&|I_{7}(s)|
\leq C\int_{Z}e^{-\frac{(r+1)^{2}+1}{2s}}
\left(1+\frac{(1+r)^{2}}{s}\right)\frac{dv_{Z}}{(1+r)^{2}}
\leq C\int_{0}^{\infty}\frac{e^{-\frac{(1+r)^{2}+1}{4s}}r\,dr}{(1+r)^{2}}
\leq C\,e^{-\frac{1}{8s}}.
\endaligned
\tag 6.30
$$
By (6.10), (6.17) and (6.24), we get
$\int_{r\leq\sqrt{T}}(|I_{1}(s,z)|+|I_{2}(s,z)|)\,dv_{Z}\leq C\sqrt{s}$
and $\quad |I_{3}(s)+I_{4}(s)|\leq C\,s$ for any $s\in (0,1]$
which, together with (6.30) and Lemma 6.6, yields the 
assertion.
\qed
\enddemo

\proclaim{Proposition 6.1}
For any $T>1$, one has
$$
\left|\int_{0}^{T}\frac{ds}{s}\left
\{\int_{B(\sqrt{T})}K(s,x,\iota x)dv_{X}-
\int_{r\leq\sqrt{T}}\left(\frac{dv_{Z}}{16\pi s}+
q_{0;1\bar{1}}(z,0)\frac{dv_{Z}}{16\pi}\right)\right\}\right|
\leq C(\log T+1).
$$
\endproclaim

\demo{Proof}
By Lemmas 6.3-7 together with (6.30), we get
$$
\hbox{L.H.S.}
\leq C\,(\log T+1)+\int_{0}^{T}\frac{ds}{s}(|I_{5}(s)|+|I_{6}(s)|+
|I_{7}(s)|)\leq
C\,(\log T+1).\qed
\tag 6.31
$$
\enddemo

\subsubhead
\bf{6.2 Singularity of Type $\pmb{(2,1)}$ and Asymptotics of 
$\pmb{\tau_{M}}$}\rm
\endsubsubhead
Let us consider the same situation as in $\S 5.2$. Here, we assume that 
$(Y_{0},\iota,o)$ is of type (2,1) in the sense of (2.11). 
By Proposition 2.3, $Z_{0}$ has only one node at $o$. In the sequel, 
we use the same notations as in $\S 5.2$.

\proclaim{Lemma 6.8}
There exists a function $I(T)$ defined for $T>1$ such that, as $t\to 0$,
$$
\log\tau(Y_{t}/\iota_{t},\omega_{Y_{t}})=I(|t|^{-1})+O(1),\quad
|I(T)|\leq C\,(\log T+1).
$$
\endproclaim

\demo{Proof}
Let $K_{t}(s,x,y)$ be the heat kernel of $(Y_{t},\omega_{Y_{t}})$ and
$$
K_{t}(s,x,\iota_{t}x)\sim(a_{1}(x,t)\,s^{-1}+a_{0}(x,t)+O(s))
\delta_{Z_{t}}(x)\quad(s\to 0)
\tag 6.32
$$
its pointwise asymptotic expansion ([B-G-V, Theorem 6.11]) where 
$\delta_{Z_{t}}$ is the Dirac $\delta$-function supported along $Z_{t}$ and 
$a_{i}(z,t)$ is a smooth function on $Z_{t}$. By Lemma 3.1, we get in the 
same manner as (5.25)
$$
\aligned
\log\tau(Y_{t}/\iota_{t},\omega_{Y_{t}})
&=\int_{0}^{1}\frac{ds}{s}\left\{\int_{Y_{t}}K_{t}(s,x,\iota x)\,dv_{Y_{t}}-
\int_{Z_{t}}\left(\frac{a_{1}(z,t)}{s}+a_{0}(z,t)\right)dv_{Z_{t}}\right\}\\
&\quad+\int_{1}^{\infty}\frac{ds}{s}\left\{\int_{Y_{t}}K_{t}(s,x,\iota x)
dv_{Y_{t}}-1\right\}+a_{-1}(t)-\Gamma'(1)(a_{0}(t)-1)
\endaligned
\tag 6.33
$$
where $a_{i}(t):=\int_{Y_{t}}a_{i}(x,t)\,dv_{Y_{t}}$. By the same 
argument as in the proof of Theorem 5.4 using (5.23) and Duhamel's 
principle, we get (5.24). Put 
$$
I(T):=\int_{0}^{T}\frac{ds}{s}\left\{\int_{B(\sqrt{T})}
K(s,x,\iota x)dv_{X}-\frac{1}{s}\int_{r\leq\sqrt{T}}\frac{dv_{Z}}{16\pi}-
\int_{r\leq\sqrt{T}}q_{0;1\bar{1}}(z,0)\frac{dv_{Z}}{16\pi}\right\}.
$$
By Lemma 6.7 and (5.24), the integrand of $ds/s$ in $I(T)$ should
coincides with the asymptotic expansion (6.33) on $W\cap Y_{t}$ and we get 
$\int_{Y_{t}\cap W}a_{0}(z,t)dv_{Z_{t}}=\int_{r\leq|t|^{-\frac{1}{2}}}
dv_{Z}/16\pi$ and
$\int_{Y_{t}\cap W}a_{1}(z,t)dv_{Z_{t}}=\int_{r\leq|t|^{-\frac{1}{2}}}
q_{0;1\bar{1}}(z,0)dv_{Z}/16\pi$
which, together with (6.33), (5.24) and the definition of $I(T)$, yields
$$
\int_{0}^{1}\frac{ds}{s}\int_{Y_{t}\cap W}\left\{K_{t}(s,x,\iota x)-
\frac{a_{1}(x,t)}{s}-a_{0}(x,t)\right\}dv_{Y_{t}}
=I(|t|^{-1})+O\left(\int_{0}^{1}\frac{ds}{s}e^{-\frac{c}{s}}\right).
\tag 6.34
$$
In the same way as $\S 5.2$, we get (5.27) which, together 
with (6.33-34) and Proposition 6.1, yields the assertion.\qed
\enddemo

\proclaim{Theorem 6.1}
As $t\to 0$, one has the following asymptotic formula:
$$
\log\tau_{M}(Y_{t},\iota_{t})=-\frac{1}{8}\,\log|t|+O(1).
$$
\endproclaim

\demo{Proof}
Let $\omega_{Z_{t}}:=\omega_{Y_{t}}|_{Z_{t}}$ be the induced metric on
$Z_{t}$. In the similar way as [B-B, Th\'eor\`eme 6.2], there exists an
universal constant $\beta\in\Bbb R$ such that one has
$\log\tau(Z_{t},\omega_{Z_{t}})=\beta\,\log|t|+O(\log\log|t|^{-1})$
as $t\to0$, which, together with Lemma 6.8 and Theorem 5.1, yields
$$
\log\tau_{M}(Y_{t},\iota_{t})=\tilde{I}(|t|^{-1})+O(\log\log|t|^{-1}),
\quad|\tilde{I}(T)|\leq C\,(\log T+1)
\tag 6.35
$$
where $\tilde{I}(T)=I(T)+\beta\log T$.
Let $\eta(t)$ be a relative canonical form such that $\eta(t)\not=0$. 
Let $\omega_{1}(t),\cdots,\omega_{g}(t)$ be a basis of
$\pi_{*}\Omega^{1}_{Y/D}$. We may assume that $\omega_{1}(0)$ has at most 
logarithmic pole at $o$ and $\omega_{2}(0),\cdots,\omega_{g}(0)$ are 
holomorphic as $Z_{0}$ has only one node. By Theorem 3.5, one has
the following on $D\backslash\{0\}$;
$$
\frac{i}{2\pi}\bar{\partial}\partial\left\{\log\tau_{M}(Y_{t},\iota_{t})+
\frac{r(M)-6}{8}\log\|\eta(t)\|^{2}+
\frac{1}{2}\log\det
\left(\int_{Z_{t}}\omega_{i}(t)\wedge\bar{\omega}_{j}(t)\right)\right\}=0.
\tag 6.36
$$
Since $\log\|\eta(t)\|^{2}=O(1)$ and 
$\log\det\left(\int_{Z_{t}}\omega_{i}(t)\wedge\bar{\omega}_{j}(t)\right)=
O(\log\log|t|^{-1})$ as $t\to 0$, by the same argument as 
[B-B, Proposition 10.1] together with (6.35-36), there exits 
$\alpha\in\Bbb R$ such that
$$
\log\tau_{M}(Y_{t},\iota_{t})=\alpha\,\log|t|^{2}+O(1)\quad(t\to 0).
\tag 6.37
$$
By (6.35), $\alpha=\lim_{T\to\infty}\tilde{I}(T)/\log T$ is an invariant 
of the instanton $(X,\omega_{X})$.
\par
To determine $\alpha$, let us compute an example. 
Take $M=\II_{1,1}\oplus E_{8}(-2)$. Let $\Phi$ be Borcherds's 
$\Phi$-function of weight 12 over $\Omega$, the Hermitian domain of 
type IV associated to the lattice 
$\II_{1,1}\oplus\II_{1,1}\oplus\Lambda_{24}(-1)$, with zero divisor 
$\Cal D$ (the discriminant locus) where $\Lambda_{24}$ is 
the 24-dimensional Leech lattice. (For $\Phi$, see [B2,3].). Put
$\Psi_{M^{\perp}}:=\Phi|_{\Omega_{M}}$ for the restriction. As the 
orthogonal compliment of $E_{8}(2)$ in $\Lambda_{24}$ is the 
16-dimensional Barnes-Wall lattice $\Lambda_{16}$ ([B2], [C-S]) which is 
free from roots, $\Psi_{M}$ is a nonzero holomorphic modular from of 
weight 12. Let 
$d\in\Delta(\II_{1,1}\oplus\II_{1,1}\oplus\Lambda_{24}(-1))$, and 
$\pi(d)\in\Lambda_{16}^{\lor}(-1)$ be the orthogonal projection. Then,
$H_{d}\cap\Omega_{M}\not=\emptyset$ if and only if 
$-2<\langle\pi(d),\pi(d)\rangle\leq 0$. Computing 
$\theta_{\Lambda_{16}^{\lor}}(\tau)$ (cf. (8.19)), we find that there is 
no norm $1$ or $3/2$ element in $\Lambda_{16}^{\lor}$, and thus 
$\pi(d)=0$. Namely, the zero divisor of $\Psi_{M^{\perp}}$ coincides 
with $\Cal D_{M}$. Since $j_{M}$ takes its value in $S^{2}(\Cal A_{1})$ by 
Theorem 2.5 and $j_{M}^{*}(\Delta_{1}\Delta_{2})$ vanishes of order $2$ 
along $\Cal D_{M}$ (because the quotient map $\Omega_{M}\to\Cal M_{M}$
blanches of order $2$ along $\Cal D_{M}$) where $\Delta_{i}$ is the
Jacobi $\Delta$-function in the $i$-th variable, 
$\Delta':=\Psi_{M}\otimes j_{M}^{*}(\Delta_{1}\Delta_{2})$ vanishes of 
order $3$ along $\Cal D_{M}$. As $\Cal D_{M}$ is the divisor of 
type $(2,1)$ by Theorem 2.5, it follows from Theorem 3.3 and (3.37)
together with  [B-B, Proposition 10.1] that
$$
\frac{i}{2\pi i}\bar{\partial}\partial\log
\left[\tau_{M}\,\|\Delta'\|^{\frac{1}{12}}\right]=
-\left(\alpha+\frac{1}{8}\right)\delta_{\Cal D_{M}}
\tag 6.38
$$
outside of subvarieties of codimension$\geq2$ where $\delta_{\Cal D_{M}}$
is the current $\int_{\Cal D_{M}}$. Since $\tau_{M}$ and $\|\Delta'\|$ 
are $\Gamma_{M}$-invariant, (6.38) can be regarded as an equation of 
currents over $\Cal M_{M}$ by Hartogus's extension theorem. From the 
residue theorem, it follows that $\alpha=-1/8$.\qed
\enddemo

\proclaim{Corollary 6.1}
If $\Lambda=\II_{1,1}\oplus E_{8}(-2)$, there exists a constant 
$C_{\Lambda}\not=0$ such that 
$\tau_{\Lambda}=\|\Delta_{\Lambda}\|^{-1/4}$ where
$\varDelta_{\Lambda}=C_{\Lambda}\,\Phi|_{\Omega_{\Lambda}}\otimes
j_{\Lambda}^{*}(\Delta_{1}\Delta_{2})$ and $\Phi$ is the denominator
function of the fake monster Lie algebra.
\endproclaim

\beginsection
$\S7$. Identification of $\pmb{\tau_{M}}$ with an Automorphic From

\par
Let $\Cal A_{g}$ be the Siegel modular variety and $\Cal A_{g}^{*}$ be the 
Satake compactification. Let
$j_{M}:\Omega_{M}\dashrightarrow\Cal A_{g(M)}$ be the rational map 
as in Proposition 2.4, and
$\hat{\Omega}_{M}=\overline{1\times j_{M}(\Omega_{M})}$ be the closure 
of the graph of $j_{M}$ in $\Omega_{M}\times\Cal A_{g(M)}$. 
Let $p_{1}:\hat{\Omega}_{M}\to\Omega_{M}$ and
$p_{2}:\hat{\Omega}_{M}\to\Cal A_{g}$ be the morphisms induced by
the projections. We regard $j_{M}^{*}\omega_{\Cal A_{g}}$ (which is 
originally defined on $\Omega_{M}^{0}$) as a current on 
$\Omega_{M}^{0}\cup\Cal D_{M}$ by putting
$j_{M}^{*}\omega_{\Cal A_{g}}=p_{1*}p_{2}^{*}\omega_{\Cal A_{g}}$ where
$\Cal D_{M}^{0}:=\bigcup_{\delta\in\Delta(N)}H_{\delta}^{0}$
($H_{\delta}^{0}=H_{\delta}\backslash\bigcup_{d\not=\pm\delta}H_{d}$).

\proclaim{Theorem 7.1}
One has the following equation of currents on $\Omega_{M}^{0}\cup
\Cal D_{M}^{0}$;
$$
\frac{i}{2\pi}\bar{\partial}\partial\log\tau_{M}=
\frac{1}{8}\delta_{\Cal D_{M}}-\frac{r(M)-6}{8}\omega_{M}-
\frac{1}{2}j_{M}^{*}\omega_{\Cal A_{g(M)}}.
$$
\endproclaim

\demo{Proof}
Theorems 3.5, 5.2, 6.1, together with Bismut-Bost's extension argument
([B-B, Proposition 10.2]), yields the assertion.\qed
\enddemo

Let $\Cal F_{g}$ ($g=g(M)$) be the sheaf of Siegel modular forms of 
weight $1$ over $\Cal A_{g}$. Set $\lambda_{M}:=
i_{*}\Cal O_{\Omega_{M}}((j_{M}^{0})^{*}\Cal F_{g})$ 
where $i:\Omega_{M}^{0}\cup\Cal D_{M}^{0}\hookrightarrow \Omega_{M}$ 
is the inclusion. It is an invertible sheaf on $\Omega_{M}$ because 
$j_{M}^{0}:\Omega_{M}^{0}\cup\Cal D_{M}^{0}\to\Cal A_{g}$ is regular and 
is defined out side of subvarieties of codimension $2$. Since $j_{M}^{0}$ 
is independ of markings, $\lambda_{M}$ is invariant under $\Gamma_{M}$. 
Namely, it is a $\Gamma_{M}$-module. Let $\chi:\Gamma_{M}\to\Bbb C^{*}$ be 
a character.

\proclaim{Definition 7.1}
$f\in H^{0}(\Omega_{M}, \lambda_{M}^{\otimes q})$ 
is said to be an automorphic form of weight $(p,q)$ with character $\chi$ if
$$
f(\gamma\cdot z)=\chi(\gamma)\,j(\gamma,z)^{p}\gamma^{*}f(z)
$$
for any $z\in\Omega_{M}$ and $\gamma\in\Gamma_{M}$ where
$j(\gamma,z):=\langle\gamma\cdot z,l_{M}\rangle/\langle z,l_{M}\rangle$
is an automorphic factor. The Petersson norm, $\|f\|$, is defined by
$$
\|f(z)\|^{2}:=K_{M}(z,\bar{z})^{p}\det\Im(j_{M}(z))^{q}|f(z)|^{2}.
$$
Here, $K_{M}(z,\bar{z})$ is the Bergman kernel and $l_{M}\in N_{\Bbb C}$ 
is the same vector as in $(3.4)$.
\endproclaim

In the sequel, we often omit character. Thus, an automorphic form is 
rigorously speaking an automorphic form with some character.
Since $\log j_{M}^{*}\|\psi\|$ becomes a locally integrable function on 
$\Omega_{M}$ for any meromorphic Siegel modular form $\psi$, the curvature 
current of $(\lambda_{M},\|\cdot\|)$ ($\|\cdot\|$ is the Petersson norm) 
can be defined in the usual manner, and coincides with
$c_{1}(\lambda_{M},\|\cdot\|)=j_{M}^{*}\omega_{\Cal A_{g}}$ on
$\Omega_{M}^{0}\cup\Cal D_{M}^{0}$.

\proclaim{Theorem 7.2}
Suppose $r(M)\leq 17$. Then, there exists a modular form $\varDelta_{M}$ 
of weight $(r(M)-6,4)$ such that $\tau_{M}=\|\varDelta_{M}\|^{-1/4}$ and
$\div(\varDelta_{M})=\Cal D_{M}$.
\endproclaim

\demo{Proof}
Take a non-zero meromorphic modular form $\phi$ of weight $(r-6,4)$ 
such that $\Cal D_{M}$ is not contained in the zero and polar 
locus of $\phi$. Put $F:=\tau_{M}^{8}\,\|\phi\|^{2}$. By Theorem 7.1 and
Hartogus's theorem, one has
$\frac{i}{2\pi}\bar{\partial}\partial\log\,F=\delta_{\Cal D_{M}}-
\delta_{\roman{div}(\phi)}$ on $\Omega_{M}$.
Therefore, $\partial\log\,F$ is a $\Gamma_{M}$-invariant meromorphic 
$1$-form on $\Omega_{M}$ with at most logarithmic poles, and thus
$G(y):=\exp\left(\int_{*}^{y}\partial\log\,F\right)$ 
($*$ is a reference point in $\Omega_{M}$) is a meromorphic function on 
$\Omega_{M}$ such that $\div(G)=\Cal D_{M}-\div(\phi)$. 
Let $\gamma\in\Gamma_{M}$. Let $[\gamma]$ be a simple closed real
curve in $\Cal M_{M}=\Omega_{M}/\Gamma_{M}$ corresponding to a path 
joining $y$ and $\gamma\cdot y$. As $\Omega_{M}$ is diffeomorphic to the 
cell, the homotopy class of $[\gamma]$ does not depend on a choice of 
$y$. Thus, $\chi(\gamma)=\exp\left(\int_{[\gamma]}\partial\log\,F\right)
=G(\gamma\cdot y)/G(y)$ is independent of $y$ and becomes a character of 
$\Gamma_{M}$. Since $\Gamma_{M}/[\Gamma_{M},\Gamma_{M}]$ is finite by 
Kazhdan's theorem ([Kz]), $\chi$ takes its values in $S^{1}=U(\Bbb C)$. 
Therefore, $\log\,|G|^{2}$ is a $\Gamma_{M}$-invariant 
pluriharmonic function whose divisor is the same as $\log F$. Since the 
Satake-Baily-Borel boundary of $\Cal M_{M}$ has $\codim\geq2$ when 
$r(M)\leq 17$, Hartogus's theorem implies that there exists a constant 
$C\not=0$ such that $F=C^{2}\,|G|^{2}$. Thus, we get 
$\tau_{M}^{-8}=C^{2}\,|G|^{2}\,\|\phi\|^{2}$, and 
$\varDelta_{M}:=C\,G\cdot\phi$ is the desired form.\qed
\enddemo

Let $\delta\in\Delta(N)$. If $\langle M\oplus\delta\rangle$ denotes the 
smallest $2$-elementary lattice generated by $M,\delta$, then
$\Omega_{\langle M\oplus\delta\rangle}$ can naturally be identified with
$H_{\delta}$. As $\varDelta_{M}$ is a section of 
$\lambda_{M}^{\otimes4}$ which vanishes of order one along $H_{\delta}$, 
it follows from Proposition 2.5 that
$\varDelta_{M}\cdot\langle z,l_{M}\rangle/\langle z,\delta\rangle$ 
restricted to $H_{\delta}^{0}$ is a section of
$\lambda_{\langle M\oplus\delta\rangle}^{\otimes4}$. Note that if 
$g(\langle M\oplus\delta\rangle)=g(M)-1$, $\Cal A_{g-1}$ is considered to 
be one of the boundary components of $\Cal A_{g}^{*}$ and the restriction
map $S:\Cal F_{g}^{\otimes 4}\to\Cal F_{g-1}^{\otimes 4}$ coincides with 
the Siegel operator ([F]).

\proclaim{Theorem 7.3}
Under the identification $\Omega_{\langle M\oplus\delta\rangle}=
H_{\delta}$, one has
$$
\varDelta_{\langle M\oplus\delta\rangle}(y)=
C_{\langle M\oplus\delta\rangle}\,
\lim_{z\to y}\frac{\langle z,l_{M}\rangle}{\langle z,\delta\rangle}\,
\varDelta_{M}(z)=C(M,\delta)\,\left.\frac{\langle\cdot,l_{M}\rangle}
{\langle\cdot,\delta\rangle}\,\varDelta_{M}\right|_{H_{\delta}}(y).
$$
Here $C_{\langle M\oplus\delta\rangle}$ is a nonzero constant.
\endproclaim

\demo{Proof}
Put $\delta:=\delta_{0}$. We separate the proof into two cases.
\newline{\bf{Case (1)}}
Assume $\langle\delta_{0},\delta_{1}\rangle=-1$ for some 
$\delta_{1}\in\Delta(N)$. With this condition, $M\oplus\Bbb Z\delta_{0}$ 
is primitive in $L_{K3}$ and
$\langle M\oplus\delta_{0}\rangle=M\oplus\Bbb Z\delta_{0}$ which implies 
$g(\langle M\oplus\delta_{0}\rangle)=g(M)-1$. By the condition, 
$\Bbb Z\delta_{0}\oplus\Bbb Z\delta_{1}=A_{2}(-1)$ where $A_{2}$ is the
$A_{2}$-root lattice whose roots are
$\pm\delta_{0},\pm\delta_{1},\pm\delta_{2}$ where
$\delta_{2}:=\delta_{0}-\delta_{1}$. 
\par{\it Step $1$.}
Take $\omega_{0}\in(H_{\delta_{0}}\cap H_{\delta_{1}})^{0}:=H_{\delta_{0}}
\cap H_{\delta_{1}}\backslash\cup_{d\not=\pm\delta_{i}}H_{d'}$. 
Choose $\kappa_{0}\in C(M)$ such that 
$\langle\kappa_{0},d\rangle\not=0$ for any $d\in\Delta(M)$.
By the surjectivity of the period map, there exists a marked $K3$ surface
$(X_{0},\phi)$ with nef and big line bundle $L_{0}$ such that
$\pi(X_{0},\phi)=\omega_{0}$ and $\phi(c_{1}(L_{0}))=\kappa_{0}$. 
Choosing $0<\epsilon\ll1$, we may suppose that there exists a K\"ahler 
class whose image by $\phi$ is $\kappa:=\kappa_{0}-\epsilon\delta_{0}$. 
Let $C_{\delta_{i}}$ be the cycle corresponding to $\delta_{i}$.
As $\langle\kappa,\delta_{i}\rangle>0$, any $C_{\delta_{i}}$ is effective. 
We choose $\epsilon$ so small that $C_{\delta_{1}},C_{\delta_{2}}$ are 
irreducible $-2$-curves and 
$C_{\delta_{0}}=C_{\delta_{1}}\cup C_{\delta_{2}}$.
\par
Let $U$ be a small neighborhood of $(X_{0},\phi)$ in $\tilde{\Omega}_{M}$, 
and $(X,\phi)\to U$ be the universal family. By construction, there 
exists a relatively nef and big line bundle $L\to X$ which restricted 
to $X_{0}$ is $L_{0}$. Choosing $m\gg1$, let $\Phi_{|mL|}:X\to\Bbb P^{N}$ 
be the morphism associated to the complete linear system $|mL|$, 
and $Y:=\Phi_{|mL|}(X)\to U$ be the image. As $\kappa_{0}$ does not 
intersect $\delta_{i}$, the cycle $C_{\delta_{i}}$ corresponding to 
$\delta_{i}$ is an irreducible $-2$-curve in $X_{t}$ if 
$t\in H_{\delta_{i}}^{0}$ (cf. Lemma 2.3). By above construction, 
$C_{\delta_{0}}$ splits into two components if 
$t\in(H_{\delta_{0}}\cap H_{\delta_{1}})^{0}$; 
$C_{\delta_{0}}=C_{\delta_{1}}\cup C_{\delta_{2}}$. By Mayer's theorem, 
$\Phi_{|mL|}:X_{t}\to Y_{t}$ is the minimal resolution whose exceptional 
locus is $C_{\delta_{i}}$ if $t\in H_{i}^{0}$ $(i=0,1,2)$,
$C_{\delta_{0}}=C_{\delta_{1}}\cup C_{\delta_{2}}$ if
$t\in(H_{\delta_{0}}\cap H_{\delta_{1}})^{0}$, and empty if outside of
the discriminant locus. Put $o\in Y_{0}$ for the image of $C_{\delta_{0}}$. 
As $\Bbb Z\delta_{1}\oplus\Bbb Z\delta_{2}$ is the $A_{2}$-root system, 
$(Y_{0},o)$ is a $K3$ surface with a $A_{2}$-singularity. Set 
$$
\Cal M_{A_{2}}:=\{\alpha\in\Bbb C^{3};\,\sum\alpha_{i}=0\},\quad
Z:=\{(x,\alpha)\in\Bbb C^{3}\times\Cal M_{A_{2}};\,
x_{1}x_{2}-\prod(x_{3}+\alpha_{i})=0\}.
\tag 7.1
$$ 
Consider the deformation of $A_{2}$-singularity 
$(Z,0)\to(\Cal M_{A_{2}},0)$ on which acts $S_{3}=W(A_{2})$ by the
permutation of coordinates $\alpha$. Then, 
$(Z/S_{3},0)\to(\Cal M_{A_{2}}/S_{3},0)$ is the semiuniversal 
deformation of $A_{2}$-singularity. The discriminant locus of 
$Z\to\Cal M_{A_{2}}$ is $D=D_{0}\cup D_{1}\cup D_{2}$ where
$D_{i}=\{\alpha\in\Cal M_{A_{2}};\,\alpha_{k}-\alpha_{j}=0\}$ 
$(\{i,j,k\}=\{0,1,2\})$. By the versality, there exist maps
$F:(Y,o)\to(Z,0)$ and $f:(U,0)\to(\Cal M_{A_{2}},0)$ which commute
with the projections such that $(Y,o)=F^{*}(Z,0)$ and 
$f(H_{\delta_{i}})=D_{i}$. 
\par
Let $\phi'=I_{M}\circ\phi$ (cf. (2.6)) be another marking, and
$(X',\phi')\to U'$ be the universal family such that $U'=I\circ U$ 
(cf. (2.8)). Let $L'\to X'$ be the relatively nef and big line bundle such 
that $\phi'(c_{1}(L'_{s}))=\kappa_{0}$. By the similar construction as 
before, we get a family $Y'\to U'$ such that $X'\to Y'$ is the 
simultaneous resolution. Since $\phi^{-1}(\kappa_{0})$ and 
$(\phi')^{-1}(\kappa_{0})$ are weak polarizations ([Mo, pp.318]) of 
$(X_{t},\phi)$ and $(X'_{t},\phi')$ respectively, we can define the weakly
polarized period map ([Mo, pp.318]) $U\to V$ and $U'\to V$ by sending 
$([\omega_{t}],\kappa_{t})$ to $([\omega_{t}],\kappa_{0})$ where 
$[\omega_{t}]$ is the period and $\kappa_{t}$ is the K\"ahler class.
Note that the weak polarized period domain in our situation is
$\Omega_{M}\times C(M)^{+}$ and $V$ is its subset. As $U\to V$ 
(resp. $U'\to V$) is an isomorphism, we may regard $Y$ and $Y'$
are families over $V$. Since there exists the universal marked family of
generalized $K3$ surfaces over $\Omega_{M}\times C(M)^{+}$ ([Mo, pp.321]), 
we get an identification $e:Y'\cong Y$. Let $p:X\to Y$ and $p':X'\to Y'$
be the simultaneous resolution, and $\iota_{M}:X\to X'$ be the isomorphism
as in (2.8). Then, $\iota:=e\circ p'\circ\iota_{M}\circ p^{-1}$ is a
rational automorphism over $Y$. By the weakly polarized global Torelli
theorem ([Mo, pp.319]), $\iota_{t}$ is an anti-symplectic involution on
$Y_{t}$ for any $t\in V$ and $\iota$ is holomorphic everwhere. By 
appropriate normalizations, we may suppose $(Y,o)=(Z,0)$ and
$\iota(x_{0},x_{1},x_{2})=(x_{1},x_{0},x_{2})$ in (7.1). 
\par{\it Step $2$.}
Let $C\to V$ be the family of fixed curves. By Theorem 2.5, we get the
decomposition $C=C^{(g)}+\sum E_{i}$ where $C^{(g)}_{t}$ is a smooth 
irreducible curve of genus $g=g(M)$ for generic $t\in V$ and $E_{i}$ is a
family of smooth rational curves. As $C^{(g)}_{0}$ is a fixed locus, 
$\Sing\,C^{(g)}_{0}=\Sing\,Y_{0}=o$. Moreover it is a Cartier divisor 
because
$(C^{(g)}_{0},o)=(Z_{0},0)\cap\{x_{0}-x_{1}=0\}=\{(x_{2},x_{3});\,
x_{2}^{2}-x_{3}^{3}=0\}$. As $C^{(g)}_{0}\cdot C^{(g)}_{0}=2(g-1)\geq0$, 
$C^{(g)}_{0}$ is a nef and effective divisor. Since $(Y_{0},o)$ is an 
$A_{2}$-singularity, $o$ can not belong to the fixed part of the linear 
system $|C^{(g)}_{0}|$ by Saint-Donat's theorem. Therefore, we can pick up 
$f_{1},\cdots,f_{g}$, nonconstant meromorphic functions on $Y$, such that 
$\div(f_{i})=P_{i}-C^{(g)}$ with $o\not\in P_{g}$ but 
$\Sing\,Y_{t}\in P_{i}$ for $1\leq i<g$. Let $\omega_{Y/V}$ be a nowhere 
vanishing relative $2$-form on $Y$. Let $\varphi_{i}(t)$ be the $1$-form 
on $C^{(g)}_{t}$ defined by
$\varphi_{i}(t):=\Res_{C^{(g)}_{t}}f_{i}\omega_{Y/V}$. As $\Sing\,Y_{t}$ 
is contained in the zero locus of $f_{j}$ for $j<g$, $\varphi_{j}(t)$ 
$(j<g)$ becomes a holomorphic $1$-form on the normalization of 
$C^{(g)}_{t}$ when $t\in H_{\delta_{0}}$. Thus, if $j<g$,
$$
\int_{C_{t}^{(g)}}\varphi_{i}(t)\wedge\overline{\varphi_{j}(t)}\to
\int_{C_{0}^{(g)}}\varphi_{i}(0)\wedge\overline{\varphi_{j}(0)}<\infty
\quad(t\to H_{\delta_{0}}).
\tag 7.2
$$
Let $A_{1},B_{1},\cdots,A_{g},B_{g}$ be the symplectic basis of 
$H^{1}(C_{t}^{(g)},\Bbb Z)$. As $(C^{(g)}_{0},0)$ is the cusp, we may 
assume that $A_{g}$ and $B_{g}$ are vanishing cycles, and all other cycles    
converges to the symplectic basis of the normalization of $C^{(g)}_{0}$.
Since we have fixed the symplectic basis, $j_{M}$ takes its value in
$\frak S_{g}$:
$$
\tau(j_{M}(t))=\left(\int_{A_{i}}\varphi_{j}(t)\right)^{-1}
\left(\int_{B_{i}}\varphi_{j}(t)\right)\in\frak S_{g}.
\tag 7.3
$$
As the bidegree $(1,1)$ part of the diagonal in 
$C^{(g)}_{t}\times C^{(g)}_{t}$ is homologous to
$\sum A_{i}^{(1)}B_{i}^{(2)}-B_{i}^{(1)}A_{i}^{(2)}$ 
in $H^{2}(C^{(g)}_{t}\times C^{(g)}_{t};\Bbb Z)$, it follows from (7.2) 
that
$$
\int_{C^{(g)}_{t}}\varphi_{g}(t)\overline{\varphi_{g}(t)}=
\Im\tau_{gg}(j_{M}(t))\,\left|\int_{A_{g}}\varphi_{g}(t)\right|^{2}+O(1).
\tag 7.4
$$
Let $E_{\beta}:=\{(x_{2},x_{3});\,
x_{2}^{2}=(x_{3}+1)(x_{3}+\beta)(x_{3}-1-\beta)\}$ be an elliptic curve
and $\gamma_{1},\gamma_{2}$ be the symplectic basis such that 
$\gamma_{1}$ converges to the cycle $|x_{3}+1|=\epsilon$ as
$\beta\to0$. Since $o\not\in P_{g}$, multiplying a constant if necessary, 
we may assume $f_{g}=1/\{x_{2}^{2}-\prod(x_{3}+\alpha_{i})\}+O(1)$ 
on a neighborhood of $o$. When $t\to H_{\delta_{0}}$, 
$\alpha(t)\to D_{0}$ and $\alpha_{2}(t)\to\alpha_{1}(t)$. 
Putting $\beta=\alpha_{1}/\alpha_{2}$, we may suppose $A_{g}$ is 
identified with $\gamma_{1}$ and get
$$
\aligned
\int_{A_{g}}\varphi_{g}(t)
&=\int_{A_{1}}\frac{dz}
{\sqrt{(z+\alpha_{0})(z+\alpha_{1})(z+\alpha_{2})}}+O(1)\\
&=\alpha_{2}^{-\frac{1}{2}}\int_{\gamma_{1}}\frac{dz}
{\sqrt{(z+1)(z+\beta)(z-1-\beta)}}+O(1)
=\frac{2\pi}{\sqrt{3}}\alpha_{2}^{-\frac{1}{2}}+O(1).
\endaligned
\tag 7.5
$$
Let $W$ be a small neighborhood of 
$\omega_{0}\in(H_{\delta_{0}}\cap H_{\delta_{1}})^{0}$ in 
$H_{\delta_{0}}^{0}$. Since 
$$
\left.
\left[\frac{\det\left(
\int_{C^{(g)}}\varphi_{i}\overline{\varphi_{j}}\right)}
{\int_{C^{(g)}}\varphi_{g}\overline{\varphi_{g}}}\right]\right|_{W}=
\left.\det\left(
\int_{C^{(g)}}\varphi_{i}\overline{\varphi_{j}}\right)_{i,j<g}\right|_{W}
\tag 7.6
$$ 
by (7.2-5), it follows that
$$
\aligned
\,&\frac{i}{2\pi}\bar{\partial}\partial\left[\log\left.
\frac{\det\left(\int_{C^{(g)}}\varphi_{i}\overline{\varphi_{j}}\right)}
{j_{M}^{*}\Im\tau_{gg}}\right|_{W}\right]\\
&=\frac{i}{2\pi}\bar{\partial}\partial\left[\left.
\log\left|\int_{A_{g}}\varphi_{g}\right|^{2}\right|_{W}\right]
+\frac{i}{2\pi}\bar{\partial}\partial\left[\left.\log\det
\left(\int_{C^{(g)}}\varphi_{i}\overline{\varphi_{j}}\right)_{i,j<g}
\right|_{W}\right]\\
&=\frac{1}{2}\delta_{H_{\delta_{0}}\cap H_{\delta_{1}}\cap W}+
j_{\langle M\oplus\delta_{0}\rangle}^{*}\omega_{\frak S_{g-1}}|_{W}.
\endaligned
\tag 7.7
$$
Since $H_{d}\cap H_{\delta_{0}}\cap W\not=\emptyset$ $(d\in\Delta(N))$
iff $d=\pm\delta_{0},\pm\delta_{1},\pm\delta_{2}$, it follows from 
Theorem 7.1 that
$$
\aligned
\,&
\frac{i}{2\pi}\bar{\partial}\partial\left.\log\left[\tau_{M}\,\det\left(
\int_{C^{(g)}}\varphi_{i}\overline{\varphi_{j}}\right)^{\frac{1}{2}}
\left(\frac{|\langle w,\delta_{0}\rangle|^{2}}{\langle w,\bar{w}\rangle}
\right)^{\frac{1}{8}}\right|_{W}\right]\\
&=\frac{1}{8}\sum_{d\in\Delta(N)\backslash\{\pm\delta_{0}\}/\pm1}
\delta_{H_{d}\cap H_{\delta_{0}}\cap W}-
\frac{r(M)-5}{8}\omega_{\langle M\oplus\delta_{0}\rangle}\\
&=\frac{1}{4}\delta_{H_{\delta_{0}}\cap H_{\delta_{1}}\cap W}-
\frac{r(\langle M\oplus\delta_{0}\rangle)-6}{8}
\omega_{\langle M\oplus\delta_{0}\rangle}
\endaligned
\tag 7.8
$$
which, together with (7.7), yields
$$
\aligned
\,&
\frac{i}{2\pi}\bar{\partial}\partial\left.\log\left[\tau_{M}\,
(j_{M}^{*}\Im\tau_{gg})^{\frac{1}{2}}
\left(\frac{|\langle w,\delta_{0}\rangle|^{2}}{\langle w,\bar{w}\rangle}
\right)^{\frac{1}{8}}\right|_{W}\right]\\
&=-\frac{r(\langle M\oplus\delta_{0}\rangle)-6}{8}
\omega_{\langle M\oplus\delta_{0}\rangle}
-\frac{1}{2}j_{\langle M\oplus\delta_{0}\rangle}^{*}
\omega_{\frak S_{g-1}}|_{W}.
\endaligned
\tag 7.9
$$
\par{\it Step $3$.}
Let $\delta_{3}\in\Delta(N)$ such that 
$\langle\delta_{3},\delta_{0}\rangle=0$. 
Take $\omega_{0}\in(H_{\delta_{0}}\cap H_{\delta_{3}})^{0}$ and choose 
$\kappa_{0}\in C_{M}$ as in Step 1. Similarly as before, there exists 
a marked $K3$ surface $(X_{0},\phi)$ with nef and big line bundle $L_{0}$ 
such that $\pi(X_{0},\phi)=\omega_0$ and $\phi(c_{1}(L_{0}))=\kappa_{0}$. 
Let $C_{\delta_{i}}$ be the cycle corresponding to $\delta_{i}$. Then, 
$C_{\delta_{0}}$ and $C_{\delta_{4}}$ are mutually disjoint irreducible 
$-2$-curves. Let $U$ be a small neighborhood of $(X_{0},\phi)$ in 
$\tilde{\Omega}_{M}$, $(X,\phi)\to U$ the universal family, and $L\to X$ 
the relatively nef and big line bundle whose restriction to $X_{0}$ is 
$L_{0}$. Using $|mL|$, set $Y:=\Phi_{|mL|}(X)\to U$ for $m\gg1$. As 
$\langle\kappa_{0},\delta_{i}\rangle=0$, the cycle $C_{\delta_{i}}$ 
corresponding to $\delta_{i}$ is an irreducible $-2$-curve in $X_{t}$ if 
$t\in H_{\delta_{i}}^{0}$. Thus, $\Phi_{|mL|}:X_{t}\to Y_{t}$ is the 
minimal resolution whose exceptional locus is $C_{\delta_{i}}$ if 
$t\in H_{i}^{0}$ $(i=0,3)$, $C_{\delta_{0}}\cup C_{\delta_{3}}$ if
$t\in(H_{\delta_{0}}\cap H_{\delta_{3}})^{0}$, and empty if outside of the
discriminant locus. Put $o_{i}$ for the image of $C_{\delta_{i}}$. 
By construction, $(Y_{0},o_{0},o_{3})$ is a $K3$ surface with two
$A_{1}$-singularities. Let $\iota:Y\to Y$ be the anti-symplectic involution
constructed as  before, and $C^{(g)}\to U$ be the family of fixed curves of
maximal genus. By construction, $C^{(g)}\to U$ is f.s.o. in the sense of 
[B-B] (though the total space admits two nodes here). Since $C^{(g)}_{t}$
has a node for $t\in H_{\delta_{0}}^{0}$ (because 
$g(\langle M\oplus\delta_{0}\rangle)=g(M)-1$), we find 
$o_{0}\in C^{(g)}_{0}$. Let $\varphi_{1}(t),\cdots,\varphi_{g}(t)$ be a 
basis of relative $1$-forms for $C\to U$ such that $\varphi_{g}(t)$ has a 
logarithmic pole at $\Sing\,C^{(g)}_{t}$ if $t\in H_{\delta_{0}}$ and that 
$\varphi_{1}(t),\cdots,\varphi_{g-2}(t)$ are regular on the normalization 
of $C_{t}^{(g)}$ for any $t\in U$. If $o_{3}\in C^{(g)}_{0}$, we may 
suppose that $\varphi_{g-1}(0)$ has a logarithmic pole at $o_{3}$ and
$\varphi_{g-1}(t)$ is regular on the normalization of $C^{(g)}_{0}$ for
$t\in H_{\delta_{0}}^{0}$. If $o_{3}\not\in C^{(g)}_{0}$, we may assume
that $\varphi_{g-1}(t)$ is regular on the normalization of $C^{(g)}_{t}$
for any $t\in U$. Let $W$ be a small neighborhood of
$\omega_{0}\in(H_{\delta_{0}}\cap H_{\delta_{3}})^{0}$ in
$H_{\delta_{0}}^{0}$. From [B-B, Propositon 13.3] together with 
Proposition 2.5, it follows that
$$
\aligned
\frac{i}{2\pi}\bar{\partial}\partial\left[\log\left.
\frac{\det\left(\int_{C^{(g)}}\varphi_{i}\overline{\varphi_{j}}\right)}
{j_{M}^{*}\Im\tau_{gg}}\right|_{W}\right]
&=\frac{i}{2\pi}\bar{\partial}\partial\left[\left.\log\det
\left(\int_{C^{(g)}}\varphi_{i}\overline{\varphi_{j}}\right)_{i,j<g}
\right|_{W}\right]\\
&=j_{\langle M\oplus\delta_{0}\rangle}^{*}\omega_{\frak S_{g-1}}|_{W}.
\endaligned
\tag 7.10
$$
Since $H_{d}\cap H_{\delta_{0}}\cap W\not=\emptyset$ $(d\in\Delta(N))$
iff $d=\pm\delta_{0},\pm\delta_{3}$, comparing (7.8) and (7.10), we get 
$$
\aligned
\,&
\frac{i}{2\pi}\bar{\partial}\partial\left.\log\left[\tau_{M}\,
(j_{M}^{*}\Im\tau_{gg})^{\frac{1}{2}}
\left(\frac{|\langle w,\delta_{0}\rangle|^{2}}{\langle w,\bar{w}\rangle}
\right)^{\frac{1}{8}}\right|_{W}\right]\\
&=\frac{1}{8}\delta_{H_{\delta_{0}}\cap H_{\delta_{3}}\cap W}
-\frac{r(\langle M\oplus\delta_{0}\rangle)-6}{8}
\omega_{\langle M\oplus\delta_{0}\rangle}
-\frac{1}{2}j_{\langle M\oplus\delta_{0}\rangle}^{*}
\omega_{\frak S_{g-1}}|_{W}.
\endaligned
\tag 7.11
$$
\par{\it Step $4$.}
In view of the proof of Theorem 7.2, we may write 
$\varDelta_{M}=\psi\otimes j_{M}^{*}E$ where $\psi$ is a meromorphic
modular form over $\Omega_{M}$ of weight $g(r-6)$ and $E$ is an
Eisenstein series of weight $4g$ ([F]) whose divisor does not contain 
the boundary component $\Cal A_{g}^{*}\backslash\Cal A_{g}$. For 
$\tau=(\tau_{ij})_{i,j\leq g}\in\frak S_{g}$, write
$\tau'=(\tau_{ij})_{i,j<g}\in\frak S_{g-1}$. Let 
$S:H^{0}(\Cal A_{g},\Cal F_{g}^{k})\to 
H^{0}(\Cal A_{g-1},\Cal F_{g-1}^{k})$ be the Siegel operator.
Since $j_{M}(H_{\delta_{0}})\subset\Cal A_{g-1}$, it follows from
Proposition 2.5 that
$$
\left.\left(\frac{\langle\cdot,l_{M}\rangle}
{\langle\cdot,\delta_{0}\rangle}\,
\varDelta_{M}\right)^{g}\right|_{H_{\delta_{0}}}=
\left.\left(\frac{\langle\cdot,l_{M}\rangle}
{\langle\cdot,\delta_{0}\rangle}\right)^{g}\,
\psi\right|_{H_{\delta_{0}}}\otimes 
j_{\langle M\oplus\delta_{0}\rangle}^{*}S(E).
\tag 7.12
$$
From the definition of Petersson norm, it follows that
$$
\aligned
\left\|\left.\left(
\frac{\langle\cdot,l_{M}\rangle}{\langle\cdot,\delta_{0}\rangle}\,
\varDelta_{M}\right)^{g}\right|_{H_{\delta_{0}}}\right\|^{2}(w)
&=K_{M}^{-g(r-5)}\left|\left.\left(\frac{\langle z,l_{M}\rangle}
{\langle z,\delta_{0}\rangle}\right)^{g}\,\psi\right|_{z=w}\right|^{2}\,
j_{M}^{*}(\det\Im\tau')^{4g}|S(E)|^{2}\\
&=\left.\|\varDelta_{M}\|^{2g}(j_{M}^{*}\Im\tau_{gg})^{-4g}
\left(\frac{|\langle z,\delta_{0}\rangle|^{2}}
{\langle z,\bar{z}\rangle}\right)^{-g}\right|_{z=w}\\
&=\left.\tau_{M}^{-8g}\,(j_{M}^{*}\Im\tau_{gg})^{-4g}
\left(\frac{|\langle z,\delta_{0}\rangle|^{2}}{\langle z,\bar{z}\rangle}
\right)^{-g}\right|_{z=w}
\endaligned
\tag 7.13
$$
where $K_{M}$ is the Bergman kernel as in (3.4). Comparing (7.9) and 
(7.11) with (7.13), $\|\varDelta_{M}\cdot\langle\cdot,l_{M}\rangle/
\langle\cdot,\delta_{0}\rangle|_{H_{\delta_{0}}}\|^{2}$ satisfies
the same $\bar{\partial}\partial$-equation as 
$\|\varDelta_{\langle M\oplus\delta_{0}\rangle}\|^{2}$ over 
$\Omega_{\langle M\oplus\delta_{0}\rangle}^{0}\cup
\Cal D_{\langle M\oplus\delta_{0}\rangle}^{0}$. In view of the proof of
Theorem 7.2, solution of the $\bar{\partial}\partial$-equation in 
Theorem 7.1 is unique up to constant. Thus, we get the assertion. 
\par{\bf Case (2)}
Suppose that there is no $d\in\Delta(N)$ such that 
$|\langle d,\delta_{0}\rangle|=1$. Thus, for any $d\in\Delta(N)$,
$H_{d}\cap H_{\delta_{0}}\not=\emptyset$ iff $\langle d,\delta\rangle=0$.
Now, we can prove the assertion in the same way as Step 3 and 4 in
Case (1), and details are left to the reader.\qed
\enddemo

\beginsection
$\S8$. Borcherds's Products Arizing from $\pmb{\varDelta_{M}}$

\par
\subsubhead
{\bf 8.1 Borcherds's Product}\rm
\endsubsubhead
For a lattice $M$ of signature $(2,b^-)$, let 
$\{e_{\gamma}\}_{\gamma\in M^{\lor}/M}$ be the standard unitary basis of 
the group ring $\Bbb C[M^{\lor}/M]$. We denote by $\rho_{M}$ the Weil
representation of the metaplectic  group $Mp_{2}(\Bbb Z)$;
$$
\rho_{M}(T)\,e_{\gamma}=e^{\pi i\langle\gamma,\gamma\rangle}\,e_{\gamma},
\quad
\rho_{M}(S)\,e_{\gamma}=\frac{\sqrt{i}^{b^{-}-2}}{\sqrt{|M^{\lor}/M|}}
\sum_{\delta\in M^{\lor}/M}e^{2\pi i\langle\gamma,\delta\rangle}\,
e_{\delta}
\tag 8.1
$$
where $T=(\binom{1\,1}{0\,1},1)$, $S=(\binom{0-1}{1\,\,0},\sqrt{\tau})$ 
are the generators of $Mp_{2}(\Bbb Z)$. In this section, we suppose
$M=\II_{1,1}(N)\oplus K$ ($N=1,2$ and $K$ is a hyperbolic 2-elementary 
lattice.) Let 
$F(\tau)=\sum_{\gamma\in M^{\lor}/M}e_{\gamma}
\sum_{k\in\Bbb Q}c_{\gamma}(k)\,q^k$ be a nearly holomorphic modular form 
of type $\rho_{M}$ with weight $1-b^{-}/2$ which has the integral Fourier 
coefficients at the cusp. We denote by $\Psi_{M}(z,F)$ Borcherds's
product attached to $F(\tau)$ ([Bo5]);
$$
-\log\|\Psi_{M}(z,F)\|^{4}=\int_{SL_{2}(\Bbb Z)\backslash\Bbb H}
\bar{\Theta}_{M}(\tau;z)\,F(\tau)\,y\,dx\,dy/y^2
+c_0(0)(\Gamma'(1)+2\log\sqrt{2\pi}).
\tag 8.2
$$
$\Psi_{M}(z,F)$ is a function on the Grassmannian $G(2,b^-)$ which is
isomorphic to the tube domain $K_{\Bbb R}+\sqrt{-1}C(K)^+$. The 
following theorem is due to Borcherds.

\proclaim{Theorem 8.1 ([Bo5, Theorem 13.3])}
\flushpar{$(1)$}
$\Psi_{M}(z,F)$ is an automorphic form on $G(2,b^-)$ for some arithmetic
subgroup of $O(M)$ of weight $c_{0}(0)/2$.
\flushpar{$(2)$}
The zeros or poles of $\Psi_{M}(z,F)$ lies on the divisor 
$\lambda^{\perp}$ $(\lambda\in M$, $\lambda^2<0)$ of order
$\sum_{0<x\in\Bbb R,\,x\lambda\in M^{\lor}}c_{x\lambda}(x^2\lambda^2/2)$.
\flushpar{$(3)$}
$\Psi_{M}(z,F)$ admits the following holomorphic infinite product expansion
near the cusp and $z\in K_{\Bbb R}+\sqrt{-1}W$;
$$
\Psi_{M}(z,F)=e^{2\pi i\langle\rho(K,W,F_{K}),z\rangle}
\prod_{\lambda\in K^{\lor}\cap W^{\lor}}\prod_{n\in\Bbb Z/N\Bbb Z}
(1-e^{2\pi i\langle\lambda,z+n/N\rangle})^{c_{\lambda+nf'/N}(\lambda^2/2)}
$$
where $\rho(K,W,F_{K})$ is the Weyl vector, $W$ is a Weyl chamber, 
$W^{\lor}$ the dual cone of $W$, and $f,f'$ are the generators of
$\II_{1,1}(N)$ such that $f\cdot f=f'\cdot f'=0$, $f\cdot f'=N$.
\endproclaim

\subsubhead
{\bf 8.2 $\pmb{2}$-Elementary $\pmb{K3}$ Surfaces with
$\pmb{(g,\delta)=(0,1)}$}\rm
\endsubsubhead
By Nikulin's table ([Ni4]), all the primitive 2-elementary hyperbolic 
lattices in $L_{K3}$ with $(g,\delta)=(0,1)$ are isometric to one of 
the following $S_{k}$ ($1\leq k\leq9$).
\par
Let $h,\delta_{0},\cdots,\delta_{8}$ be the basis of $\I_{1,9}(2)$
such that $h^{2}=2$, $h\cdot\delta_{i}=0$,
$\delta_{i}\cdot\delta_{j}=-2\delta_{ij}$. Then, 
$\I_{1,8}(2)=A_{1}\oplus A_{1}(-1)^{\oplus8}=
\Bbb Z\,h\oplus\Bbb Z\,\delta_{1}\oplus\cdots\oplus\Bbb Z\,\delta_{8}=
\I_{1,9}(2)\cap\delta_{0}^{\perp}$. Put
$$
\rho:=(3h-\delta_{0}-\cdots-\delta_{8})/2\in\I_{1,9}(2)^{\lor},\quad
\kappa:=3h-(\delta_{1}+\cdots+\delta_{8}).
\tag 8.3
$$
$\Lambda$ is a sublattice of $\I_{1,9}(2)^{\lor}$ defined by
$\Lambda=\Bbb Z\rho\oplus\Bbb Z\delta_{0}\oplus
\rho^{\perp}\cap\delta_{0}^{\perp}\cap\I_{1,9}(2)$. Note that
$\kappa^{\perp}\cap\I_{1,8}(2)\cong E_{8}(-2)$ (cf. [Man]).
As $\delta_{0}^{\perp}\cap\Lambda=\Bbb Z h\oplus\Bbb Z\delta_{1}\oplus
\cdots\oplus\Bbb Z\delta_{8}=\I_{1,8}(2)$, 
above $\delta_{0},\cdots,\delta_{8}$ satisfy
$$
\delta_{0}^{\perp}\cap\cdots\cap\delta_{k-1}^{\perp}\cap\Lambda
=\Bbb Z\,h\oplus\Bbb Z\,\delta_{k}\oplus\cdots\oplus\Bbb Z\,\delta_{8}
=\I_{1,9-k}(2),\quad
\delta_{i}\cdot\rho=1.
\tag 8.4
$$
Define 2-elementary lattices $\Lambda_{k}$, $S_{k}$, $T_{k}$ 
($1\leq k\leq9$) by 
$$
\Lambda_{k}:=\delta_{0}^{\perp}\cap\cdots\cap\delta_{k-1}^{\perp}=
\I_{1,9-k}(2),\quad 
T_{k}:=\II_{1,1}(2)\oplus\Lambda_{k}\cong\I_{2,10-k}(2),
\quad S_{k}:=T_{k}^{\perp}.
\tag 8.5
$$ 
Here, the first orthogonal compliment is considered in $\Lambda_{k}$, 
and the last one in $L_{K3}$. Then,
$(r(S_{k}),l(S_{k}),\delta(S_{k}))=(10+k,12-k,1)$.
Let $\pi_{k}:\Lambda\to\Lambda_{k}^{\lor}$ be the orthogonal projection. 
As is well known (cf. [Man]), the Weyl vector of $\Lambda_{k}$ is given by
$$
\rho_{k}:=\pi_{k}(\rho)=(3h-\delta_{k}-\cdots-\delta_{8})/2.
\tag 8.6
$$
We denote by $\Am(\Lambda_{k,\Bbb R})$ the Weyl chamber containing 
$\rho_{k}$ and by $\NE(\Lambda_{k,\Bbb R})$ the dual cone of 
$\Am(\Lambda_{k,\Bbb R})$ i.e., 
$\NE(\Lambda_{k,\Bbb R})=\{r\in\Lambda_{\Bbb R};\,
\langle r,\Am(\Lambda_{k,\Bbb R})\rangle>0\}$.
\par
Put $M:=T_{k}=\II_{1,1}(2)\oplus\I_{1,9-k}(2)$. Vectors in $M$ are denoted
by $(m,n,\lambda)$ ($m,n\in\Bbb Z$, $\lambda\in\I_{9-k}(2)$) whose norm is
$4mn-\langle\lambda,\lambda\rangle$. The period domain $\Omega_{S_{k}}$ is
isomorphic to the tube domain $\Lambda_{k,\Bbb R}+\sqrt{-1}C(\Lambda_{k})$
via the following map;
$$
\Lambda_{k,\Bbb R}+\sqrt{-1}C(\Lambda_{k})^{\pm}\ni v\to
\left(1/2,-\langle v,v\rangle/2,v\right)\in\Omega_{S_{k}}^{\pm}.
\tag 8.7
$$
We define $e_{0},e_{1},v_{0},v_{1},v_{2},v_{3}\in\Bbb C[M^{\lor}/M]$ by
$$
e_{0}:=e_{(0,0,0)},\quad
e_{1}:=e_{(0,0,\rho_{k})},\quad
v_{i}:=\sum_{2\langle\delta,\delta\rangle\equiv i\mod4}e_{\delta}.
\tag 8.8
$$
Put $q=e^{2\pi i\tau}$.
Let $\theta_{A_{1}+\delta/2}(\tau):=\sum_{k\in\Bbb Z}q^{(k+\delta/2)^{2}}$
$(\delta\in\{0,1\})$ be the theta function of $A_{1}$-lattice. Define
$f_{0}(\tau)$, $f_{1}(\tau)$, and $\{c_{k,0}(l)\}_{l\in\Bbb Z}$,
$\{c_{k,1}(l)\}_{l\in\Bbb Z+1/4}$ by
$$
\align
f_{0}(\tau):&=\frac{\eta(2\tau)^{8}\theta_{A_{1}}(\tau)^{k}}
{\eta(\tau)^8\eta(4\tau)^8}=\sum_{l\in\Bbb Z}c_{k,0}(l)\,q^{l}=
q^{-1}+8+2k+O(q),
\tag 8.9\\
f_{1}(\tau):&
=-16\frac{\eta(4\tau)^8\theta_{A_{1}+1/2}(\tau)^k}{\eta(2\tau)^{16}}
=\sum_{l\in1/4+\Bbb Z}2c_{k,1}(l)\,q^{l}.
\tag 8.10
\endalign
$$
Then, $f_{0}(\tau)$ is a modular form of weight $(k-8)/2$ for 
$\Gamma_{0}(4)$ with the same character as that of $\theta_{A_1}(\tau)$. 
Define $g_{i}(\tau)$ $(i\in\Bbb Z/4\Bbb Z)$ by
$$
g_{i}(\tau)=\sum_{l\equiv i\mod 4}c_{k,0}(l)\,q^{l/4},
\quad\sum_{i\in\Bbb Z/4\Bbb Z}g_{i}(\tau)=
\frac{\eta(\tau/2)^8\theta_{A_1}(\tau/4)^k}{\eta(\tau)^8\eta(\tau/4)^8}=
f_{0}(\tau/4),
\tag 8.11
$$
and a modular form of type $\rho_{M}$ by
$$
F_{k}(\tau):=f_{0}(\tau)\,e_0+f_{1}(\tau)\,e_1
+\sum_{i\in\Bbb Z/4\Bbb Z}g_{i}(\tau)\,v_{i}.
\tag 8.12
$$

\proclaim{Theorem 8.2}
If $k<8$, one has $\varDelta_{S_{k}}(z)^2=
C_{k}\,\Psi_{\roman{I}_{2,10-k}(2)}(z,F_k)$ and the following infinite
product expansion;
$$
\varDelta_{S_{k}}(z)=C_{k}\,e^{2\pi i\langle\rho_{k},z\rangle}
\prod_{\delta\in\{0,1\}}\prod_{r\in\Pi^{+}_{\delta}(\Lambda_{k})}
(1-e^{2\pi i\langle r,z\rangle})^{c_{k,\delta}(r\cdot r/2)}
$$
where $\Pi^{+}_{\delta}(\Lambda_{k}):=
(\delta\rho_{k}+\Lambda_{k})\cap\NE(\Lambda_{k,\Bbb R})$.
\endproclaim

Let
$M\Gamma_{0}(4)=\{(\binom{a\,b}{c\,d},\sqrt{c\tau+d})\in Mp_2(\Bbb Z);\,
c\equiv0\mod4\}$ be a subgroup of $Mp_{2}(\Bbb Z)$. Put $Z:=S^2$. 

\proclaim{Lemma 8.1}
There exists a character $\chi$ of $M\Gamma_{0}(4)$ such that
$\rho_{M}(g)\,e_{0}=\chi(g)^{8-k}e_{0}$.
\endproclaim

\demo{Proof}
From (8.1), it follows that
$\ker\rho_{M}\backslash M\Gamma_{0}(4)=\{Z^{a}T^{b};\,0\leq a,b\leq3\}$ 
and $\rho_{M}(g)\,e_{\gamma}=
i^{a(8-k)}e^{b\pi i\langle\gamma,\gamma\rangle}\,e_{\gamma}$ if 
$g\in\ker\rho_{M}\cdot Z^{a}T^{b}$. Putting $\chi(g)=i^{a}$, we get the
assertion.\qed
\enddemo

For $\gamma=(\binom{a\,b}{c\,d},\sqrt{c\tau+d})\in Mp_{2}(\Bbb Z)$, put
$j(\gamma,\tau):=\sqrt{c\tau+d}$. For $g\in Mp_2(\Bbb Z)$, we set 
$f|_{g}(\tau):=f(g\cdot\tau)j(g,\tau)^{-l}$ where $l$ is the weight of 
$f$. The following construction of modular forms of type $\rho_{M}$ from 
a scalar valued modular form of higher level, is due to Borcherds 
([Bo6,7]). (It is a special case of his construction.)

\proclaim{Proposition 8.1 ([Bo6,7])}
Let $\phi(\tau)$ be a modular form for $\Gamma_{0}(4)$ with character 
$\chi^{8-k}$. Then, 
$\tilde{\phi}(\tau):=\sum_{g\in\Gamma_{0}(4)\backslash\Gamma(1)}
f|_{g}(\tau)\rho_{M}(g^{-1})\,e_{0}$ is well defined and
becomes a modular form of type $\rho_{M}$.
\endproclaim

\proclaim{Definition 8.1}
If the level of $M$ is $4$ and $F=\tilde{\phi}$ for some modular form
for $\Gamma_{0}(4)$ in Theorem 8.1, we write $\Psi_{M}(z,\phi)$ in stead
of $\Psi_{M}(z,\tilde{\phi})$. $\Psi_{M}(z,\phi)$ is said to be Borcherds's
product attached to $M$ and $\phi$.
\endproclaim

\proclaim{Lemma 8.2}
$F_{k}(\tau)$ is a modular form of weight $(k-8)/2$ of type $\rho_{M}$.
\endproclaim

\demo{Proof}
Pick up the representatives
$\Gamma_{0}(4)\backslash\Gamma(1)=\{1,S,ST,ST^2,ST^3,V\}$. Here we put
$V:=S^{-1}T^{2}S$. From (8.1), it follows that
$$
\rho_{M}((ST^{l})^{-1})\,e_{0}=
i^{\frac{k-8}{2}}2^{-\frac{12-k}{2}}\sum_{j=0}^{3}i^{-lj}\,v_{j},
\quad\rho_{M}(V^{-1})\,e_{0}=e_{1}.
\tag 8.13
$$
Take $\phi(\tau)=f_{0}(\tau)$ in Proposition 8.1. As 
$\theta_{A_{1}}(\tau)$ is a modular form of weight $1/2$ with character
$\chi^{-1}$ for $\Gamma_{0}(4)$ (cf. [Bo5, Theorem 4.1]), 
$f_{0}(\tau)$ is a modular form of weight $(k-8)/2$ for $\Gamma_{0}(4)$ 
with character $\chi^{8-k}$. (Note that $\chi^8\equiv1$.) Since
$f_{0}|_{ST^{l}}(\tau)=2^{(8-k)/2}i^{-k/2}f_{0}((\tau+l)/4)$, comparing 
with (8.13), we get
$$
\sum_{l=0}^{3}f_{0}|_{ST^{l}}(\tau)\cdot\rho_{M}((ST^{l})^{-1})\,e_{0}=
\sum_{i=0}^{3}g_{i}(\tau)\,v_{i}.
\tag 8.14
$$
Similarly, since $f_{0}|_{V}(\tau)=f_{1}(\tau)$, comparing with (8.13),
$$
f_{0}|_{V}(\tau)\cdot\rho_{M}(V^{-1})\,e_{0}=f_{1}(\tau)\,e_{1}.
\tag 8.15
$$
Thus, $F_{k}(\tau)=\tilde{f_{0}}(\tau)$. As the weight of $f_{0}(\tau)$ 
is $(k-8)/2$, so is $F_{k}(\tau)$.\qed
\enddemo

\demo{Proof of Theorem $8.2$}
Since it follows from (8.9-12) that
$$
\aligned
F_{k}(\tau)
&=(q^{-1}+8+2k+O(q))\,e_{0}+(8+2k+O(q))\,v_{0}+O(q^{1/4})\,v_{1}+
O(q^{1/2})\,v_{2}\\
&\quad+(q^{-1/4}+O(q^{3/4}))\,v_{3}+O(1)\,e_{1},
\endaligned
\tag 8.16
$$
$\Psi_{T_{k}}(z,f_{0}):=\Psi_{T_k}(z,F_{k})$ is a modular form for some 
arithmetic subgroup of $O(T_{k})$ of weight $8+2k$. Note that $v_{0}$
contains $e_0$ with multiplicity one. By Theorem 8.1 (2), 
the zero of $\Psi_{T_k}(z,f_{0})$ consists of all the hyperplanes 
perpendicular to the root of $T_{k}$ with multiplicity $2$. Thus, 
$\varDelta_{S_{k}}(z)^{2}$ and $\Psi_{T_k}(z,f_{0})$ have the same zero 
and weight which prove the first assertion. By [Bo5, Theorem 10.4], 
the Weyl vector of $\Psi(z,f_{0})$ is $2\rho_{k}$. Let 
$f_{0}(\tau/4)=\sum_{l\in4^{-1}\Bbb Z}c(l)\,q^{l}$ be the Fourier 
expansion at the cusp. Using 
$c(\langle r,r\rangle/2)=c_{0,k}(\langle2r,2r\rangle/2)$ for any
$r\in\Lambda_{k}^{\lor}$, it follows from Theorem 8.1 (3) that
$$
\aligned
\Psi_{T_k}(z,f_{0})
&=e^{4\pi i\langle\rho_{k},z\rangle}
\prod_{r\in\Lambda_{k}\cap\roman{NE}(\Lambda_{k,\Bbb R})}
(1-e^{2\pi i\langle r,z\rangle})^{c_{0,k}(r^2/2)}\\
&\quad\times\prod_{r\in(\rho_{k}+\Lambda_{k})\cap 
\roman{NE}(\Lambda_{k,\Bbb R})}
(1-e^{2\pi i\langle r,z\rangle})^{2c_{1,k}(r^2/2)}\\
&\quad\times
\prod_{r\in\Lambda_{k}^{\lor}\cap\roman{NE}(\Lambda_{k,\Bbb R})}
(1-e^{2\pi i\langle r,z\rangle})^{c(r^2/2)}
(1+e^{2\pi i\langle r,z\rangle})^{c(r^2/2)}\\ 
&=e^{4\pi i\langle\rho_{k},z\rangle}\prod_{r\in\Pi_{0}^{+}(\Lambda_{k})}
(1-e^{2\pi i\langle r,z\rangle})^{2c_{0,k}(r^2/2)}
\prod_{r\in\Pi_{1}^{+}(\Lambda_{k})} 
(1-e^{2\pi i\langle r,z\rangle})^{2c_{1,k}(r^2/2)}\\ 
&=\left[e^{2\pi i\langle\rho_{k},z\rangle}\prod_{\delta\in\{0,1\}}
\prod_{r\in\Pi_{\delta}^{+}(\Lambda_{k})}
(1-e^{2\pi i\langle r,z\rangle})^{c_{\delta,k}(r^2/2)}\right]^{2}.
\qed
\endaligned
\tag 8.16
$$
\enddemo

\subsubhead
\bf{8.3 Nikulin's $\pmb{K3}$ Surfaces}\rm
\endsubsubhead
For $\Lambda=\II_{1,1}\oplus E_{8}(-2)$, a $\Lambda$-2-elementary $K3$ 
surface is one of the exceptional type ($(r,l,\delta)=(10,8,0)$) in 
Theorem 2.5 discovered by Nikulin ([Ni4]). 
It is an elliptic $K3$ surface whose fixed locus consists of two smooth 
fibers. Let $\Psi_{\Lambda^{\perp}}:=\Phi|_{\Omega_{\Lambda}}$ be the same 
automorphic form as in Corollary 6.1.
\par
Put $\theta_{\Lambda_{16}+\delta}(\tau)=
\sum_{\lambda+\delta\in\Lambda_{16}}q^{\lambda^{2}/2}$ for the theta series 
of $\Lambda_{16}$ where $\delta\in\Lambda_{16}^{\lor}/\Lambda_{16}$. 
By [Bo2], there exists an involution $I$ on $\Lambda_{24}$ (the Leech 
lattice) whose fixed lattice is $E_{8}(2)$ and anti-fixed lattice is
$\Lambda_{16}$. As $\Lambda_{24}$ is self-dual, there exists a canonical
identification 
$\Lambda_{16}^{\lor}/\Lambda_{16}\cong E_{8}(2)^{\lor}/E_{8}(2)
=\Lambda^{\lor}/\Lambda$ via which we identify the generators of the 
standard basis of the group rings 
$\Bbb C[\Lambda_{16}^{\lor}/\Lambda_{16}]$ and that of
$\Bbb C[\Lambda^{\lor}/\Lambda]$ (cf. [Bo5, $\S4$]). Under this 
identification, $\rho_{\Lambda_{16}}=\rho_{\Lambda}$ because
$\delta(\Lambda_{16})=\delta(\Lambda)=0$. Let 
$\Theta_{\Lambda_{16}}(\tau)=
\sum_{\delta\in\Lambda_{16}^{\lor}/\Lambda_{16}}
\theta_{\Lambda_{16}+\delta}(\tau)\,e_{\delta}$ 
be the theta series of $\Lambda_{16}$. Let $\{c_{\delta}(k)\}$ be the
Fourier coefficients of the following modular form;
$$
\Theta_{\Lambda_{16}}(\tau)/\Delta(\tau)=
\sum_{\delta\in\Lambda_{16}^{\lor}/\Lambda_{16}}e_{\delta}\,
\theta_{\Lambda_{16}+\delta}(\tau)/\Delta(\tau)=
\sum_{\delta\in\Lambda_{16}^{\lor}/\Lambda_{16}}
\sum_{k\in\Bbb Z}c_{\delta}(k)\,q^{k/2}\,e_{\delta}.
\tag 8.17
$$ 

\proclaim{Theorem 8.3}
$\Psi_{\Lambda^{\perp}}(z)=
\Psi_{\Lambda\oplus\roman{II}_{1,1}}(z,\Theta_{\Lambda_{16}}/\Delta)$.
\endproclaim

\demo{Proof}
From [C-S, Chap.4] (note that $q=e^{\pi i\tau}$ in [C-S] although
$q=e^{2\pi i\tau}$ here), it follows that 
$\theta_{\Lambda_{16}}(\tau)=1+O(q^{2})$, 
$1/\Delta(\tau)=q^{-1}+24+O(q)$, and
$$
\aligned
\theta_{\Lambda_{16}^{\lor}}(\tau)
&=\sqrt{\det\Lambda_{16}}(i/\tau)^{8}\theta_{\Lambda_{16}}(-1/\tau)\\
&=2^{4}\cdot2^{-9}[\{\theta_{3}(\tau)^{2}-\theta_{2}(\tau)^{2}\}^{8}+
\{\theta_{3}(\tau)^{2}+\theta_{2}(\tau)^{2}\}^{8}+
\theta_{3}(\tau)^{8}\theta_{2}(\tau)^{8}\\
&\quad+30\{\theta_{3}(\tau)^{2}-\theta_{2}(\tau)^{2}\}^{4}\cdot
\{\theta_{3}(\tau)^{2}+\theta_{2}(\tau)^{2}\}^{4}]\\
&=2^{-5}[(1-2q^{1/4})^{16}+(1+2q^{1/4})^{16}+
(4q^{1/4})^{4}(1+4q^{1/2})^{4}\\
&\quad+30(1-2q^{1/4})^{8}(1+2q^{1/4})^{8}]+O(q)=1+O(q).
\endaligned
\tag 8.18
$$
In particular, $\delta\not=0$ implies that 
$\theta_{\Lambda_{16}+\delta}(\tau)=O(q)$ and thus
$$
\Theta_{\Lambda_{16}}(\tau)/\Delta(\tau)=(q^{-1}+24+O(q))\,e_{0}+
\sum_{\delta\not=0}O(1)\,e_{\delta}.
\tag 8.19
$$
From Thereom 8.1, it follows that
$\Psi_{\Lambda\oplus\roman{II}_{1,1}}(z,F)$ is a holomorphic modular 
form of weight $12$ whose zero divisor coincides with the discriminant.
Comparing weight and zeros, we get the assertion.\qed
\enddemo

\beginsection
$\S9$. GKM Superalgebras Arizing from $\pmb{\varDelta_{M}}$

\par
\subsubhead
\bf{9.1 Generalized Kac-Moody Superalgebras and $\pmb{K3}$ Surfaces}\rm
\endsubsubhead
Following [G-N2], let us recall generalized Kac-Moody superalgebras
(GKM superalgebra for short) associated to an algebraic $K3$ surface.
(For the general theory of GKM (super)algebras, see [Bo1,2] and [G-N1-3].)
\par
Let $X$ be an algebraic $K3$ surface and $S:=\Pic_{X}$ its Picard lattice.
Let $\Exc(S)$ be the set of all $-2$-curves in $X$. 
Let $C(S)=\{v\in S_{\Bbb R};\,\langle v,v\rangle>0\}$ be the light cone 
of $S$ and $C^{+}(S)$ be the connected component containing the ample class. 
Let $\Am(S_{\Bbb R})$ be the ample cone:
$\Am(S_{\Bbb R}):=\{l\in C^{+}(S);\,\langle l,\delta\rangle>0,\,
\forall\delta\in \Exc(S)\}$. Put 
$\overline{\Am}(S):=\Am(S_{\Bbb R})\cap S^{\lor}$ where closure is 
considered in $C^{+}(S)$. Let $W(S)$ be the Weyl group of $S$.
Then, $\Delta(S)=W(S)(\Exc(S))$ and $\overline{\Am}(S_{\Bbb R})$ is the 
fundamental domain for the action of $W(S)$ on $C^{+}(S)$.
\par
To define a Lie superalgebra associated to $X$, we need the set
of simple roots. Let us put ${}_{s}\Delta^{re}:=\Exc(S)$ for the set of all 
simple real roots. Let ${}_{s}\Delta^{im}_{\bar0}$ (resp.
${}_{s}\Delta^{im}_{\bar1}$) be a sequences of elements in 
$\overline{\Am}(S)$ such that any $a\in\overline{\Am}(S)$ can appear in 
${}_{s}\Delta^{im}_{\bar*}$ finitely many times;
$$
{}_{s}\Delta^{im}_{\bar*}=\{m(a)_{\bar*}a;\,a\in\overline{\Am}(S),\,
\langle a,a\rangle>0\}\cup\{\tau(a)_{\bar*}a;\,a\in\overline{\Am}(S),\,
\langle a,a\rangle=0\}.
\tag 9.1
$$
Here $m(a)_{\bar*},\tau(a)_{\bar*}\in\Bbb Z$, and $m(a)_{\bar*}a$ 
(resp. $\tau(a)_{\bar*}a$) implies that $a$ appears $m(a)_{\bar{*}}$ 
(resp. $\tau(a)_{\bar{*}}$) times in ${}_{s}\Delta^{im}_{\bar*}$. An 
element of ${}_{s}\Delta^{im}_{\bar0}$ (resp. ${}_{s}\Delta^{im}_{\bar1}$) 
is called a simple even (resp. odd) imaginary root. Put 
${}_{s}\Delta^{im}:={}_{s}\Delta^{im}_{\bar0}\cup{}_{s}\Delta^{im}_{\bar1}$ 
for the set of all simple imaginary roots, and 
${}_{s}\Delta:={}_{s}\Delta^{re}\cup{}_{s}\Delta^{im}$ for the set of
all simple roots. Let us write ${}_{s}\Delta=\{h_{i}\}_{i\in I}$ 
($h_{i}\in S^{\lor}$). Let $A=(a_{ij})_{i,j\in I}$ be the Gramm matrix.
Namely, $a_{ij}=\langle h_{i},h_{j}\rangle$. As $A$ satisfies the
conditions of generalized Cartan matrix ([Bo1]):
\newline
(1) $i\not=j$ implies $a_{ij}\geq0$, \quad
(2) $h_{i}\in{}_{s}\Delta^{re}$ implies $a_{ii}=-2$ and $a_{ij}\in\Bbb Z$,
\newline
we get a GKM superalgebra $\frak g(S,{}_{s}\Delta^{im}):=\frak g'(A)$.
(See [Bo1,2], [G-N1,2] for details.) Let 
$\Pi^{+}:=\sum_{\alpha\in{}_{s}\Delta}\Bbb Z_{+}\alpha\backslash\{0\}
\subset\NE(S):=\NE(S_{\Bbb R})\cap S^{\lor}$ be the set of positive roots
where $\NE(S_{\Bbb R})$ is the dual cone (Mori cone) to the ample cone 
$\Am(S_{\Bbb R})$. Let $\frak g_{\alpha}$ be the root space attached to 
$\alpha\in\Pi^{+}\cup\{0\}\cup-\Pi^{+}$.
Then, $\frak g(S,{}_{s}\Delta^{im})$ admits root space decomposition:
$\frak g(S,{}_{s}\Delta^{im})=(\oplus_{\alpha\in\Pi^{+}}\frak g_{\alpha})
\oplus\frak g_{0}\oplus(\oplus_{\alpha\in\Pi^{+}}\frak g_{-\alpha})$,
$\frak g_{0}=S_{\Bbb R}$.
According to the decomposition into even and odd part, we get
$\frak g_{\alpha}=\frak g_{\alpha,\bar0}\oplus\frak g_{\alpha,\bar1}$. 
Multiplicity of $\alpha\in\Pi^{+}$ is defined by $\mult(\alpha):=
\dim\frak g_{\alpha,\bar0}-\dim\frak g_{\alpha,\bar1}\in\Bbb Z$.
\par
In view of Borcherds's works, it is the denominator formula that connects 
GKM superalgebras and automorphic forms. We recall it when $S$ admits a 
Weyl vector.

\proclaim{Definition 9.1}
$\rho\in S_{\Bbb Q}$ is a Weyl vector if $\langle\rho,\delta\rangle=1$ 
for all $\delta\in\Exc(S)$.
\endproclaim

For $a\in{}_{s}\Delta^{im}$, define $m(a)\in\Bbb Z$ by
(1) $m(a)=m(a)_{\bar0}-m(a)_{\bar1}$ if $\langle a,a\rangle>0$, 
and (2) the following formal series if $\langle a,a\rangle=0$ and $a$ is 
primitive,
$$
\prod_{n=1}^{\infty}(1-q^{n})^{\tau(na)}=1-\sum_{k=1}^{\infty}m(ka)\,q^{k}
\quad(\tau(na):=\tau(na)_{\bar0}-\tau(na)_{\bar1}).
\tag 9.2
$$

\proclaim{Theorem 9.1 ([Bo1,2] (cf. [G-N1,2]))}
Suppose that $S$ has a Weyl vector $\rho$. 
\flushpar{$(1)$}
For a GKM superalgebra $\frak g(S,{}_{s}\Delta^{im})$, one has the 
following identity;
$$
\aligned
\Phi_{\frak g(S,{}_{s}\Delta^{im})}(z):
&=\sum_{w\in W(S)}\det(w)
\{e^{2\pi i\langle w(\rho),z\rangle}-
\sum_{r\in{}_{s}\Delta^{im}}m(r)\,
e^{2\pi i\langle w(\rho+r),z\rangle}\}\\
&=e^{2\pi i\langle\rho,z\rangle}\prod_{\alpha\in\Pi^{+}}
(1-e^{2\pi i\langle\alpha,z\rangle})^{\roman{mult}(\alpha)}.
\endaligned
$$
The formal series $\Phi_{\frak g(S,{}_{s}\Delta^{im})}(z)$ is said to be 
the denominator function.
\flushpar{$(2)$}
Let $\Psi(y)$ be a formal series with the following integral Fourier 
expansion:
$$
\Psi(z)=\sum_{w\in W(S)}\det(w)\{e^{2\pi i\langle w(\rho),z\rangle}-
\sum_{r\in\overline{\roman{Am}}(S)}
m(r)\,e^{2\pi i\langle w(\rho+r),z\rangle}\}.
$$ 
Define $\tau(nr)\in\Bbb Z$ by $(8.2)$ for primitive norm zero
$r\in\overline{\Am}(S)$  and $n\in\Bbb N$. Let 
$\frak g(S,{}_{s}\Delta^{im})$ be the GKM superalgebra whose simple 
imaginary roots are
$$
\aligned
\,&
{}_{s}\Delta^{im}_{\bar0}=\{m(r)r;\,m(r)>0,\,\langle r,r\rangle>0\}\cup
\{\tau(r)r;\,\tau(r)>0,\,\langle r,r\rangle=0\},\\
&{}_{s}\Delta^{im}_{\bar1}=\{-m(r)r;\,m(r)<0,\,\langle r,r\rangle>0\}\cup
\{-\tau(r)r;\,\tau(r)<0,\,\langle r,r\rangle=0\}.
\endaligned
$$
Then, $\Phi_{\frak g(S,{}_{s}\Delta^{im})}(z)=\Psi(z)$.
\endproclaim

Among all the GKM superalgebras associated to $K3$ surfaces, one of the 
fake monster Lie algebras constructed by Borcherds is the most beautiful 
and interesting.
\par 
Let $\Lambda,S$ and $T$ be the 2-elementary lattices defined by
$\Lambda:=\II_{1,1}\oplus E_{8}(-2)$, $T=\II_{1,1}(2)\oplus\Lambda$,
$S=T^{\perp}=\II_{1,1}(2)\oplus E_{8}(-2)$. A 2-elementary $K3$ surface 
of type $S$ is the universal cover of an Enriques surface and $T$ is the
transcendental lattice. The period domain $\Omega_{S}$ is realized as 
the tube domain $\Lambda_{\Bbb R}+\sqrt{-1}C(\Lambda)$ as before 
(cf. (8.7)). Set $\rho:=(0,1,0)$, $\rho':=(1,0,0)$. Then, $\rho$ is a 
Weyl vector of $\Lambda$. 

\proclaim{Theorem 9.2 ([Bo2,4])}
There exists a GKM superalgebra $\frak g(\Lambda,{}_{s}\Delta^{im})$
(one of the fake monster Lie algebras) whose denominator 
function, $\Phi$, is the automorphic form over $\Omega_{S}$ of weight $4$ 
with zero divisor $\Cal D_{S}$. For $\Im y\gg0$, one has the following;
$$
\aligned
\Phi(z)
&=\sum_{w\in W(\Lambda)}\det(w)\,e^{2\pi i\langle\rho,w(z)\rangle}
\prod_{n>0}(1-e^{2\pi in\langle\rho,w(z)\rangle})^{(-1)^{n}8}\\
&=e^{2\pi i\langle\rho,z\rangle}\prod_{r\in\Pi^{+}}
(1-e^{2\pi i\langle r,z\rangle})^{(-1)^{r\cdot(\rho-\rho')}c(r^2/2)}
\endaligned
$$
where $\Pi^{+}=\Bbb N\rho\cup\{r\in\Lambda;\,\rho\cdot r>0\}$ and
$\sum_{n\geq-1}c(n)\,q^{n}=
\eta(\tau)^{-8}\eta(2\tau)^{8}\eta(4\tau)^{-8}$.
\endproclaim

\proclaim{Theorem 9.3}
There exists a constant $C_{S}\not=0$ such that 
$\varDelta_{S}=C_{S}\,\Phi$. $\Phi(z)^2$ is Borcherds's product attached 
to the modular form $F_0$ of $(8.12)$.
\endproclaim

\demo{Proof}
Comparing the weight and zero of $\varDelta_{S}$ and $\Phi$, we get the 
first assertion. Regarding $\I_{1,9}(2)$ as a sublattice of $\Lambda$,
we can prove the second assertion in the same manner as Theorem 8.1.
\qed
\enddemo

\remark{Remark}
We remark that $\I_{2,10}(2)$ (not the transcendental lattices 
$\II_{1,1}\oplus\II_{1,9}(-2)$) is used in (3) to construct Borcherds's 
product. In [A], [Kn], Allcock and Kond\B o uses the similar relation 
between $\Lambda$ and $\I_{1,9}$ to study the moduli space of Enriques 
surfaces.
\endremark

\subsubhead
\bf{9.2 GKM Superalgebras Arising from $\pmb{\varDelta_{M}}$}
\rm
\endsubsubhead
We keep the notations in $\S8.1$, and study two classes of 2-elementary 
$K3$ surfaces as before.

\proclaim{Theorem 9.4}
There exists a GKM superalgebra attached to 
$\Lambda_{k}$ whose denominator function is (up to a constant) 
$\varDelta_{S_{k}}(z)$.
\endproclaim

\demo{Proof}
In view of Theorem 9.1 (2), it is enough to show that 
$\varDelta_{S_{k}}(z)$ admits the following integral Fourier expansion 
at the cusp;
$$
\varDelta_{S_{k}}(z)=C_{S_{k}}\,\sum_{w\in W(\Lambda_{k})}\det(w)
\{e^{2\pi i\langle w(\rho_{k}),z\rangle}-
\sum_{r\in\overline{\roman{Am}}(\Lambda_{k})}
m(r)\,e^{2\pi i\langle w(\rho_{k}+r),z\rangle}\}.
\tag 9.3
$$
By Theorem 8.2, $\varDelta_{S_{k}}(z)$ has the following Fourier 
expansion at the cusp;
$$
\varDelta_{S_{k}}(z)=C_{S_{k}}\,
\sum_{r\in\Lambda_{k}^{\lor}\cap C(\Lambda_{k})^{+}}n(r)\,
e^{2\pi i\langle r,z\rangle}
\tag 9.4
$$
where $n(r)\in\Bbb Z$. By (8.4), (8.5), Theorems 7.3 and 9.2, one has
$$
\varDelta_{S_{k}}(w(z))=\lim_{v\to z}\left(\frac{i}{2\pi}\right)^{k}
\frac{\Phi(w(v))}{\prod_{i=0}^{k-1}\langle w(v),\delta_{i}\rangle}
=\det(w)\,\varDelta_{S_{k}}(z)\quad(w\in W(\Lambda_k))
\tag 9.5
$$
because $w(\delta_{i})=\delta_{i}$ for $i\leq k-1$. Namely, 
$n(w(r))=\det(w)\,n(r)$ for any $r\in\Lambda_{k}^{\lor}$ and 
$w\in W(\Lambda_{k})$. Together with the same argument as 
[G-N1, Theorem 2.3 (a)], one has 
$r-\rho_{k}\in\overline{\Am}(\Lambda_{k})$ if $n(r)\not=0$. 
Putting $m(r):=-n(r+\rho_{k})$, (9.4) becomes
$$
\varDelta_{S_{k}}(z)=C_{S_{k}}\,\sum_{w\in W(\Lambda_{k})}\det(w)
\sum_{r\in\overline{\roman{Am}}(\Lambda_{k})}
-m(r)\,e^{2\pi i\langle w(\rho_{k}+r),z\rangle}.
\tag 9.6
$$
By Theorem 8.2, $-m(0)=1$ which, together with (9.6), yields the assertion.
\qed
\enddemo

To study the case of $2$-elementary $K3$ surface of type $\Lambda$, 
Borcherds's $\Phi$-function of rank $26$ is crucial. Let $\Lambda_{24}$ 
be the Leech lattice, and put $L:=\II_{1,1}\oplus\Lambda_{24}(-1)$. 
Let $\rho=(0,1,0)$ be a Weyl vector of $L$.

\proclaim{Theorem 9.5 ([Bo2,3])}
The denominator function $\Phi$ of the fake monster Lie algebra is the 
automorphic form over $\Omega_{\roman{II}_{1,1}\oplus L}$ of weight $12$ 
with only simple zero along the discriminant. The denominator formula 
becomes
$$
\Phi(z)
=\sum_{w\in W(L)}\sum_{n>0}\det(w)\tau(n)\,
e^{2\pi i\langle w(\rho),z\rangle}
=e^{2\pi i\langle\rho,z\rangle}\prod_{r\in\Pi^{+}}
(1-e^{2\pi i\langle r,z\rangle})^{p_{24}(1-r^2/2)}
$$
where $\tau(n)$ is the Ramanujan $\tau$-function, $p_{24}(n)$ the number 
of partitions of $n$ into $24$ colors, and 
$\Pi^{+}=\Bbb N\rho\cup\{r\in\Lambda;\,\rho\cdot r>0\}$.
\endproclaim

\proclaim{Theroem 9.6}
$\Psi_{\Lambda^{\perp}}$ is the denominator function of a GKM 
superalgebra. 
\endproclaim

\demo{Proof}
By Theorem 8.3, $\Psi_{\Lambda^{\perp}}(z)$ has the similar Fourier 
expansion as (9.4). Since $\Psi_{\Lambda^{\perp}}=
\Phi|_{\Omega_{\roman{II}_{1,1}\oplus\Lambda}}$, it follows from 
Theorem 9.5 that 
$\Psi_{\Lambda^{\perp}}(w(z))=\det(w)\,\Psi_{\Lambda^{\perp}}(z)$. 
Thus, we get the assertion in the same manner as Theorem 9.4.\qed
\enddemo

\beginsection
$\S10$. An Explicit Formula for $\pmb{\varDelta_{A_{1}}}$ and
$\pmb{\varDelta_{\roman{II}_{1,1}(2)}}$

\par
We determine an explicit formula for $\varDelta_{A_{1}}$ and
$\varDelta_{\roman{II}_{1,1}(2)}$ in this section.
\par
\subsubhead
\bf{10.1 An Explicit Formula for $\pmb{\varDelta_{A_{1}}}$}\rm
\endsubsubhead
Let $(X,\iota)$ be a $2$-elementary $K3$ surface of type $A_{1}$.
Then, $X/\iota=\Bbb P^2$ and the fixed curve $X^{\iota}$ is a smooth plane 
sextic curve where the quotient map $X\to\Bbb P^2$ is the morphism 
associated to the complete linear system of ample line bundle of degree $2$.
As $X^{\iota}\subset\Bbb P^2$ has an ambiguity of $PGL(3,\Bbb C)$ (arizing
from a choice of $3$ sections), there exists a morphism 
$i:\Cal M^{0}_{A_{1}}\ni[(X,\iota)]\to
[X^{\iota}]\in H^{sm}_{6}/PGL(3,\Bbb C)$ where 
$H^{sm}_{6}=\Bbb P(\Sym^{6}\Bbb C^{3})^{\lor}\backslash D_6$ is the set 
of all $smooth$ sextic curves in $\Bbb P^2$, and $D_6$ is the discriminant
locus of universal plane sextic curves 
$\pi:\Cal C_{6}\to H_{6}=\Bbb P(\Sym^6\Bbb C^3)^{\lor}$. Note that
$\xi=(\xi_{I})\in\Bbb P(\Sym^6\Bbb C^3)^{\lor}$ corresponds to the curve
$C_{\xi}=\{x\in\Bbb P^2;\,\sum_{|I|=6}\xi_{I}x^{I}=0\}$. Conversely, if
$\xi,\xi'\in H^{sm}_6$ are in the same orbit of $PGL(3,\Bbb C)$, $C_{\xi}$
and $C_{\xi'}$ are projectively equivalent and thus the double covers of
$\Bbb P^2$ blanching along $C_{\xi}$ and $C_{\xi'}$ are isomorphic. 
In this way, one verify 
$i:\Cal M^{0}_{A_{1}}\cong H^{sm}_{6}/PGL(3,\Bbb C)$ ([Sha]).
\par
We denote by $\lambda(\Cal C_{6}/H^{sm}_{6})$ the determinant of cohomology
in the sense of $\S3.1$. Let $\Jac:H^{sm}_{6}\ni\xi\to[\Jac(C_{\xi})]
\in\Cal A_{10}$ be the Torelli map. As $j_{A_{1}}=\Jac\circ i$ and 
$(\Jac)^{*}\Cal F_{10}=\lambda(\Cal C_{6}/H^{sm}_{6})^{GL(3,\Bbb C)}$, 
we find $\lambda_{A_{1}}=
i^{*}\lambda(\Cal C_{6}/H^{sm}_{6})^{GL(3,\Bbb C)}$ where
$\lambda(\Cal C_{6}/H^{sm}_{6})^{GL(3,\Bbb C)}$ is the sheaf of 
$GL(3,\Bbb C)$-invariant sections of $\lambda(\Cal C_{6}/H^{sm}_{6})$. 
Moreover, this identification is an isometry if $\lambda_{A_{1}}$ is 
equipped with the Petersson metric and if
$i^{*}\lambda(\Cal C_{6}/H^{sm}_{6})^{GL(3,\Bbb C)}$ with the 
$L^2$-metric. Our first task is to construct a section of 
$\lambda_{A_{1}}^{\otimes15}$ arizing from the discriminant of 
plane sextics. In the sequel, we put $F(x;\xi)=\sum_{|I|=6}\xi_{I}x^{I}$. 
By the Poincar\'e residue sequence, we get the following.

\proclaim{Lemma 10.1}
If we denote by $x,y,z$ the homogeneous coordinates of $\Bbb P^2$, then
$$
H^{0}(C_{\xi},\Omega^{1}_{C_{\xi}})=
\bigoplus_{i+j+k=3}\Bbb C\,\Res_{C_{\xi}}
\frac{x^{i}y^{j}z^{k}(x\,dy\wedge dz-y\,dx\wedge dz+z\,dx\wedge dy)}
{F(x,y,z;\xi)}.
$$
\endproclaim

We define a local section of $\lambda(\Cal C_{6}/H^{sm}_{6})$ by
$$
\omega(\xi):=\bigwedge_{i+j+k=3}\Res_{C_{\xi}}
\frac{x^{i}y^{j}z^{k}(x\,dy\wedge dz-y\,dx\wedge dz+z\,dx\wedge dy)}
{F(x,y,z;\xi)}\in\det H^{0}(C_{\xi},\Omega^{1}).
\tag 10.1
$$
To be precise, if $p:(\Sym^6\Bbb C^{3})^{\lor}\backslash\{0\}\to
\Bbb P(\Sym^6\Bbb C^3)^{\lor}$ is the natural projection, then
$\omega$ is a section of $p^{*}\lambda(\Cal C_{6}/H^{sm}_{6})$ over
$(\Sym^6\Bbb C^3)^{\lor}$.

\proclaim{Lemma 10.2}
For any $g\in GL(3,\Bbb C)$, $g^{*}\omega=\det(g)^{20}\omega$ where
$GL(3,\Bbb C)$ acts on $(\Sym^6\Bbb C^3)^{\lor}$ via the induced 
representation.
\endproclaim

\demo{Proof}
From the definition of representation of $GL(3,\Bbb C)$ on 
$(\Sym^6\Bbb C^3)^{\lor}$, it follows that 
$F(g\cdot(x,y,z);g\cdot\xi)=F(x,y,z;\xi)$. Since $GL(3,\Bbb C)$ acts on
$\det\Sym^3\Bbb C^3$ by $\det(\cdot)^{10}$, we get
$g\cdot\wedge_{i+j+k=3}x^{i}y^{j}z^{k}=
\det(g)^{10}\wedge_{i+j+k=3}x^{i}y^{j}z^{k}$. By computation, we get 
$g^{*}(x\,dy\wedge dz-y\,dx\wedge dz+z\,dx\wedge dz)=
\det(g)(x\,dy\wedge dz-y\,dx\wedge dz+z\,dx\wedge dz)$.
Together with all of these, we get the assertion: 
$g^{*}(\omega(g\cdot\xi))=\det(g)^{20}\omega(\xi)$.
\qed
\enddemo

\proclaim{Lemma 10.3}
There exists a homogeneous polynomial $D_{6}(\xi)$ of degree $75$ 
in the $\xi$-variable such that the discriminant locus of the plane
sextics is $\div(D_{6})$. Moreover, for any $g\in GL(3,\Bbb C)$,
$D_{6}(g\cdot\xi)=\det(g)^{-150}D_{6}(\xi)$.
\endproclaim

\demo{Proof}
Let $v_6:\Bbb P(\Bbb C^3)\hookrightarrow\Bbb P(\Sym^6\Bbb C^3)$ be the
Veronese embedding. Put $X:=v_6(\Bbb P(\Bbb C^3))$ and 
$X^{\lor}\subset\Bbb P(\Sym^6\Bbb C^3)^{\lor}$ for the projective dual
variety of $(X,\Cal O_{\Bbb P(\roman{Sym}^6\Bbb C^3)}(1))$. By Katz's 
formula ([Kt, Cor. 5.6]), $X^{\lor}$ is a projective hypersurface of 
degree
$$
\deg X^{\lor}
=\int_{\Bbb P(\Bbb C^3)}\frac{c(\Bbb P(\Bbb C^3))}
{(1+v^{*}_{6}c_{1}(\Bbb P(\Sym^6\Bbb C^3)))^2}
=\int_{\Bbb P(\Bbb C^3)}\frac{1+3H+3H^2}{(1+6H)^2}=75
\tag 10.2
$$
where $H=c_{1}(\Bbb P(\Bbb C^3))$ is the hyperplane section of 
$\Bbb P(\Bbb C^3)$. Take $D_{6}(\xi)$ as a defining equation of 
$X^{\lor}$. As $X^{\lor}$ is invariant under the action of $GL(3,\Bbb C)$,
there exists $l\in\Bbb Z$ such that
$D_{6}(g\cdot\xi)=\det(g)^{l}\,D_{6}(\xi)$. Putting 
$g=\lambda\,I_{3}$ $(\lambda\in\Bbb C^{\times})$, we find $l=-150$.\qed
\enddemo

\proclaim{Proposition 10.1}
$\Delta_{6}^{2}(\xi):=D_{6}(\xi)^{2}\cdot\omega(\xi)^{\otimes15}$ is a 
$GL(3,\Bbb C)$-invariant section of $p^{*}\lambda(\Cal C_{6}/H^{sm}_{6})$.
In particular, we may regard $\Delta_{6}^{2}\in 
H^{0}(\Cal M_{A_{1}}^{0},\lambda_{A_{1}}^{\otimes15})$.
\endproclaim

\demo{Proof}
The first assertion follows from Lemmas 10.2 and 10.3, and the second
from $H^{0}(H^{sm}_{6},
\lambda(\Cal C_{6}/H^{sm}_{6})^{\otimes 15})^{GL(3,\Bbb C)}=
H^{0}(\Cal M_{A_{1}}^{0},\lambda_{A_{1}}^{\otimes15})$.\qed
\enddemo

Let $\delta_{1},\delta_{2}\in\Delta(\II_{2,18}\oplus A_{1}(-1))$
such that $\delta_{1}/2\in(\II_{2,18}\oplus A_{1}(-1))^{\lor}$ and
$\delta_{2}/2\not\in(\II_{2,18}\oplus A_{1}(-1))^{\lor}$.
Since $\Gamma_{A_{1}}=O(\II_{2,18}\oplus A_{1}(-1))$,
it follows from a theorem of Nikulin ([Ni1, Proposition 1.15.1]) that
$\Delta(\II_{2,18}\oplus A_{1}(-1))/\Gamma_{A_{1}}=
\{\delta_{1},\delta_{2}\}$. In particular, 
$\Cal D_{A_{1}}/\Gamma_{A_{1}}$ consists of $2$
irreducible components: 
$\Cal D_{A_{1}}/\Gamma_{A_{1}}=
H_{\delta_{1}}+H_{\delta_{2}}$.

\proclaim{Lemma 10.4}
$\Delta_{6}^{2}$ vanishes of order $2$ along $H_{\delta_{2}}$.
\endproclaim

\demo{Proof}
Since $\langle A_{1}\oplus\delta_{2}\rangle=
A_{1}\oplus A_{1}(-1)$, $2$-elementary $K3$ surfaces over
$H_{\delta_{2}}$ are those of type $A_{1}\oplus A_{1}(-1)$.
Take a generic $2$-elementary K3 surface $(X,\iota)$ of this type. Then, 
$\Pic(X/\iota)=\I_{1,1}$ and $X/\iota$ is 
a blow-up of $\Bbb P^2$ at one point. This implies that the complete linear 
system associated to the nef and big line bundle of degree $2$ maps $X$ 
to $\Bbb P^2$ and $X^{\iota}$ to a nodal sextic curve with one singular 
point. Conversely, given such a nodal sextic curve, the minimal resolution 
of the double cover of $\Bbb P^2$ blanching along it is a $2$-elementary 
$K3$ surface of type $A_{1}\oplus A_{1}(-1)$. This extends the
isomorphism $i:\Cal M_{A_{1}}^{0}\cup H_{\delta_{2}}^{0}\to 
H^{nod}_{6}/PGL(3,\Bbb C)$ and $\lambda_{A_{1}}=
i^{*}\lambda(\Cal C_{6}/H^{nod}_{6})^{GL(3,\Bbb C)}$ where $H^{nod}_{6}$ 
is the set of all plane sextic curves with {\it at most one node}\rm. 
By Proposition 10.1, we find $\Delta_{6}^{2}\in H^{0}(H^{nod}_{6},
\lambda(\Cal C_{6}/H^{nod}_{6})^{\otimes15})^{GL(3,\Bbb C)}=
H^{0}(\Cal M_{A_{1}}^{0}\cup H^{0}_{\delta_{2}},
\lambda_{A_{1}}^{\otimes15})$. To compute the vanishing order of 
$\Delta_{6}^{2}$, take a generic point $p\in H_{\delta_{2}}$ and 
its small neighborhood $U$. As $\omega$ is regarded to be a nowhere 
vanishing section of $\lambda_{A_{1}}$ on $U$, the assertion follows 
from Proposition 10.1 because $H_{\delta_{2}}=\div(D_{6})$.\qed
\enddemo

\proclaim{Lemma 10.5}
$\Delta_{6}^{2}$ has an algebraic singularity along $H_{\delta_{1}}$.
\endproclaim

\demo{Proof}
By Shah's result ([Sha]), $H_{\delta_{1}}$ contracts to a 
point in $H^{ss}_{6}/PGL(3,\Bbb C)$ corresponding to a smooth conic of 
multiplicity $3$ where $H^{ss}_{6}$ is the set of all semi-stable sextic
curves. If $\{C_{\xi_{t}}\}_{|t|<1}$ is a degenerating family of sextic 
curves such that $C_{\xi_{0}}$ is a smooth triple conic, then 
$\|\Delta_{6}^{2}(\xi_{t})\|_{L^{2}}\sim|t|^{a}(\log|t|)^{b}$ 
$(t\to0)$ with some $a,b\in\Bbb R$. If $(U,s)$ is a holomorphic
disc which meets transversally to $H_{\delta_{1}}$ at $s=0$, then
$$
\|i^{*}\Delta_{6}^{2}(s)\|\sim|s|^{l}(\log|s|)^{m}\quad(s\to0)
\tag 10.3
$$
with some $l,m\in\Bbb R$ because $s$ and $t$ are algebraically related. 
Since $\langle A_{1},\delta_{1}\rangle=\II_{1,1}$ and 
$g(\II_{1,1})=10=g(A_{1})$, there exists a section, $\sigma$, 
of $\lambda_{A_{1}}^{\otimes15}$ over $U$ such that
$\|\sigma(s)\|\sim1$ $(s\to0)$ because $j_{A_{1}}(U)$ is away
from the boundary component of the Siegel modular variety 
$\Cal A_{10}$. Since $i^{*}\Delta_{6}^{2}(s)/\sigma(s)$ is a holomorphic
function on $U\backslash\{0\}$, it must be meromorphic on $U$ because of
(10.3). In particular, $m=0$ and $\Delta_{6}^{2}$ has zero (or pole) of 
order $l\in\Bbb Z$ along $H_{\delta_{1}}$.\qed
\enddemo

Let $E_{1},E_{2}$ be divisors on $\Omega_{A_{1}}$ defined by
$E_{1}=\sum_{\delta\sim\delta_{1}}H_{\delta}$ and
$E_{2}=\sum_{\delta\sim\delta_{2}}H_{\delta}$ where $\delta$ runs over
$\Delta(\roman{II}_{2,18}\oplus A_{1}(-1))$ and $\delta\sim\delta'$
if $\delta=\gamma\cdot\delta'$ for some
$\gamma\in\Gamma_{A_{1}}$. Identify $\Delta_{6}^{2}$ to be a 
$\Gamma_{A_{1}}$-invariant section of $\lambda_{A_{1}}^{\otimes15}$ over
$\Omega_{A_{1}}$, namely $\Delta_{6}^{2}\in H^{0}(\Omega_{A_{1}},
\lambda_{A_{1}}^{\otimes15})^{\Gamma_{A_{1}}}$. Since the quotient
map $\Omega_{A_{1}}\to\Cal M_{A_{1}}$ blanches of
order $2$ along $\Cal D_{A_{1}}=E_{1}+E_{2}$, we get 
the following from Lemmas 10.4 and 10.5.

\proclaim{Lemma 10.6}
$\Delta_{6}^{2}$ vanishes of order $2l$ $(l\in\Bbb Z)$ along $E_{1}$ and
of order $4$ along $E_{2}$.
\endproclaim

\proclaim{Proposition 10.2}
$\Psi_{\roman{II}_{2,18}\oplus A_{1}(-1)}(z,\Theta_{E_{7}}/\Delta)$
is a modular form of weight $75$ with zero divisor $57E_{1}+E_{2}$ where
$\Theta_{E_{7}}(\tau)$ is the theta series of the $E_{7}$-lattice.
\endproclaim

\demo{Proof}
Since $\Theta_{E_{7}}(\tau)/\Delta(\tau)$ is a modular form of type 
$\rho_{\roman{II}_{2,18}\oplus A_{1}(-1)}$ of weight $-17/2$, we can
construct Borcherds's product 
$\Psi_{\roman{II}_{2,18}\oplus A_{1}(-1)}(z,\Theta_{E_{7}}/\Delta)$
by Theorem 8.1. Since $\theta_{E_{7}}(\tau)=1+126q+O(q^2)$,
$\theta_{E_{7}+1/2}(\tau)=56q^{3/4}+576q^{7/4}+O(q^2)$ and
$1/\Delta(\tau)=q^{-1}+24+O(q)$ (cf. [C-S]), we get
$$
\Theta_{E_{7}}(\tau)/\Delta(\tau)=(q^{-1}+150+O(q))\,e_{0}+
(56q^{-1/4}+O(q^{3/4}))\,e_{1}.
\tag 10.4
$$
By Theorem 8.1, 
$\Psi_{\roman{II}_{2,18}\oplus A_{1}(-1)}(z,\Theta_{E_{7}}/\Delta)$ 
has weight $c_{0}(0)/2=75$ whose zero divisor is
$\sum_{\delta\in\Delta(\roman{II}_{2,18}\oplus A_{1}(-1))}H_{\delta}
+56\sum_{\delta\in\Delta(\roman{II}_{2,18}\oplus A_{1}(-1)),
\delta/2\in(\roman{II}_{2,18}\oplus A_{1}(-1))^{\lor}}H_{\delta}=
57E_{1}+E_{2}$.\qed
\enddemo

\proclaim{Theorem 10.1}
$\varDelta_{A_{1}}^{15}(z)=C_{A_{1}}
\Delta_{6}^{8}(z)/
\Psi_{\roman{II}_{2,18}\oplus A_{1}(-1)}(z,\Theta_{E_{7}}/\Delta)$.
\endproclaim

\demo{Proof}
Put $F(z):=\Delta_{6}^{8}(z)/[\varDelta_{A_{1}}(z)^{15}\cdot
\Psi_{\roman{II}_{2,18}\oplus A_{1}(-1)}(z,\Theta_{E_{7}}/\Delta)]$.
Since $\varDelta_{A_{1}}$ is an automorphic form of weight $(-5,4)$,
$\Delta_{6}^{2}$ of $(0,15)$ and
$\Psi_{\roman{II}_{2,18}\oplus A_{1}(-1)}(z,\Theta_{E_{7}}/\Delta)$
of $(75,0)$ by Theorem 7.2, Propositions 10.1 and 10.2, $F$ is an
automorphic form of weight $(0,0)$. Namely, taking a higher power of $F$
if necessary, $F$ is a $\Gamma_{A_{1}}$-invariant meromorphic
function on $\Omega_{A_{1}}$, thus a meromorphic function on the
modular variety $\Cal M_{A_{1}}$. From Theorem 7.2, Lemma 10.6 
and Proposition 10.2, it follows that $\div(F)=8(l-9)E_{1}$. Since $F$ 
extends to a meromorphic function on the Satake-Baily-Borel 
compactification of $\Cal M_{A_{1}}$, the residue theorem 
implies $l=9$ and $F$ must be a (nonzero) constant.\qed
\enddemo

\remark{Remark}
It seems that it is in [B-K-P-S-B] that 
$\Psi_{\roman{II}_{2,18}\oplus A_{1}(-1)}(z,\Theta_{E_{7}}/\Delta)$
first appeared. From their construction, its Weyl vector is the 
projection of that of $\II_{1,1}\oplus E_{8}(-1)^{\oplus3}$ to
$\II_{1,1}\oplus E_{8}(-1)^{\oplus2}\oplus A_{1}(-1)$. They used this 
modular form to show that the moduli space of $K3$ surfaces of degree 
$2$ is quasi-affine. See [B-K-P-S-B, Theorem 1.3 and Example 2.1] 
for the details.
\endremark

\subsubhead
\bf{10.2 An Explicit Formula for 
$\pmb{\varDelta_{\roman{II}_{1,1}(2)}}$}\rm
\endsubsubhead
Let $(X,\iota)$ be a $2$-elementary $K3$ surface of type 
$\II_{1,1}(2)$. $(X,\iota)$ is said to be {\it generic }\rm if 
$E\cdot\iota(E)>0$ for any $-2$-curve $E$. (Equivalently, 
$E+\iota(E)$ lies in the positive cone: $(E+\iota(E))^{2}\geq0$.)
Let $e,f$ be a basis of $\II_{1,1}(2)$ such that $e^{2}=f^{2}=0$ and 
$e\cdot f=2$. Let $\phi$ be a marking of $(X,\iota)$. If it is generic,
we may suppose that $C_{e}:=\phi^{-1}(e)$ and $C_{f}:=\phi^{-1}(f)$ 
are nef divisors because $2C_{e}\cdot E=C_{e}\cdot(E+\iota(E))\geq0$ 
for any $-2$-curve $E$. Thus, $X$ has two elliptic fibrations associated 
to the linear systems $|C_{e}|$ and $|C_{f}|$. Since $C_{e}+C_{f}$ is 
ample, the linear system $|C_{e}+C_{f}|$ induces a finite surjective
morphism $\Phi:X\to\Bbb P^{1}\times\Bbb P^{1}(\subset\Bbb P^{3})$,
because $H^{0}(X,\phi^{-1}(e+f))=H^{0}(X,\phi^{-1}(e))\otimes
H^{0}(X,\phi^{-1}(f))$. By the Lefschetz formula, these spaces consist 
of $\iota$-invariant sections. Thus $\Phi$ induces a map;
$\Phi/\iota:X/\iota\to\Bbb P^{1}\times\Bbb P^{1}$. As it induces an
isomorphism of the Picard lattice, it must be an isomorphism. Thus,
any generic $2$-elementary $K3$ surfaces of type $\II_{1,1}(2)$ is
realized as a double cover of $\Bbb P^{1}\times\Bbb P^{1}$ blanching
along a smooth curve of bidegree $(4,4)$. By this description, we can 
construct the discriminant of fixed curves as that of curves of 
bidegree $(4,4)$ on $\Bbb P^{1}\times\Bbb P^{1}$.
\par
Let $H_{4,4}=\Bbb P(\Sym^4\Bbb C^2\otimes\Sym^4\Bbb C^2)^{\lor}$ be the 
set of all curves of bidegree $(4,4)$ on $\Bbb P^{1}\times\Bbb P^{1}$. 
As before, $\xi=(\xi_{ij,kl})_{i+j=4,k+l=4}\in H_{4,4}$ represents the 
curve defined by $C_{\xi}:=\{((x,y),(z,w))\in\Bbb P^{1}\times\Bbb P^{1};\,
\sum\xi_{ij,kl}x^{i}y^{j}z^{k}w^{l}=0\}$. Let
$\pi:\Cal C_{4,4}\to H_{4,4}$ be the universal family i.e.
$\pi^{-1}(\xi)=C_{xi}$, on which acts the group of projective 
transformations $PG:=P(GL(2,\Bbb C)\times GL(2,\Bbb C))$. Then, above 
description implies an isomorphism $i:\Cal M_{\roman{II}_{1,1}(2)}^{00}
\ni[(X,\iota)]\to[X^{\iota}]\in H_{4,4}^{sm}/PG$ where $H_{4,4}^{sm}$ is 
the set of all smooth curves in $H_{4,4}$ and 
$\Cal M_{\roman{II}_{1,1}(2)}^{00}$ is the set of all isomorphism classes 
of {\it generic }\rm $2$-elementary $K3$ surfaces of type $\II_{1,1}(2)$. 
Let $\lambda(\Cal C_{4,4}/H_{4,4})$ be the determinant of cohomology on 
which acts $G:=GL(2,\Bbb C)\times GL(2,\Bbb C)$. Note that the action of 
$G$ is not effective. As before, we get an identification
$\lambda_{\roman{II}_{1,1}(2)}=i^{*}
\lambda(\Cal C_{4,4}/H^{nod}_{4,4})^{G}$ on 
$\Cal M^{00}_{\roman{II}_{1,1}(2)}\cup\Cal D_{\roman{II}_{1,1}(2)}^{0}$
where $H^{nod}_{4,4}$ is the set of all curve of bidegree $(4,4)$ with 
at most {\it one node }\rm. Put $F(x,y,z,w;\xi):=
\sum_{i+j=k+l=4}\xi_{ij,kl}\,x^{i}y^{j}z^{k}w^{l}$. We regard $(x,y)$ and
$(z,w)$ as the homogeneous coordinates of $\Bbb P^{1}$. By the Poincar\'e
residue sequence, we get the following.

\proclaim{Lemma 10.7}
$$
H^{0}(C_{\xi},\Omega^{1}_{C_{\xi}})=\bigoplus_{i+j=k+l=2}\Bbb C\,
\Res_{C_{\xi}}\frac{x^{i}y^{j}z^{k}w^{l}\,(x\,dy-y\,dx)\wedge
(z\,dw-w\,dz)}{F(x,y,z,w;\xi)}.
$$
\endproclaim

As before, let us define a local section of 
$\lambda(\Cal C_{4,4}/H^{nod}_{4,4})$ by
$$
\omega(\xi):=\bigwedge_{i+j=k+l=2}
\Res_{C_{\xi}}\frac{x^{i}y^{j}z^{k}w^{l}\,(x\,dy-y\,dx)\wedge
(z\,dw-w\,dz)}{F(x,y,z,w;\xi)}\in\det H^{0}(C_{\xi},\Omega^{1}).
\tag 10.5
$$

\proclaim{Lemma 10.8}
For any $(g,h)\in G$, $(g,h)^{*}\omega=(\det(g)\det(h))^{18}\,\omega$.
\endproclaim

\demo{Proof}
Since the proof is similar to that of Lemma 10.2, we leave it to the 
reader.\qed
\enddemo

\proclaim{Lemma 10.9}
There exists a homogeneous polynomial $D_{4,4}(\xi)$ of degree $68$
in the $\xi$-variable such that the discriminant locus of 
$\pi:\Cal C_{4,4}\to H_{4,4}$ coincides with $\div(D_{4,4})$. Moreover,
for any $(g,h)\in G$, $D_{4,4}((g,h)\cdot\xi)=(\det(g)\det(h))^{-136}
D_{4,4}(\xi)$.
\endproclaim

\demo{Proof}
Let $Q=\Bbb P^{1}\times\Bbb P^{1}(\subset\Bbb P^{3})$ be a hyperquadric.
Since the line bundle of bidegree $(4,4)$ is $-2K_{Q}$ ($K_{Q}$ is the 
canonical bundle), let us consider the Veronese embedding
$v_{4,4}:Q\hookrightarrow\Bbb P(H^{0}(Q,-2K_{Q}))^{\lor}$. Let 
$X:=v_{4,4}(Q)$ be the image and $X^{\lor}$ be the projective dual variety
of $(X,\Cal O_{\Bbb P(H^{0}(Q,-2K_{Q}))^{\lor}}(1))$. By Katz's formula
([Kt, Cor. 5.6]), $X^{\lor}$ is a projective hypersurface of degree
$$
\deg X^{\lor}=\int_{Q}\frac{c(Q)}{(1+v_{4,4}^{*}
c_{1}(\Bbb P(H^{0}(Q,-2K_{Q}))^{\lor}))^{2}}=
\int_{Q}\frac{1+2H+2H^{2}}{(1+4H)^{2}}=68
\tag 10.6
$$
where $H=\Cal O_{\Bbb P^{3}}(1)$. Let $D_{4,4}(\xi)$ be a defining equation
of $X^{\lor}$. Since $X^{\lor}$ is stable under the action of $G$, 
there exists $l\in\Bbb Z$ such that
$D_{4,4}((g,h)\cdot\xi)=(\det(g)\det(h))^{l}\,D_{4,4}(\xi)$. Putting
$g=\lambda\,I_{2}$, $h=I_{2}$, we find $l=-136$.\qed
\enddemo

\proclaim{Proposition 10.3}
$\Delta_{4,4}^{9}(\xi):=D_{4,4}(\xi)^{9}\cdot\omega(\xi)^{\otimes68}$ 
becomes a $G$-invariant section of $\lambda(\Cal C_{4,4}/H_{4,4}^{nod})$. 
In particular, we may regard $\Delta_{4,4}(\xi)^{9}\in H^{0}(
\Cal M_{\roman{II}_{1,1}(2)}^{00},
\lambda_{\roman{II}_{1,1}(2)}^{\otimes68})$.
\endproclaim

\demo{Proof}
The assertion follows from Lemmas 10.8 and 10.9.\qed
\enddemo

\proclaim{Lemma 10.10}
The discriminant locus $H:=\Cal D_{\roman{II}_{1,1}(2)}/
\Gamma_{\roman{II}_{1,1}(2)}$ is an irreducible divisor of 
$\Cal M_{\roman{II}_{1,1}(2)}$. $\Delta_{4,4}^{9}$ vanishes of order $9$
along $H$.
\endproclaim

\demo{Proof}
The first assertion follows from Nikulin's theory 
([Ni, Proposition 1.15.1]), and the second from Proposition 10.3.\qed
\enddemo

Let us discuss {\it special }\rm $2$-elementary $K3$ surfaces of type
$\II_{1,1}(2)$. A $2$-elementary $K3$ surface of type $\II_{1,1}(2)$, 
$(X,\iota)$, is said to be special if it has a $-2$-curve $E$ such that
$E$ and $\iota(E)$ is disjoint. (This definition is analogous to that of
Enriques surfaces ([B-P-V-V, Theorem 18.2], [Na, Remark 4.6]).)
Let $\phi$ be a marking of $(X,\iota)$ and 
$I:=\phi\circ\iota^{*}\phi^{-1}$ be the involution induced by $\iota$.
Put $\delta:=\phi(E)\in\Delta(L_{K3})$. Then, 
$\langle\delta,I(\delta)\rangle=0$ by assumption. 
Let $T:=(\II_{1,1}(2))^{\perp}=\II_{1,1}(2)\oplus\I_{1,17}$
be the transcendental lattice. Let $p:L_{K3}\to T^{\lor}$ be the 
orthogonal projection. We set $d:=p(\delta)=(\delta-I(\delta))/2\in 
T^{\lor}$. Then, $d^2=-1$ and $d\equiv(e+f)/2\mod T$ where $\{e,f\}$ is 
the basis of $\II_{1,1}(2)$ as before. Conversely, if $d'\in T^{\lor}$ is 
a vector such that $(d')^2=-1$ and $d'\equiv(e+f)/2\mod T$, one can easily 
verify that there exists $\delta'\in\Delta(L_{K3})$ such that 
$\langle\delta',I(\delta')\rangle=0$ and $d'=p(\delta')$. 

\proclaim{Lemma 10.11}
$(X,\iota)$ is a special $2$-elementary $K3$ surface of type $\II_{1,1}(2)$
iff its period lies on the divisor
$\Cal D'_{\roman{II}_{1,1}(2)}:=\sum_{d\in T^{\lor},d^2=-1}H_{d}$.
\endproclaim

\demo{Proof}
Suppose that $(X,\iota)$ is special. Let $\omega$ be a canonical form of 
$(X,\iota)$. Let $E$ be a $-2$-curve as above. Since 
$\langle\omega,E\rangle=\langle\phi(\omega),d\rangle=0$, 
we get $\phi(\omega)\in H_{d}$. Conversely, let
$(X,\iota,\phi)$ be a marked $2$-elementary $K3$ surface of type 
$\II_{1,1}(2)$ whose period lies on $H_{d}\subset\Cal D'$. Let 
$\delta\in\Delta(L_{K3})$ such that $\langle\delta,I(\delta)\rangle=0$, 
$p(\delta)=d$. Let $E=\phi^{-1}(\delta)$ be an effective divisor and 
$E=\sum m_{i}\,C_{i}+\sum n_{j}\,E_{j}$ be the irreducible 
decomposition where $C_{i}^{2}\geq0$, $E_{j}^{2}=-2$ and $m_{i},n_{j}>0$.
By assumption, $\langle E,\iota(E)\rangle=0$. Thus, 
$0\geq\langle\sum n_{j}E_{j},\sum n_{j}\iota(E_{j})\rangle$. 
Suppose that there is no $E_{j}$ such that
$\langle E_{j},\iota(E_{j})\rangle=0$. Since there is no 
$\iota$-invariant $-2$-curve and $\Gamma^{2}\geq-2$ for any irreducible
curve, we get $\langle E_{j},\iota(E_{j})\rangle\geq2$. Note that
$\langle E,\iota(E)\rangle=(E+\iota(E))^{2}/2-2\equiv0\mod2$. Since
$$
\aligned
0\geq\langle\sum n_{j}E_{j},\sum n_{j}\iota(E_{j})\rangle
&\geq\sum_{E_{j}=\iota(E_{i})}\langle n_{i}E_{i}+n_{j}E_{j},
n_{i}\iota(E_{i})+n_{j}\iota(E_{j})\rangle\\
&=\sum_{E_{j}=\iota(E_{i})}\langle n_{i}E_{i}+n_{j}E_{j},
n_{i}E_{j}+n_{j}E_{i}\rangle\\
&\geq\sum_{E_{j}=\iota(E_{i})}2(n_{i}-n_{j})^{2}\geq0,
\endaligned
\tag 10.7
$$
we can rewrite $\sum n_{i}E_{i}=\sum n_{k}(E_{k}+\iota(E_{k}))$. As
$(E_{k}+\iota(E_{k}))^{2}\geq0$, it is impossible that $E^2=-2$. Thus, 
$\langle E_{j},\iota(E_{j}\rangle=0$ for some $E_{j}$ and we get 
the assertion.\qed 
\enddemo

\proclaim{Lemma 10.12}
$\Delta_{4,4}^{9}$ has an algebraic singularity along 
$H':=\Cal D'_{\roman{II}_{1,1}(2)}/\Gamma_{\roman{II}_{1,1}(2)}$.
\endproclaim

\demo{Proof}
Since special $2$-elementary $K3$ surfaces of type $\II_{1,1}(2)$ are
realized as a double cover of singular quadric surface like special 
Enriques surfaces (cf. [B-P-V-V, Theorem 18.2]), we can prove 
the assertion analogously to Lemma 10.5.\qed
\enddemo

Let us regard $\Delta_{4,4}^{9}$ as a 
$\Gamma_{\roman{II}_{1,1}(2)}$-invariant section of 
$\lambda_{\roman{II}_{1,1}(2)}^{\otimes68}$ over 
$\Omega_{\roman{II}_{1,1}(2)}$.

\proclaim{Lemma 10.13}
$\div(\Delta_{4,4}^{9})=18\,\Cal D_{\roman{II}_{1,1}(2)}+
l\,\Cal D'_{\roman{II}_{1,1}(2)}$ for some $l\in\Bbb Z$.
\endproclaim

\demo{Proof}
Since the projection map $\Omega_{\roman{II}_{1,1}(2)}\to
\Cal M_{\roman{II}_{1,1}(2)}$ ramifies of order $2$ along the discriminant
locus, we get the assertion by Lemmas 10.10, 10.11, and 10.12.\qed
\enddemo

Let $e_{00},e_{01},e_{10},e_{11}$ be the basis of $\Bbb C[T^{\lor}/T]=
\Bbb C[\II_{1,1}(2)^{\lor}/\II_{1,1}(2)]$ such that $e_{\alpha\beta}$ 
corresponds to $(\alpha\,e+\beta\,f)/2\in\II_{1,1}(2)^{\lor}$. 
Let $f(\tau):=\eta(\tau)^{-8}\eta(2\tau)^{-8}$ be a modular form for
$\Gamma_{0}(2)$. Modifying Proposition 8.1 for $\Gamma_{0}(2)$, we can 
verify that
$$
\tilde{f}(\tau):=f(\tau)\,e_{00}+8\{f(\tau/2)+
f((\tau+1)/2)\,(e_{00}+e_{01}+e_{10})+8\{f(\tau/2)-f((\tau+1)/2)\,e_{11}
\tag 10.8
$$
is a modular form of type $\rho_{T}$ of weight $-8$.

\proclaim{Proposition 10.4}
$\Psi_{T}(z,f):=\Psi_{T}(z,\tilde{f})$ is a modular form of weight $68$ 
with zero divisor $\Cal D_{\roman{II}_{1,1}(2)}+
16\Cal D'_{\roman{II}_{1,1}(2)}$.
\endproclaim

\demo{Proof}
Since $\eta(\tau)^{-8}\eta(2\tau)^{-8}=q^{-1}+8+O(q)$, we get
$$
\tilde{f}(\tau)=(q^{-1}+132+O(q))\,e_{00}+O(1)\,(e_{01}+e_{10})+
(16q^{-1/2}+O(q^{1/2}))\,e_{11},
\tag 10.9
$$
which, together with Theorem 8.1, yields the assertion.\qed
\enddemo

\proclaim{Theorem 10.2}
$\varDelta_{\roman{II}_{1,1}(2)}^{17}=C_{\roman{II}_{1,1}(2)}\,
\Delta_{4,4}^{9}(z)/\Psi_{\roman{II}_{1,1}(2)\oplus\roman{II}_{1,17}}
(z,\eta(\tau)^{-8}\eta(2\tau)^{-8})$.
\endproclaim

\demo{Proof}
Put $F(z):=\Delta_{4,4}^{9}(z)/[\varDelta_{\roman{II}_{1,1}(2)}(z)^{17}
\cdot\Psi_{\roman{II}_{1,1}(2)\oplus\roman{II}_{1,17}}
(z,\eta(\tau)^{-8}\eta(2\tau)^{-8})]$. Since 
$\varDelta_{\roman{II}_{1,1}(2)}$ is an automorphic form of weight 
$(-4,4)$, $\Delta_{4,4}^{9}$ of $(0,68)$, and $\Psi_{T}(z,f)$ of $(68,0)$,
$F(z)$ is a $\Gamma_{\roman{II}_{1,1}(2)}$-invariant meromorphic function
on $\Omega_{\roman{II}_{1,1}(2)}$. From Theorem 7.2, Lemma 10.3 and 
Proposition 10.4, it follows that 
$\div(F)=(l-16)\Cal D'_{\roman{II}_{1,1}(2)}$. Thus, $F$ descends to 
a meromorphic function on $\Cal M_{\roman{II}_{1,1}(2)}$ with divisor 
$(l-16)H'/2$. By the residue theorem, we get $l=16$ and $F$ must be a 
nonzero constant.\qed
\enddemo

\remark{Remark}
As in the case of Borcherds's $\Phi$-function, we can show
$$
\Psi_{\roman{II}_{1,1}(2)\oplus\roman{II}_{1,17}}
(z,\eta(\tau)^{-8}\eta(2\tau)^{-8})^{2}=\Psi_{\roman{I_{2,18}(2)}}
(z,\eta(\tau)^{-16}\eta(2\tau)^{16}\eta(4\tau)^{16})
$$ 
if we regard
$\II_{1,1}(2)=\{m\,h+n\,\delta\in\I_{1,1};\,m\equiv n\mod 2\}$ and
$\II(2)\oplus\I_{1,17}$ as a sublattice of $\I_{2,18}$. Here, $h^2=1$, 
$\delta^2=-1$, and $h\cdot\delta=0$ is the basis of $\I_{1,1}$. It is 
interesting that $\eta(\tau)^{-16}\eta(2\tau)^{16}\eta(4\tau)^{16}$ is 
the square of the modular form used to get Borcherds's $\Phi$-function.
\endremark

\beginsection
$\S11$. A Theta Product and Borcherds's Product

\par
The Nikulin type of the lattice $S_5$ in (8.5) is 
$(r,l,\delta)=(16,6,1)$ and a 2-elementary $K3$ surface of type $S_{5}$ 
is the minimal resolution of the double cover of $\Bbb P^{2}$ blanching 
along generic $6$ lines. Therefore, the moduli space of 2-elementary 
$K3$ surfaces of type $S_{5}$ is isomorphic to the configuration space 
$\Bbb X(3,6)$ of six points in $\Bbb P^{2}$, and the period map induces 
a morphism $Prd:\Bbb X(3,6)\to\Cal M_{S_{5}}$. In [Ma], Matsumoto 
described the inverse of the period map by using the theta function. 
Let $\Bbb H_{2}$ be the domain defined by
$\Bbb H_{2}:=\{W\in M(2,2;\Bbb C);\,(W-W^{*})/2i>0\}$ where 
$W^{*}={}^{t}\overline{W}$. Identification of $\Bbb H_{2}$ and
$\Omega_{S_{5}}\cong\Lambda_{5}+\sqrt{-1}C_{\Lambda_{5}}$ is given by
$$
\Bbb H_{2}\ni y=\left(\matrix
y_{0}+y_{1}&\frac{y_{0}+y_{1}+y_{2}+iy_{3}}{1+i}\\
\frac{y_{0}+y_{1}+y_{2}-iy_{3}}{1-i}&y_{0}+y_{2}
\endmatrix\right)\to
(1:-\det y:y_{0}:y_{1}:y_{2}:y_{3})\in\Omega_{S_{5}}
\tag 11.1
$$
under which $\Gamma_{M}(1+i)$ (an arithmetic subgroup of 
$\hbox{Aut}(\Bbb H_{2})$) is identified with $\Gamma_{S_{5}}$ by
[Ma, Proposition 1.5.1]. Note that the quadratic form attached to 
$\Lambda_{5}$ is given by $q(y)=4\det(y)$. On $\Bbb H_{2}$ exist theta 
functions, ten of which are represented by
$$
\Theta_{a,b}(W)=
\sum_{m\in\Bbb Z[\sqrt{-1}]^{2}}
\exp\pi i\left\{(m+a)^{*}W(m+a)+2\Re\,b^{*}m\right\}
\tag 11.2
$$
where $a,b\in\{0,(1+i)/2\}^{2}$ with $a^{*}b\in\Bbb Z$.
Any one of these ten theta functions is said to be an even theta 
function in this paper.

\proclaim{Theorem 11.1}
Via the identification of $\Omega_{S_{5}}$ with $\Bbb H_{2}$,
$2^{12}\varDelta_{S_{5}}=\prod_{\roman{even}}\Theta_{a,b}$. 
\endproclaim

\demo{Proof}
By [Ma, 1.4, Proposition 1.5.1 and Lemma 2.3.1], 
$\div(\Theta_{a,b})$ is the $\Gamma_{S_{5}}$-orbit of
$H_{\alpha(a,b)}+H_{\alpha'(a,b)}$ where
$\alpha(a,b),\alpha'(a,b)\in\Delta(T_{5})$ are roots defined in 
[Ma, 1.4]. (By [Ma, 2.3], $(a,b)$ corresponds uniquely to a 
partition of $\{1,\cdots,6\}$ into $\{i,j,k\}\cup\{l,m,n\}$. 
Then, $\alpha(a,b):=\alpha(i,j,k)$ and
$\alpha'(a,b):=\alpha(l,m,n)$ in [Ma, 1.4].) In particular, 
on $\Omega_{S_{5}}$, one has
$$
\div(\prod\Theta_{a,b})=\sum H_{\alpha(a,b)}+H_{\alpha'(a,b)}.
\tag 11.3
$$
By [Ma, Proposition 3.1.1], $(\prod\Theta_{a,b})^{2}$ 
becomes a modular form of weight $20$ relative to $\Gamma_{M}(1+i)$. 
Consider the function 
$\varDelta_{S_{5}}^{2}/(\prod\Theta_{a,b})^{2}$ (or its higher power 
if necessary) which descents to a meromorphic function on 
$\Cal M_{S_{5}}$ by the automorphic property. Compared with
Theorem 7.2, it has no pole and thus is a constant. Namely, there 
exists a constant $C$ such that $C\,\varDelta_{S_{5}}=\prod\Theta_{a,b}$. 
Comparing the first non-zero Fourier coefficient ([G-N1, (1.7)]), 
we get $C=2^{12}$.\qed
\enddemo

Put $h=\left(\matrix 1&\frac{1}{1+i}\\\frac{1}{1-i}&1\endmatrix\right)$,
$\delta_{6}=
\left(\matrix 1&\frac{1}{1+i}\\\frac{1}{1-i}&0\endmatrix\right)$,
$\delta_{7}=
\left(\matrix 0&\frac{1}{1+i}\\\frac{1}{1-i}&1\endmatrix\right)$,
$\delta_{8}=
\left(\matrix 0&\frac{i}{1+i}\\\frac{-i}{1-i}&0\endmatrix\right)$, and
$\delta_{6}'=h-\delta_{7}-\delta_{8}$, 
$\delta_{7}'=h-\delta_{6}-\delta_{8}$,
$\delta_{8}'=h-\delta_{6}-\delta_{7}$. Then, the set of simple roots of
$\Lambda_{5}$ is $\{\delta_{6},\delta_{7},\delta_{8},
\delta_{6}',\delta_{7}',\delta_{8}'\}$, and the Weyl vector is
$2\rho_{5}=\sum_{6\leq i\leq8}\delta_{i}+\sum_{6\leq i\leq8}\delta_{i}'$.
Put $q_{i}=e^{2\pi i y_{i}}$ $(1=0,1,2,3)$,
$y=y_{0}h+y_{1}\delta_{6}+y_{2}\delta_{7}+y_{3}\delta_{8}$ and
$r=r_{0}h+r_{1}\delta_{6}+r_{2}\delta_{7}+r_{3}\delta_{8}\in\Bbb H_{2}$.
Comparing Theorem 8.1 and 11.1, we get the following.

\proclaim{Corollary 11.1}
$$
2^{-12}\prod_{(a,b)\,\roman{even}}\Theta_{a,b}(y)=
q_{0}^{3}q_{1}^{-1}q_{2}^{-1}q_{3}^{-1}\prod_{\epsilon\in\{0,1\}}
\prod_{r\in\Pi^{+}_{\epsilon}(\Lambda_{5})}
(1-q_{0}^{r_{0}}q_{1}^{-r_{1}}q_{2}^{-r_{2}}q_{3}^{-r_{3}}
)^{c_{5,\epsilon}(2\det r)}
$$
where $\Pi^{+}_{\epsilon}(\Lambda_{5})=
(\sum_{6\leq i\leq8}(\Bbb Z+\frac{\epsilon}{2})_{\geq0}\delta_{i}+
\sum_{6\leq i\leq8}(\Bbb Z+\frac{\epsilon}{2})_{\geq0}\delta'_{i})
\backslash\{0\}$, and $\{c_{5,\epsilon}(m)\}$ is the same as in 
Theorem $8.2$.
\endproclaim


\Refs
\widestnumber\no{9999999999}

\ref
\no [A]\by Allcock, D.
\paper The period lattice for Enriques surfaces
\jour math.AG/9905166\yr 1999
\endref

\ref
\no [B-P-V-V]\by Barth, W., Peters, C., Van de Ven, A.
\book Compact Complex Surfaces
\publ Springer
\endref

\ref
\no [Be]\by Beauville, A.
\paper Application aux \'espace de module
\jour Ast\'erisque\vol 126\yr 1985\pages 141-152
\endref

\ref
\no [B-G-V]\by Berline, N., Getzler, E., Vergne, M.
\book Heat Kernels and Dirac Operators
\publ Springer
\endref

\ref
\no[B-C-O-V]\by Bershadsky, M., Cecotti, S., Ooguri, H., Vafa, C.
\paper Kodaira-Spencer theory of gravity and exact results for quantum
string amplitudes
\jour Comm. Math. Phys.\vol 165\yr 1994
\endref

\ref
\no [Bi]\by Bismut, J.-M.
\paper Equivariant immersions and Quillen metrics
\jour J. Diff. Geom.\vol 41\yr 1995\pages 53-157
\endref

\ref
\no [B-B]\by Bismut, J.-M., Bost, J.-B.
\paper Fibr\'es d\'eterminants, m\'etrique de Quillen et 
d\'eg\'en\'erescence des courbes
\jour Acta Math.\vol 165\yr 1990\pages 1-103
\endref

\ref
\no [B-G-S]\by Bismut, J.-M., Gillet, H., Soul\'e, C.
\paper Analytic torsion and holomorphic determinant bundles I, II, III
\jour Comm. Math. Phys.\vol 115\yr 1988\pages 49-78, 79-126, 301-351
\endref

\ref
\no [Bc]\by Borcea, C.
\paper $K3$ surfaces with involution and mirror pairs of Calabi-Yau
manifolds, ``Essays on Mirror Manifolds II''\publ International Press
\yr 1996
\endref

\ref
\no [Bo1]\by Borcherds, R.E.
\paper Generalized Kac-Moody algebras
\jour J. Alg.\vol 115\yr 1988\pages 501-512
\endref

\ref
\no [Bo2]\by{---------}
\paper Monstrous moonshine and monstrous Lie superalgebras
\jour Invent. Math.\vol 109\yr 1992\pages 405-444
\endref

\ref
\no [Bo3]\by{---------}
\paper Automorphic forms on $O_{s+2,2}(\Bbb R)$ and infinite products
\jour Invent. Math.\vol 120\yr 1995\pages 161-213
\endref

\ref
\no [Bo4]\by{---------}
\paper The moduli space of Enriques surfaces and the fake monster Lie 
superalgebra
\jour Topology\vol 35\yr 1996\pages 699-710
\endref

\ref
\no [Bo5]\by{---------}
\paper Automorphic forms with singularities on Grassmanians
\jour Invent. Math.\vol 132\yr 1998\pages 491-562
\endref

\ref
\no [Bo6]\by{---------}
\paper An automorphic form related to cubic surfaces
\jour preprint
\endref

\ref
\no [Bo7]\by{---------}
\paper private communications (e-mails to K.-I. Yoshikawa)
\endref

\ref
\no [B-K-P-S-B]\by Borcherds, R., Katzarkov, L., Pantev, T., 
Shepherd-Barron, N.
\paper Families of $K3$ surfaces
\jour J. Alg. Geom.\vol 7\yr 1998\pages 183-193
\endref

\ref
\no [B-C]\by Bott, R., Chern, S.S.
\paper Hermitian vector bundles and the equidistribution of the zeros
of their holomorphic sections
\jour Acta Math.\vol 114\yr 1968\pages 71-112
\endref

\ref
\no [B-R]\by Burns, D., Rapoport, M.
\paper On the Torelli problems for K\"ahlerian $K3$ surfaces
\jour Ann. Sci. Ec. Norm. Sup.\vol IV Ser.8\yr 1975\pages 235-274
\endref

\ref
\no [C-L]\by Cheng, S.-Y., Li, P.
\paper Heat kernel estimates and lower bound of eigenvalues
\jour Comment. Math. Helv.\vol 56\yr 1981\pages 327-338
\endref

\ref
\no [C-L-Y]\by Cheng, S.-Y., Li, P., Yau, S.-T.
\paper On the upper estimate of the heat kernel of a complete Riemannian
manifold
\jour Amer. J. Math.\vol 103\yr 1981\pages 1021-1063
\endref

\ref
\no [C-S]\by Conway, J.H., Sloane, N.J.A.
\book Sphere Packings, Lattices and Groups
\publ Springer
\endref


\ref
\no [F]\by Freitag, E.
\book Siegelshe Modulfunktionen
\publ Springer
\endref

\ref
\no [G-T]\by Gilberg, D., Trudinger, N.S.
\book Elliptic Partial Differential Equations of Second Order
\publ Springer
\endref

\ref
\no [G-N1]\by Gritsenko, V., Nikulin, V.
\paper Siegel automorphic form corrections of some Lorentzian Kac-Moody 
Lie algebras
\jour Amer. J. Math.\vol 119\yr 1997\pages 181-224
\endref

\ref
\no [G-N2]\by{---------}
\paper $K3$ surfaces, Lorentzian Kac-Moody algebras and mirror symmetry
\jour Math. Res. Lett.\vol 3\yr 1996\pages 211-229
\endref

\ref
\no [G-N3]\by{---------}
\paper Automorphic forms and Lorentzian Kac-Moody algebras, I, II
\jour Intern. J. Math.\vol 9\yr 1998\pages 153-199, 201-275
\endref

\ref
\no [H-M]\by Harvey, J., Moore, G.
\paper Exact gravitational threshold correction in the FHSV model
\jour Phys. Rev. D\vol 57\yr 1998\pages 2329-2336
\endref

\ref
\no [J-K]\by Jorgenson, J., Kramer, J.
\paper Star products of Green currents and automorphic forms
\jour preprint\yr 1998
\endref

\ref
\no [J-T1]\by Jorgenson, J., Todorov, A.
\paper A conjectured analogue of Dedekind's eta function for $K3$ surfaces
\jour Math. Res. Lett.\vol 2\yr 1995\pages 359-376
\endref

\ref
\no [J-T2]\by{---------}
\paper Enriques surfaces, analytic discriminants, and Borcherds's 
$\Phi$-function
\jour Comm. Math. Phys.\vol 191\yr 1998\pages 249-264
\endref

\ref
\no [J-T3]\by{---------}
\paper A correction for the paper ``Enriques surfaces ...''
\jour preprint\yr 1998
\endref

\ref
\no [K-R]\by K\"ohler, K., Roessler, D.
\paper A fixed point formula of Lefschetz type in Alakelov geometry 
I,II,III
\jour preprint\yr 1998, 1999
\endref

\ref
\no [Kt]\by Katz, N.M.
\paper Pinceaux de Lefschetz: th\'eor\`eme d'existence, SGA $7$,II
\jour Lect. Notes Math.\vol340\yr 1973\pages 213-253
\endref

\ref
\no [Kz]\by Kazhdan, D.
\paper Connection of the dual space of a group and the structure of its 
closed subgroups
\jour Funct. Anal. Appl.\vol 1\yr 1967\pages 63-65
\endref

\ref
\no [Ko]\by Kobayashi, R.
\paper Moduli of Einstein metrics on a $K3$ surface and degeneration of 
type I
\jour Adv. Study Pure Math.\vol 18-II\yr 1990\pages 257-311
\endref

\ref
\no [K-T]\by Kobayashi, R., Todorov, A.
\paper Polarized period map for generalized $K3$ surfaces and the moduli of
Einstein metrics
\jour Tohoku Math. J.\vol 39\yr 1987\pages 341-363
\endref

\ref
\no [Kn]\by Kond\B o, S.
\paper The rationality of the moduli space of Enriques surfaces
\jour Compositio Math.\vol 91\yr 1994\pages 159-173
\endref

\ref
\no [Kr]\by Kronheimer, P.B.
\paper The construction of ALE spaces as hyper-K\"ahler quotients
\jour J. Diff. Geom.\vol 29\yr 1989\pages 665-683
\endref

\ref
\no [L-Y]\by Li, P., Yau, S.-T.
\paper On the parabolic kernel of the Schr\"odinger operator
\jour Acta Math.\vol 156\yr 1986\pages 153-201
\endref

\ref
\no [Ma]\by Matsumoto, K.
\paper Theta function on the bounded symmetric domain of type $\I_{2,2}$ 
and the period map of a $4$-parameter family of $K3$ surfaces
\jour Math. Ann.\vol 295\yr 1993\pages 383-409
\endref

\ref
\no [Man]\by Manin, Y.
\book Cubic Forms
\publ North Holland
\endref

\ref
\no [Mo]\by Morrison, D.
\paper Some remarks on the moduli of $K3$ surfaces
\jour Progress in Math.\vol 39\yr 1983\pages 303-332
\endref

\ref
\no [Na]\by Namikawa, Y.
\paper Periods of Enriques surfaces
\jour Math. Ann.\vol 270\yr 1985\pages 201-222
\endref

\ref
\no [Ni1]\by Nikulin, V.V.
\paper Integral symmetric bilinear forms and some of their applications
\jour Math. USSR Izvestija\vol 14\yr 1980\pages 103-167
\endref

\ref
\no [Ni2]\by Nikulin, V.V.
\paper Finite automorphism groups of K\"ahler $K3$ surfaces
\jour Trans. Moscow Math. Soc.\vol 38\yr 1980\pages 71-135
\endref

\ref
\no [Ni3]\by{---------}
\paper Involution of integral quadratic forms and their applications to 
real algebraic geometry
\jour Math. USSR Izvestiya\vol 22\yr 1984\pages 99-172
\endref

\ref
\no [Ni4]\by{---------}
\paper On the quotient groups of the automorphism group of hyperbolic 
forms by the subgroup generated by 2-reflections
\jour J. Soviet Math.\vol 22\yr 1983\pages 1401-1476
\endref

\ref
\no [P-S-S]\by Piatetskii-Shapiro, I., Shafarevich, I.R.
\paper A Torelli theorem for algebraic surfaces of type $K3$
\jour Math. USSR Izv.\vol 35\yr 1971\pages 530-572
\endref

\ref
\no [R-S]\by Ray, D.B., Singer, I.M.
\paper Analytic torsion for complex manifolds
\jour Ann. of Math.\vol 98\yr 1973\pages 154-177
\endref

\ref
\no [Sc]\by Schumacher, G.
\paper On the geometry of moduli spaces
\jour Manuscripta Math.\vol50\yr 1985\pages 229-267
\endref

\ref
\no [Sha]\by Shah, J.
\paper A complete moduli space for $K3$ surface of degree $2$
\jour Ann. of Math.\vol 112\yr 1980\pages 485-510
\endref

\ref
\no [Si]\by Siu, U.-T.
\book Lectures on Hermitian-Einstein Metrics for Stable Bundles and 
K\"ahler-Einstein Metrics\publ Birkh\"auser
\endref

\ref
\no [Ti]\by Tian, G.
\paper Smoothness of the universal deformation space of compact Calabi-Yau
manifolds and its Petersson-Weil metric, ``Mathematical Aspects of String 
Theory, $1987$''
\jour Adv. Studies in Math. Physics\vol 1\yr 1987
\endref

\ref
\no [To1]\by Todorov, A.
\paper Applications of the K\"ahler-Einstein-Calabi-Yau metric to moduli of
$K3$ surfaces
\jour Invent. Math.\vol 81\yr 1980\pages 251-266
\endref

\ref
\no [To2]\by{---------}
\paper The Weil-Pertersson geometry of the moduli space of $SU(n\geq3)$
(Calabi-Yau) manifolds I
\jour Comm. Math. Phys.\vol 126\yr 1989\pages 325-346
\endref

\ref
\no [V]\by Voisin, C.
\paper Miroirs et involutions sur les surfaces $K3$
\jour Asterisque\vol 218\yr 1993\pages 273-322
\endref

\ref
\no [Ya]\by Yau, S.-T.
\paper On the Ricci curvature of a compact K\"ahler manifold and 
the complex Monge-Amp\`ere Equation, I
\jour Comm. Pure Appl. Math.\vol 31\yr 1978\pages 339-411
\endref


\ref
\no [Yo1]\by Yoshikawa, K.-I.
\paper Degeneration of algebraic manifolds and the spectrum of Laplacian
\jour Nagoya Math. J.\vol 146\yr 1997\pages 83-129
\endref

\ref
\no [Yo2]\by{---------}
\paper Discriminant of theta divisors and Quillen metrics
\jour math.AG/9904092 (to appear in J. Diff. Geom.)\yr 1999
\endref

\endRefs
\enddocument